# THE SIZES OF THE PIONEERING, LOWEST CROSSING AND PIVOTAL SITES IN CRITICAL PERCOLATION ON THE TRIANGULAR LATTICE

By G. J. Morrow and Y. Zhang

*University of Colorado*

Let $L_n$ denote the lowest crossing of a square $2n \times 2n$ box for critical site percolation on the triangular lattice imbedded in $\mathbf{Z}^2$. Denote also by $F_n$ the pioneering sites extending below this crossing, and $Q_n$ the pivotal sites on this crossing. Combining the recent results of Smirnov and Werner [*Math. Res. Lett.* **8** (2001) 729–744] on asymptotic probabilities of multiple arm paths in both the plane and half-plane, Kesten's [*Comm. Math. Phys.* **109** (1987) 109–156] method for showing that certain restricted multiple arm paths are probabilistically equivalent to unrestricted ones, and our own second and higher moment upper bounds, we obtain the following results. For each positive integer $\tau$, as $n \to \infty$:
1. $E(|L_n|^\tau) = n^{4\tau/3+o(1)}$.
2. $E(|F_n|^\tau) = n^{7\tau/4+o(1)}$.
3. $E(|Q_n|^\tau) = n^{3\tau/4+o(1)}$.
These results extend to higher moments a discrete analogue of the recent results of Lawler, Schramm and Werner [*Math. Res. Lett.* **8** (2001) 401–411] that the frontier, pioneering points and cut points of planar Brownian motion have Hausdorff dimensions, respectively, 4/3, 7/4 and 3/4.

## 1. Introduction.

Consider site percolation on the triangular lattice. Each vertex of the lattice is open with probability $p$ and closed with probability $1-p$ and the sites are occupied independently of each other. We will realize the triangular lattice with vertex set $\mathbf{Z}^2$. For a given $(x,y) \in \mathbf{Z}^2$, its nearest neighbors are defined as $(x \pm 1, y)$, $(x, y \pm 1)$, $(x+1, y-1)$ and $(x-1, y+1)$. Bonds between neighboring or adjacent sites therefore correspond to vertical or horizontal displacements of one unit, or diagonal displacements between two nearest vertices along a line making an angle of 135° with the positive $x$-axis.

Recall that the triangular lattice may also be viewed with sites as hexagons in a regular hexagonal tiling of the plane. This point of view is convenient to describe the fact that when $p = 1/2$ (critical percolation) and the hexagonal mesh tends to zero, the percolation cluster interface has a conformally invariant scaling limit, namely, the stochastic Loewner evolution process $\mathrm{SLE}_6$ [11]. Smirnov and Werner [12] combine the convergence of the cluster interface with recent results on the probabilities of crossings of annular and semi-annular regions by $\mathrm{SLE}_6$ calculated by Lawler, Schramm and Werner [6, 7, 9] to obtain corresponding probabilities for the critical site percolation on the triangular lattice.

We will use the Smirnov and Werner [12] estimates in the case of three-arm, two-arm and four-arm paths to establish results, respectively, on the length of the lowest crossing, the size of the pioneering sites extending below this crossing and the number of pivotal sites on this crossing of a square box with sides parallel to the coordinate axes in $\mathbf{Z}^2$. Here and throughout the paper we will be working with the critical percolation model. To illustrate how our work fits in with known results for planar Brownian motion, we describe various geometric features of the Brownian paths as follows. Define the hull $\mathbf{K}_t$ at time $t$ of a planar Brownian motion $\beta_s$, $s \geq 0$, as the union of the Brownian path $\beta[0,t] := \{\beta_s, 0 \leq s \leq t\}$ with the bounded components of its complement $\mathbf{R}^2 \backslash \beta[0,t]$. The frontier or outer boundary of $\beta[0,t]$ is defined as the boundary of $\mathbf{K}_t$. By contrast, a pioneer point of the Brownian path is defined as any point $\beta_s$ at some time $s \leq t$ such that $\beta_s$ is in the boundary of $\mathbf{K}_s$, that is, such that $\beta_s$ is on the frontier of $\beta[0,s]$. A point $\beta_s$ for some $0 < s < t$ is called a cut point of $\beta[0,t]$ if $\beta[0,s] \cap \beta(s,t] = \varnothing$. Lawler, Schramm and Werner [8] have shown that the frontier, pioneer points and cut points of a planar Brownian motion have Hausdorff dimensions, respectively, 4/3, 7/4 and 3/4. We answer an open question ([12], question 2) to find an analogue of this result in the case of critical percolation on the triangular lattice. Indeed, we asymptotically evaluate all moments of the sizes of the corresponding lowest crossing, pioneering sites and pivotal sites that we define below.

It turns out that our method requires a more careful analysis in the four-arm case than in the two and three-arm cases. As pointed out in [12], the probability estimates of annular crossings of multiple-arm paths [see (2.3) below] lead naturally to a prediction of our first moment results. Only in the pivotal (four-arm) case do we need to apply the estimate of Smirnov and Werner [12] for probabilities of multiple-arm crossings of semi-annular regions, in addition to the basic annular estimates to actually establish the prediction. In all cases, however, the methods of Kesten [4] are essential to construct the probability estimates for our moment calculations. This calculation handles, in particular, the probability of four-arm paths near the boundary of the box used to define the pivotal sites.

Denote by $\mathbf{T}$ the full triangular lattice graph whose vertex set is $\mathbf{Z}^2$ and whose edges are the nearest neighbor bonds. Define $\|\mathbf{x}\| := \max\{|x|, |y|\}$ for $\mathbf{x} = (x, y) \in \mathbf{Z}^2$. For any real number $r \geq 0$, we denote the square box of vertices $B(r) := \{\mathbf{x} \in \mathbf{Z}^2 : \|\mathbf{x}\| \leq r\}$. A path is a sequence of distinct vertices connected by nearest neighbor bonds. Thus, a path is simple. Following Grimmett [2], the boundary or surface of a set $X$ of vertices is the set $\partial X$ of vertices in $X$ that are adjacent to some vertex not in $X$. A path is open (closed) from a set $X$ to a set $Y$ if each vertex of the path is open (closed) and contained in $\mathbf{Z}^2 \setminus (X \cup Y)$ except for the endpoints in $\partial X$ and $\partial Y$ which may or may not be open (closed). The interior of $X$ is defined by $\text{int}(X) = X \setminus \partial X$. A set $X$ of vertices is connected if the graph induced by $X$ is connected as a subgraph of $\mathbf{T}$. Let $R$ be a connected set of vertices lying within a finite union of rectangles with sides parallel to the coordinate axes. We say that a path is "in $R$" if its vertices remain in $R$ except possibly for its endpoints. If, in addition, $R$ is a single such rectangle, a horizontal open (closed) crossing of $R$ is an open (closed) path in $R$ from the left side of $R$ to the right side of $R$. A vertical crossing is defined similarly.

Let $n$ be a positive integer. For each $\mathbf{x} \in B(n)$, we define the event

$$\begin{aligned} \mathcal{L}(\mathbf{x}, n) := {} & \text{there exists a horizontal open crossing of } B(n) \\ & \text{containing the vertex } \mathbf{x}, \text{ and there exists a} \\ & \text{closed path in } B(n) \text{ from } \mathbf{x} \text{ to the bottom of } B(n). \end{aligned} \tag{1.1}$$

The lowest crossing for any given configuration of vertices for which a horizontal open crossing of $B(n)$ exists is known (see [2]) to be the unique horizontal open crossing $\gamma_n$ of $B(n)$ that lies in the region on or beneath any other horizontal open crossing. In fact, on any given configuration we may also represent the set of vertices in $\gamma_n$ as equal to the set $L_n := \{\mathbf{x} : \mathcal{L}(\mathbf{x}, n) \text{ occurs}\}$. Although this fact is well known, we briefly review its proof. First, any vertex $\mathbf{x}$ of $\gamma_n$ admits a closed path to the bottom of $B(n)$ [so that $\mathcal{L}(\mathbf{x}, n)$ occurs], else one could construct a crossing strictly lower than $\gamma_n$. Therefore, $\gamma_n \subset L_n$. On the other hand, to show $L_n \subset \gamma_n$, assume the event $\mathcal{L}(\mathbf{x}, n)$. The open path in this event lies above $\gamma_n$, so the closed path in this event has to cross $\gamma_n$, unless $\mathbf{x} \in \gamma_n$. Thus, the set of vertices of $\gamma_n$ is precisely $L_n$.

We introduce next the pioneering sites extending below the lowest crossing of a configuration in $B(n)$. Define the event

$$\begin{aligned} \mathcal{F}(\mathbf{x}, n) := {} & \mathbf{x} \text{ is open and there exist two open paths in } B(n) \\ & \text{started from } \mathbf{x}, \text{ one to the left side and one to the} \\ & \text{right side of } B(n), \text{ and there exists a closed path} \\ & \text{in } B(n) \text{ from } \mathbf{x} \text{ to the bottom of } B(n). \end{aligned} \tag{1.2}$$

Note that (1.1) implies (1.2). The difference is that in (1.2) the two open paths need not be disjoint, whereas in (1.1) the horizontal crossing through $\mathbf{x}$ breaks up into two disjoint open paths. We define the set of pioneering sites as the set $F_n := \{\mathbf{x} : \mathcal{F}(\mathbf{x}, n) \text{ occurs}\}$. Geometrically, $F_n$ consists of the union of the lowest crossing with the many complicated orbs and tendrils hanging from it; the vertices of these latter sets do not admit two disjoint paths to the sides of $B(n)$. Alternatively, $F_n$ is the set of open sites discovered through the exploration process that starts at the lower left corner of $B(n)$ and runs until it meets the right side, that determines the interface between the lowest spanning open cluster in $B(n)$ and the closed cluster attaching to its bottom side. This description of $F_n$ explains its correspondence to the trace of $\mathrm{SLE}_6$.

Finally, we define the pivotal sites lying along the lowest crossing. Define the event

$$\begin{aligned}\mathcal{Q}(\mathbf{x}, n) := {} & \text{there exists a horizontal open crossing of } B(n) \\ & \text{containing the vertex } \mathbf{x}, \text{ and there exist two} \\ & \text{disjoint closed paths in } B(n) \text{ started from } \mathbf{x}, \text{ one} \\ & \text{to the top side and one to the bottom side of } B(n).\end{aligned} \tag{1.3}$$

We define the set of pivotal sites as the set $Q_n := \{\mathbf{x} : \mathcal{Q}(\mathbf{x}, n) \text{ occurs}\}$. The two closed "arms" emanating from a pivotal (and therefore open) site $\mathbf{x}$ force any horizontal open crossing of $B(n)$ to pass through $\mathbf{x}$. Let $\mathcal{C}_n$ be the open cluster containing the lowest crossing whenever such a horizontal open crossing exists. That a pivotal site exists implies that this cluster also contains the highest horizontal open crossing and that the site belongs to both the highest and lowest crossing.

We can now state our main result. Here and throughout $P$ and $E$ denote, respectively, the probability and expectation for the critical percolation.

THEOREM 1. *For each positive integer $\tau$, as $n \to \infty$:*

1. $E(|L_n|^\tau) = n^{4\tau/3+o(1)}$.
2. $E(|F_n|^\tau) = n^{7\tau/4+o(1)}$.
3. $E(|Q_n|^\tau) = n^{3\tau/4+o(1)}$.

Note that probability upper bounds follow immediately by Markov's inequality from the $\tau$th moment upper bounds in Theorem 1. On the other hand, a bound on the distribution of small values of $|L_n|$ is obtained in [5]. Let $\mathcal{L}$ denote the event that there is a horizontal open crossing of $B(n)$. These authors show that there exist constants $\alpha, c > 0$ and $C_1 < \infty$ such that $P(|L_n| \le n^{1+c} | \mathcal{L}) \le C_1 n^{-\alpha}$ ([5], Theorem 2). We conjecture that this result continues to hold for the triangular lattice case if the exponent $1 + c$ is replaced by $4/3 - \delta$ for any $\delta > 0$, where now $\alpha > 0$ may depend on $\delta$.

A one arm version of Theorem 1 is obtained by Kesten [3]. He shows that there exists a limiting measure $\mu$ on configuration space, conditioned by the event that the origin is connected to $\partial B(n)$ as $n \to \infty$ such that, with respect to $\mu$, there is a unique open cluster $W$ connected to the origin with probability 1 ([3], Theorem 3). He then proves $E_\mu(|W \cap B(n)|^\tau) \asymp (n^2\pi_n)^\tau$, where $\pi_n := P[\mathbf{0}$ is connected by an open path to $(n,\infty) \times \mathbf{R}]$ ([3], Theorem 8). We do not study here the number of sites in $B(n)$ from which a five-arm path to the $\partial B(n)$ exists. By the results of Smirnov and Werner [12] [see (2.3)], the expected number of such sites is predicted to be $n^{o(1)}$.

Results analogous to the above-mentioned Hausdorff dimension properties of certain planar Brownian motion point sets but now for the stochastic Loewner evolution process $\mathrm{SLE}_6$ itself have been obtained as follows: the dimension of the $\mathrm{SLE}_6$ curve (or trace) is 7/4 [1]; the dimension of the (outer) frontier of $\mathrm{SLE}_6$ is 4/3 [7], and the dimension of the set of cut points of $\mathrm{SLE}_6$ is 3/4 [6]; see Remark 5 of [12]. Perhaps for both the Brownian and SLE processes one can obtain moment estimates on the number of disks of radius $\varepsilon > 0$ needed to cover a given one of the above point sets similar to the moment estimates presented here. Some results for the expected number of such disks have been obtained by Rhode and Schramm [10] concerning both the $\mathrm{SLE}_{\boldsymbol{\kappa}}$ curve and hull with $\boldsymbol{\kappa}$ in a range of values including the case $\boldsymbol{\kappa} = 6$ that corresponds to the critical percolation of this paper.

Finally we mention that items 1 and 2 of Theorem 1 may be proved by the same method that Kesten [3] uses to establish the one-arm case mentioned above. However, that method does not extend to the four-arm case since then the exponent in expression (2.3) becomes less than $-1$; see Section 3.1 for further details. We emphasize that the difficulty in this paper lies in the case of higher moments ($\tau \geq 2$) for pivotal sites wherein we study the organization of $\tau$ vertices in the box $B(n)$ at which four-arm events occur. We develop a disjoint boxes method (Section 4) that yields a proof of items 1 and 2 and that also lays a groundwork for the proof of item 3 of Theorem 1. Our organization of the disjoint boxes leads to two developments. First, it allows for the construction of certain horseshoe estimates governed by Lemma 5 that are critical in establishing the correct asymptotic order for even the first moment in the pivotal case. We carry out these constructions in Sections 5 and 6. Second, it allows for the analysis of groups of vertices that are closely clustered together in the analysis of second and higher moments by a separate method based on Lemmas 7 and 8 shown in Section 7. These lemmas are extensions of Kesten's [4] Lemmas 4 and 5. This latter work indeed forms the technological foundations for much of the current paper.

**2. Lower bounds.** In this section we establish lower bounds for each of the moment estimates of Theorem 1. To do this, we begin by writing down the known asymptotic probabilities of multiple-arm paths from [12]. Next Kesten's method is applied to obtain lower bounds for the probabilities of certain restricted multiple-arm path events. The expectation lower bounds then follow easily.

Note that $\mathcal{L}(\mathbf{x},n)$ is a certain sub-event of a so-called three-arm path that we now introduce. Define $B(\mathbf{x},m) := \mathbf{x} + B(m)$. Assume that $B(\mathbf{x},m)$ belongs to the interior of $B(n)$. Denote

$$A(\mathbf{x},m;n) := B(n) \setminus \operatorname{int}(B(\mathbf{x},m)). \tag{2.1}$$

The event of a three-arm path from $B(\mathbf{x},m)$ to $\partial B(n)$ is defined by

$$\begin{aligned} \mathcal{U}_3(\mathbf{x},m;n) := {} & \text{there exist two disjoint open paths in} \\ & A(\mathbf{x},m;n) \text{ from } \partial B(\mathbf{x},m) \text{ to } \partial B(n), \\ & \text{and there exists a closed path in} \\ & A(\mathbf{x},m;n) \text{ from } \partial B(\mathbf{x},m) \text{ to } \partial B(n). \end{aligned} \tag{2.2}$$

We denote $\mathcal{U}_3(\mathbf{x},n) := \mathcal{U}_3(\mathbf{x},0;n)$. We shall use the following estimate from [12].

Lemma 1. *For each fixed $m \geq 0$, $P(\mathcal{U}_3(\mathbf{0},m;n)) = n^{-2/3+o(1)}$ as $n \to \infty$.*

Proof. The proof follows by a direct translation of Theorem 4 of [12] as follows. Consider the event that there exist $\kappa$ disjoint crossings of the annulus $A(r_0,r) := \{z \in \mathbf{C} : r_0 < |z| < r\}$, not all closed nor all open, for the hexagonal tiling of fixed mesh 1 in $\mathbf{C}$. Let $\mathcal{H}_\kappa(r_0,r)$ denote generically any of the sub-events defined by a specific ordering of closed and open crossings among the $\kappa$ disjoint crossings. Then for each $\kappa \geq 2$,

$$P(\mathcal{H}_\kappa(r_0,r)) = r^{-(\kappa^2-1)/12+o(1)} \qquad \text{as } r \to \infty. \tag{2.3}$$

Choose now two open and one closed crossings in the definition of $\mathcal{H}_3(r_0,r)$. Then Lemma 1 follows by applying (2.3) for $\kappa = 3$ and noting, on account of the mild change in geometry between the hexagonal and present models for the triangular lattice, that the event $\mathcal{U}_3(\mathbf{0},m;n)$ satisfies $\mathcal{H}_3(m/2,n) \subset U_3(\mathbf{0},m;n) \subset \mathcal{H}_3(2m,n/2)$. □

Similar to (2.2), we define the events of two-arm and four-arm paths from $B(\mathbf{x},m)$ to $\partial B(n)$ by

$$\begin{aligned} \mathcal{U}_2(\mathbf{x},m;n) := {} & \text{there exist two paths in } A(\mathbf{x},m;n) \\ & \text{from } \partial B(\mathbf{x},m) \text{ to } \partial B(n), \\ & \text{one of them being open and the other closed} \end{aligned} \tag{2.4}$$

and

$$\begin{aligned} \mathcal{U}_4(\mathbf{x},m;n) := {} & \text{there exist two disjoint open paths in} \\ & A(\mathbf{x},m;n) \text{ from } \partial B(\mathbf{x},m) \text{ to } \partial B(n), \\ & \text{and there exist two disjoint closed paths in} \\ & A(\mathbf{x},m;n) \text{ from } \partial B(\mathbf{x},m) \text{ to } \partial B(n). \end{aligned} \tag{2.5}$$

We also denote $\mathcal{U}_\kappa(\mathbf{x},n) := \mathcal{U}_\kappa(\mathbf{x},0;n)$, $\kappa = 2,4$. As in the proof of Lemma 1, except now with $\kappa = 2$ and $\kappa = 4$, respectively in (2.3), we obtain the following.

LEMMA 2. *For fixed* $m \geq 0$, $P(\mathcal{U}_2(\mathbf{0},m;n)) = n^{-1/4+o(1)}$ *and* $P(\mathcal{U}_4(\mathbf{0},m;n)) = n^{-5/4+o(1)}$ *as* $n \to \infty$.

We now make precise the notion of restricted multiple-arm paths. Let $\pi$ be the probability of a given $\kappa$-arm path from a given vertex inside $B(n/4)$ to $\partial B(n)$. We restrict the $\kappa$-arm path by specifying disjoint intervals of length proportional to $n$ and separated by intervals also proportional to $n$ on the $\partial B(n)$ for the hitting sets of the various arms. Kesten [4] shows that there is only a multiplicative constant cost in the probability $\pi$ of this restricted event. In fact, Kesten shows a little more that we will describe explicitly for the Lemmas 3 and 4. To begin, define certain rectangles that sit on the four sides of $B(n)$, counting counterclockwise from the left side of $B(n)$, by $R_1 := [-n,-n/2]\times[-n/2,n/2]$, $R_2 := [-n/2,n/2]\times[-n,-n/2]$, $R_3 := [n/2,n]\times[-n/2,n/2]$ and $R_4 := [-n/2,n/2]\times[n/2,n]$. Let $R$ be a rectangle with sides parallel to the coordinate axes and sharing one side with the boundary of a box $B$. We say that a path $h$-tunnels through $R$ on its way to $\partial B$ if the intersection of the path with the smallest infinite vertical strip containing $R$ remains in $R$. Thus, the path may weave in and out of $R$ but not through the top or bottom sides of $R$, and comes finally to $\partial B$. Likewise, we say that a path $v$-tunnels through $R$ on its way to $\partial B$ if the roles of horizontal and vertical are interchanged in the preceding definition. This definition is consistent with the requirements of Kesten's [4] Lemma 4. Accordingly, for each $\mathbf{x} \in B(n/4)$, we define a certain restricted three-arm path event by

$$\begin{aligned} \mathcal{T}_3(\mathbf{x},n) := \; & \exists \text{ two disjoint open paths in } B(n) \text{ started from } \mathbf{x}, \\ & \text{one to the left side of } B(n) \text{ that } h\text{-tunnels through} \\ & R_1, \text{ and one to the right side of } B(n) \text{ that } h\text{-tunnels} \\ & \text{through } R_3, \text{ and there is a closed path in } B(n) \text{ from } \mathbf{x} \\ & \text{to the bottom of } B(n) \text{ that } v\text{-tunnels through} R_2. \\ & \text{Further, there are vertical open crossings of each of } R_1 \\ & \text{and } R_3, \text{ and there is a horizontal closed crossing of } R_2. \end{aligned} \tag{2.6}$$

By the proof of Kesten's [4] Lemma 4, we obtain the following.

LEMMA 3. *There exists a constant* $C_3$ *such that uniformly for all* $\mathbf{x} \in B(n/4)$, $P(\mathcal{U}_3(\mathbf{x},n)) \leq C_3 P(\mathcal{T}_3(\mathbf{x},n))$.

Note by Lemma 3 that, for $\mathbf{x} \in B(n/4)$, the probabilities of $\mathcal{L}(\mathbf{x},n)$, $\mathcal{U}_3(\mathbf{x},n)$ and $\mathcal{T}_3(\mathbf{x},n)$ are all comparable.

We next define certain restricted two-arm and four-arm path events. In the two-arm case we introduce rectangles that cut off the top and bottom sides of $B(n)$ by $S_2 := [-n,n] \times [-n,-n/2]$ and $S_4 := [-n,n] \times [n/2,n]$. In the four-arm case we have the similar picture as the three-arm case, except now there is one more closed arm that $v$-tunnels through $R_4$ on the way to the top of $B(n)$. We thus define for any $\mathbf{x} \in B(n/4)$,

$$\begin{aligned} \mathcal{T}_2(\mathbf{x},n) := & \text{ there exists an open path in } B(n) \text{ from } \mathbf{x} \\ & \text{ to the top of } B(n) \text{ that } v\text{-tunnels through } S_4, \\ & \text{ and a closed path in } B(n) \text{ from } \mathbf{x} \text{ to the} \\ & \text{ bottom of } B(n) \text{ that } v\text{-tunnels through } S_2. \\ & \text{ Further, there exists a horizontal open crossing} \\ & \text{ of } S_4 \text{ and a horizontal closed crossing of } S_2 \end{aligned} \tag{2.7}$$

and

$$\begin{aligned} \mathcal{T}_4(\mathbf{x},n) := & \exists \text{ two disjoint open paths in } B(n) \text{ started from } \mathbf{x}, \\ & \text{ one to the left side of } B(n) \text{ that } h\text{-tunnels through} \\ & R_1, \text{ and one to the right side of } B(n) \text{ that } h\text{-tunnels} \\ & \text{ through } R_3, \text{ and } \exists \text{ two disjoint closed paths in } B(n) \\ & \text{ from } \mathbf{x}, \text{ one to the bottom of } B(n) \text{ that } v\text{-tunnels} \\ & \text{ through } R_2 \text{ and one to the top of } B(n) \text{ that } v\text{-tunnels} \\ & \text{ through } R_4. \text{ Further, there are vertical open crossings} \\ & \text{ of each of } R_1 \text{ and } R_3, \text{ and there are} \\ & \text{ horizontal closed crossings of each of } R_2 \text{ and } R_4. \end{aligned} \tag{2.8}$$

Again, by the proof of Kesten [4], we have the following.

LEMMA 4. *There are constants $C_2$ and $C_4$ such that uniformly for $\mathbf{x} \in B(n/4)$, $P(\mathcal{U}_\kappa(\mathbf{x},n)) \le C_\kappa P(\mathcal{T}_\kappa(\mathbf{x},n))$, $\kappa = 2,4$.*

2.1. *Proof of lower bounds.* We now obtain expectation lower bounds for the sizes of the lowest crossing, pioneering sites and pivotal sites. By definition, we have $|L_n| = \sum_{\mathbf{x}\in B(n)} \mathbb{1}_{\mathcal{L}(\mathbf{x},n)}$. Thus,

$$E|L_n| = \sum_{\mathbf{x}\in B(n)} P(\mathcal{L}(\mathbf{x},n)) \ge \sum_{\mathbf{x}\in B(n/4)} P(\mathcal{L}(\mathbf{x},n)). \tag{2.9}$$

By Lemmas 1 and 3 and the inclusion $\{\mathbf{x} \text{ is open}\} \cap \mathcal{T}_3(\mathbf{x},n) \subset \mathcal{L}(\mathbf{x},n)$, we have $2P(\mathcal{L}(\mathbf{x},n)) \ge (1/C_3)P(\mathcal{U}_3(\mathbf{x},n)) \ge (1/C_3)P(\mathcal{U}_3(\mathbf{0},5n/4)) \ge n^{-2/3+o(1)}$, uniformly for $\mathbf{x} \in B(n/4)$. Therefore, summing on $\mathbf{x} \in B(n/4)$ in (2.9), we obtain

$$E|L_n| \ge (n/4)^2 n^{-2/3+o(1)} = n^{4/3+o(1)}. \tag{2.10}$$

In exactly the same way, but using now Lemmas 2 and 4 and the inclusions $\{\mathbf{x} \text{ is open}\} \cap \mathcal{T}_2(\mathbf{x}, n) \subset \mathcal{F}(\mathbf{x}, n)$ and $\{\mathbf{x} \text{ is open}\} \cap \mathcal{T}_4(\mathbf{x}, n) \subset \mathcal{Q}(\mathbf{x}, n)$, we have

$$E|F_n| \geq (n/4)^2 n^{-1/4+o(1)} = n^{7/4+o(1)} \tag{2.11}$$

and

$$E|Q_n| \geq (n/4)^2 n^{-5/4+o(1)} = n^{3/4+o(1)}. \tag{2.12}$$

Note finally that the $\tau$th moment lower bounds in Theorem 1 follow immediately from (2.10)–(2.12) and Jensen's inequality for all $\tau \geq 1$.

**3. Lowest crossing and pioneering sites.** In this section we carefully study an upper bound for the first and second moments of $|L_n|$ and $|F_n|$. We do this to establish a dyadic summation construction alternative to Kesten's [3] method that we will later incorporate in our analysis of the pivotal case in Sections 5 and 6. We introduce the following concentric square annuli of vertices in $B(n)$. Let $j_0 = j_0(n)$ be the smallest integer $j$ such that $2^{-j} n \leq 1$. Define

$$\begin{aligned} &A_0 := B(n/2), \\ &A_j := B((1 - 2^{-(j+1)})n) \backslash B((1 - 2^{-j})n), \qquad 1 \leq j < j_0, \\ &A_{j_0} := B(n) \Big\backslash \bigcup_{j=0}^{j_0 - 1} A_j = \partial B(n). \end{aligned} \tag{3.1}$$

The annuli $A_j$ become thinner as they approach the boundary of $B(n)$ such that, for $j < j_0$, the distance from a point $\mathbf{x} \in A_j$ to $\partial B(n)$ is comparable with $2^{-j} n$ and also comparable with the width of $A_j$. Since $A_{j_0} = \partial B(n)$, we will use instead the property, valid for all $j \leq j_0$, that if $\mathbf{x} \in A_j$ and $\|\mathbf{y} - \mathbf{x}\| < 2^{-(j+1)} n$, then $\mathbf{y} \in B(n)$. Notice that the annuli $A_j$ are natural for an approach based on disjoint boxes by the following reasoning. For any vertex $\mathbf{x} \in B(n)$, we will construct a box $B(\mathbf{x}, r)$ centered at $\mathbf{x}$ that is roughly as large as it can be yet stays inside $B(n)$. The collection of vertices $\mathbf{x}$ that give rise to boxes $B(\mathbf{x}, r)$ with radii $r \asymp 2^{-j} n$ correspond to the annuli $A_j$. Therefore, roughly speaking, the sizes of largest disjoint boxes may be organized by arranging the centers of the boxes in these annuli.

If the sizes of the sets $A_j$ were defined by areas of the regions between concentric squares rather than by cardinalities of subsets of vertices of $B(n)$, we would obtain an upper bound for the sizes of these sets immediately by using the fact that $(2n)^2(1 - 2^{-(j+1)})^2 - (2n)^2(1 - 2^{-j})^2 \leq 2^{-j+2} n^2$. An error in approximating $|A_j|$ by the area between concentric squares may come about due to inclusion or exclusion of a ring of vertices. However, if $j < j_0$, then the thickness of a given annulus is $2^{-(j+1)} n \geq 1/2$ so the area estimate may only be an under-estimate by a factor of at most 4. Therefore, since the boundary of $B(n)$ has cardinality $8n - 4$, we have, for all $0 \leq j \leq j_0$, that

$$|A_j| \leq 2^{-j+4} n^2. \tag{3.2}$$

3.1. *Expectation upper bound.* We write the expectation of the size of the lowest crossing as

$$E|L_n| = \sum_{\mathbf{x}\in B(n)} P(L(\mathbf{x},n)) = \sum_{j=0}^{j_0}\sum_{\mathbf{x}\in A_j} P(\mathcal{L}(\mathbf{x},n)). \tag{3.3}$$

We note that by (2.2), for $\mathbf{x}\in A_j$,

$$P(\mathcal{L}(\mathbf{x},n)) \le P(\mathcal{U}_3(\mathbf{0},2^{-(j+1)}n)). \tag{3.4}$$

Let $\varepsilon > 0$. By Lemma 1, there exists a constant $C_\varepsilon$ such that

$$P(\mathcal{U}_3(\mathbf{0},r)) \le C_\varepsilon r^{-2/3+\varepsilon} \qquad \text{for all } r \ge 1. \tag{3.5}$$

Here and in the sequel we allow $C_\varepsilon$ to be a constant depending on $\varepsilon$ that may vary from appearance to appearance. Thus, by (3.2)–(3.5), we obtain

$$E|L_n| \le C_\varepsilon \sum_{j=0}^{j_0} |A_j|(2^{-(j+1)}n)^{-2/3+\varepsilon} \le C_\varepsilon n^{4/3+\varepsilon}\sum_{j=0}^{\infty} 2^{-j/3}. \tag{3.6}$$

Since the geometric series in (3.6) converges, we obtain by (3.6) that

$$E|L_n| \le n^{4/3+o(1)}. \tag{3.7}$$

By the same argument, we construct an upper bound for $E\,|F_n|$. Indeed, let $\varepsilon > 0$. By Lemma 2, there exists a constant $C_\varepsilon$ such that

$$P(\mathcal{U}_2(\mathbf{0},r)) \le C_\varepsilon r^{-1/4+\varepsilon} \qquad \text{for all } r \ge 1. \tag{3.8}$$

Therefore, just as in (3.3) and (3.6) but using now (3.8) in place of (3.5), we find $E|F_n| \le n^{7/4+o(1)}$. The proof of the upper bounds for $\tau = 1$ of items 1 and 2 of Theorem 1 is thus complete.

We comment that the exponent $(-1/4+\varepsilon)$ that takes the place of $(-2/3+\varepsilon)$ in the upper bound for $E|F_n|$ does not affect the convergence of the geometric sum because the exponents in (3.5) and (3.8) are greater than $-1$. Note, however, that a four arm calculation similar to that shown above would require the use of an exponent $(-5/4+\varepsilon)$ so that the corresponding geometric sum would not converge. It is precisely for this reason that we must establish an alternative to Kesten's [3] method of proof to obtain our Theorem 1 for the pivotal case. The approach we have shown above for the first moment upper bound may be extended, in fact, to all moments, though we will not show the general case due to the fact mentioned earlier that Kesten's method may be applied successfully to obtain a general $\tau$th moment bound in the one, two- and three-arm cases. We only show in addition below a second moment upper bound for the lowest crossing because it demonstrates the way we extend our dyadic summation method to higher moments in all cases, including the pivotal one.

3.2. *Second moment upper bound.* We show an estimation of the second moment of $|L_n|$. Write

$$E(|L_n|^2) = \sum_{\mathbf{x}\in B(n)} \sum_{\mathbf{y}\in B(n)} P(\mathcal{L}(\mathbf{x},n)\cap\mathcal{L}(\mathbf{y},n)). \tag{3.9}$$

Recall the definition of $j_0 = j_0(n)$ and $A_j$ in (3.1). Consider first the "diagonal" contribution to (3.9) defined by

$$I := \sum_{j=0}^{j_0} \sum_{\mathbf{x}\in A_j} \sum_{\mathbf{y}\in B(\mathbf{x},2^{-(j+2)}n)} P(\mathcal{L}(\mathbf{x},n)\cap\mathcal{L}(\mathbf{y},n)). \tag{3.10}$$

Fix $j \le j_0$ and $\mathbf{x} \in A_j$ and work on the inner sum in (3.10). For this purpose, we introduce a net of concentric square annuli $a_m = a_m(\mathbf{x})$ whose union is the box $B(\mathbf{x}, 2^{-(j+2)}n)$ as follows:

$$\begin{aligned} &a_m(\mathbf{x}) := B(\mathbf{x},2^{-m}n) \setminus B(\mathbf{x},2^{-(m+1)}n), \qquad j+2 \le m \le j_0-1,\\ &a_{j_0}(\mathbf{x}) := B(\mathbf{x},2^{-j_0}n). \end{aligned} \tag{3.11}$$

Notice that $a_{j_0}$ may consist of only the single point $\mathbf{x}$. By this decomposition, we have that $\sum_{m=j+2}^{j_0} \sum_{\mathbf{y}\in a_m} P(\mathcal{L}(\mathbf{x},n)\cap\mathcal{L}(\mathbf{y},n))$ is equal to the inner sum in (3.10). By (3.11), the size of $a_m$ is easily estimated by

$$|a_m| \le 2^{-2m+2}n^2 \qquad \text{all } j+2 \le m \le j_0. \tag{3.12}$$

Furthermore, for $\mathbf{x} \in A_j$ and $\mathbf{y} \in a_m$, with $j+2 \le m < j_0$, by halving the distance between $\mathbf{x}$ and $\mathbf{y}$, we have that

$$B(\mathbf{x},2^{-(m+2)}n) \cap B(\mathbf{y},2^{-(m+2)}n) = \varnothing. \tag{3.13}$$

Also, for $\mathbf{y} \in a_m$ with $m \ge j+2$, since $\|\mathbf{y}-\mathbf{x}\| \le 2^{-(j+2)}n$ and $2^{-(m+2)}n + 2^{-(j+2)}n < 2^{-(j+1)}n$, we have that both $B(\mathbf{x},2^{-(m+2)}n)$ and $B(\mathbf{y},2^{-(m+2)}n)$ are subsets of $B(n)$. Therefore, since $\mathcal{L}(\mathbf{x},n)\cap\mathcal{L}(\mathbf{y},n)$ implies that for each of the boxes in (3.13) there exists a three-arm path from the center of the box to its boundary, we have by (3.13), independence and (3.5) that, for all $\mathbf{y} \in a_m$ with $m < j_0$,

$$P(\mathcal{L}(\mathbf{x},n)\cap\mathcal{L}(\mathbf{y},n)) \le C_\varepsilon (2^{-m}n)^{-4/3+2\varepsilon}. \tag{3.14}$$

Also, trivially, (3.14) continues to hold with $m = j_0$, since then $2^{-m}n \ge 1/2$. Thus, by (3.14), we have

$$\begin{aligned} \sum_{m=j}^{j_0} \sum_{\mathbf{y}\in a_m} P(\mathcal{L}(\mathbf{x},n)\cap\mathcal{L}(\mathbf{y},n)) &\le C_\varepsilon \sum_{m=j}^{\infty} 2^{-2m}n^2(2^{-m}n)^{-4/3+2\varepsilon}\\ &\le C_\varepsilon n^{2/3+2\varepsilon}2^{-2j/3}. \end{aligned} \tag{3.15}$$

Therefore, by (3.10), (3.15) and (3.2), we have

$$I \le C_\varepsilon n^{2/3+2\varepsilon} \sum_{j=0}^{\infty} |A_j| 2^{-2j/3} \le C_\varepsilon n^{8/3+2\varepsilon} \sum_{j=0}^{\infty} 2^{-5j/3} \le C_\varepsilon n^{8/3+2\varepsilon}. \tag{3.16}$$

We next consider the off-diagonal part of the sum (3.9) defined by

$$II := \sum_{j=0}^{j_0} \sum_{\mathbf{x}\in A_j} \sum_{k=0}^{j_0} \sum_{\mathbf{y}\in A_k} \chi_{\{\|\mathbf{x}-\mathbf{y}\|>2^{-j-2}n\}} P(\mathcal{L}(\mathbf{x},n)\cap \mathcal{L}(\mathbf{y},n)). \tag{3.17}$$

Here $\chi_A$ denotes the indicator function of the given set of vertices $A$. In the sum $II$, for all eligible vertices $\mathbf{x}\in A_j$ and $\mathbf{y}\in A_k$ with $k \ge j$, we have that the boxes $B(\mathbf{x}, 2^{-j-3}n)$ and $B(\mathbf{y}, 2^{-k-3}n)$ are disjoint and lie inside $B(n)$, while if $k < j$, the same holds true by (3.1) when we replace these boxes, respectively, by $B(\mathbf{x}, 2^{-j-3}n)$ and $B(\mathbf{y}, 2^{-k-5}n)$. Therefore, by (3.5), (3.17) and (3.2) we have

$$\begin{aligned} II &\le C_\varepsilon \sum_{j=0}^{j_0} |A_j| \sum_{k=0}^{j_0} |A_k| (2^{-j}n)^{-2/3+\varepsilon} (2^{-k}n)^{-2/3+\varepsilon} \\ &= C_\varepsilon \left( \sum_{j=0}^{j_0} |A_j| (2^{-j}n)^{-2/3+\varepsilon} \right)^2 \le C_\varepsilon n^{8/3+2\varepsilon}. \end{aligned} \tag{3.18}$$

Therefore, by definitions (3.9), (3.10) and (3.17), and by collecting the estimates (3.16) and (3.18), we obtain

$$E(|L_n|^2) \le n^{8/3+o(1)}. \tag{3.19}$$

We handle an upper bound for the second moment of the number of pioneering sites by the same method. Thus, we have established the upper bound for $\tau = 2$ of items 1 and 2 of Theorem 1. This concludes our discussion of these items.

**4. Method of disjoint boxes.** Denote $p_{n,\tau}(\mathbf{x}_1, \mathbf{x}_2, \ldots, \mathbf{x}_\tau) = P(\bigcap_{i=1}^{\tau} \mathcal{Q}(\mathbf{x}_i, n))$. Recall the definition of $j_0 = j_0(n)$ and $A_j$ in (3.1). Define, for all $j_1 \le j_2 \le \cdots \le j_\tau$,

$$\Sigma_{j_1,j_2,\ldots,j_\tau} := \sum_{\mathbf{x}_1\in A_{j_1}} \sum_{\mathbf{x}_2\in A_{j_2}} \cdots \sum_{\mathbf{x}_\tau\in A_{j_\tau}} p_{n,\tau}(\mathbf{x}_1, \ldots, \mathbf{x}_\tau). \tag{4.1}$$

By symmetry, to obtain an upper bound for the $\tau$th moment of the number of pivotal sites, it suffices to estimate the sum

$$\Sigma_0 := \sum_{j_1=0}^{j_0} \sum_{j_2=j_1}^{j_0} \cdots \sum_{j_\tau=j_{\tau-1}}^{j_0} \Sigma_{j_1,j_2,\ldots,j_\tau}. \tag{4.2}$$

Moreover, by induction on $\tau$ in Theorem 1, we may assume that all the vertices in (4.1)–(4.2) are distinct. In this section we establish a parametrization of certain boxes centered at the vertices $\mathbf{x}_1, \mathbf{x}_2, \dots, \mathbf{x}_\tau$ that are both mutually disjoint and large enough to yield convergence of the sum (4.2) in our method for estimating this sum shown in Sections 5 and 6. Indeed, we are led naturally to a certain graph $\mathbf{G}$ defined below that organizes the vertices and their relative distances from one another. Although this organization is somewhat complicated, it will allow us to introduce estimations of $p_{n,\tau}(\mathbf{x}_1, \mathbf{x}_2, \dots, \mathbf{x}_\tau)$ that refine the estimation approach based solely on disjoint boxes (illustrated in Section 3.2) because our estimation will depend also on the configuration of the graph.

We lay the groundwork for the definition of the graph $\mathbf{G}$ as follows. Let $c \geq 2$ be a positive integer depending only on $\tau$ that we will specify later. We say a vertex $\mathbf{v}$ is "near to" a vertex $\mathbf{u}$, for some $\mathbf{u} \in A_j$, if $\mathbf{v} \in B(\mathbf{u}, 2^{-j-2c}n)$, and write this (asymmetric) relation as $\mathbf{v}N\mathbf{u}$. If $\mathbf{v}$ is not near to $\mathbf{u}$, we write instead $\mathbf{v}\widetilde{N}\mathbf{u}$. Let now $\mathbf{x}_i \in A_{j_i}$, $1 \leq i \leq \tau$, with $j_1 \leq j_2 \leq \cdots \leq j_\tau$. We say that a sequence of vertices $\mathbf{x}_{f_1}, \dots, \mathbf{x}_{f_k}$ is a chain that leads from $\mathbf{x}_{f_1}$ to $\mathbf{x}_{f_k}$ if $\mathbf{x}_{f_i}N\mathbf{x}_{f_{i+1}}$ for each $i = 1, \dots, k-1$. Define $e_1 := 1$ and

$$V_1 := \{\mathbf{x}_e : e \neq e_1, \text{ and there exists a chain from } \mathbf{x}_e \text{ to } \mathbf{x}_{e_1}\}.$$

Thus, $V_1$ is the set of all vertices that lead to $\mathbf{x}_{e_1}$. Note that $\mathbf{x}_1$ may lead to $\mathbf{x}_2$, but if $\mathbf{x}_2$ does not lead to $\mathbf{x}_1$, then $\mathbf{x}_{e_2} \notin V_1$. We denote inductively, by $E_i := \{e : \mathbf{x}_e \in V_i\}$, the set of indices corresponding to vertices in $V_i$, $i = 1, 2, \dots$, that we now continue to define. Note that the cardinalities of $E_i$ and $V_i$ are the same since we have assumed the vertices $\mathbf{x}_e$ are distinct. Let $e_2$ be the smallest index with $e_2 > e_1$ such that $e_2 \notin E_1$. Define

$$V_2 := \{\mathbf{x}_e : e \notin (\{e_1, e_2\} \cup E_1), \text{ and there exists a chain from } \mathbf{x}_e \text{ to } \mathbf{x}_{e_2}\}.$$

Thus, no element of $V_2$ begins a chain that leads to $\mathbf{x}_{e_1}$. It may be that $\mathbf{x}_{e_1}$ leads to $\mathbf{x}_{e_2}$, but we leave $\mathbf{x}_{e_1}$ out of $V_2$ as defined. Continuing in this fashion, we take $e_3$ to be the smallest index with $e_3 > e_2$ such that $e_3 \notin (E_1 \cup E_2)$. Define

$$V_3 := \{\mathbf{x}_e : e \notin (\{e_1, e_2, e_3\} \cup E_1 \cup E_2), \text{ and } \exists \text{ a chain from } \mathbf{x}_e \text{ to } \mathbf{x}_{e_3}\}.$$

Finally, we obtain a disjoint collection of sets of vertices $V_1, \dots, V_r$, where some of the $V_i$ may be empty. We say that $V_i$ is the set of vertices chained to the root $\mathbf{x}_{e_i}$.

Thus, for example, if $\tau = 3$ and both $\mathbf{x}_3N\mathbf{x}_1$ and $\mathbf{x}_2N\mathbf{x}_1$, then $V_1 = \{\mathbf{x}_2, \mathbf{x}_3\}$ and $e_2$ is undefined. Also if $\mathbf{x}_3N\mathbf{x}_1$ and $\mathbf{x}_2\widetilde{N}\mathbf{x}_1$ but instead $\mathbf{x}_2N\mathbf{x}_3$, then again $V_1 = \{\mathbf{x}_2, \mathbf{x}_3\}$ and $e_2$ is undefined. If, on the other hand, $\mathbf{x}_2N\mathbf{x}_1$, $\mathbf{x}_3\widetilde{N}\mathbf{x}_1$ and $\mathbf{x}_3\widetilde{N}\mathbf{x}_2$, then $V_1 = \{\mathbf{x}_2\}$ and $e_2 = 3$ and $V_2$ is empty. Further, if $\mathbf{x}_2\widetilde{N}\mathbf{x}_1$, $\mathbf{x}_3\widetilde{N}\mathbf{x}_1$ and $\mathbf{x}_3N\mathbf{x}_2$, then $V_1$ is empty and $e_2 = 2$ and $V_2 = \{\mathbf{x}_3\}$. Finally, if $\mathbf{x}_2\widetilde{N}\mathbf{x}_1$, $\mathbf{x}_3\widetilde{N}\mathbf{x}_1$ and $\mathbf{x}_3\widetilde{N}\mathbf{x}_2$, then $e_i = i$ and $V_i$ is empty, $i = 1, 2, 3$.

Suppose now, in general, that $e_i$ is defined for $i=1,\dots,r$. Thus, $r$ is the number of root vertices. Note by definition that the vertex $\mathbf{x}_1$ is always counted among the roots. We say that a vertex $\mathbf{x}_{e_i}$ is isolated if $V_i=\varnothing$. At nonisolated roots we introduce a decomposition of the sets $V_i$ themselves by means of a local "near to" relation. It turns out that we will be able to work with one original root $\mathbf{x}_{e_i}$ and its corresponding set of vertices $V_i$ at a time in constructing the moment estimates of Sections 5 and 6, so in what follows we only write out a decomposition of $V_1$. We will represent this decomposition as a graph $\mathbf{G}_1$ below, where, in general, a connected graph $\mathbf{G}_i$ with vertex set $\{\mathbf{x}_{e_i}\}\cup V_i$ is associated with the root vertex $\mathbf{x}_{e_i}$. The graph $\mathbf{G}$ on all vertices is defined simply as the union of the component graphs $\mathbf{G}_i$.

Let $|V_1|\geq 1$. We denote $V_1:=\{\mathbf{y}_1,\mathbf{y}_2,\dots\}$, where the names of the vertices have been changed such that $\mathbf{y}_1\in a_{p_1}(\mathbf{x}_1),\mathbf{y}_2\in a_{p_2}(\mathbf{x}_1),\dots$, for $p_1\geq p_2\geq\cdots$, where $\mathbf{y}_1$ is determined such that $\mathbf{y}_1 N\mathbf{x}_1$ and such that $\mathbf{y}_1$ minimizes the distance to $\mathbf{x}_1$. Therefore, $p_1\geq j_1+2c$ by (3.11). We now say that $\mathbf{v}M_{\mathbf{w}}\mathbf{u}$ for some $\mathbf{u}\in a_m(\mathbf{w})$, if $\mathbf{v}\in B(\mathbf{u},2^{-m-2c}n)$. We call $M_{\mathbf{w}}$ a local relation, where $\mathbf{w}$ fixes the locale of the relation. We describe how to decompose $V_1$ via a collection of local relations starting with $M_{\mathbf{x}_1}$ in a way wholly similar to the decomposition of the original set of vertices $\{\mathbf{x}_1,\dots,\mathbf{x}_\tau\}$ via the $N$-relation. Indeed, set $f_1=1$, rename $\mathbf{y}_1$ as $\mathbf{w}_1$, and define $W_1$ as all the vertices of $V_1\setminus\{\mathbf{w}_1\}$ that are chained to the root $\mathbf{w}_1$ by means of a chain of relations for the relation $M_{\mathbf{x}_1}$. We rename $p_{f_1}=m_1$ so that $\mathbf{w}_1\in a_{m_1}(\mathbf{x}_1)$. Let $f_2$ be the smallest index with $f_2>f_1$ such that $\mathbf{y}_{f_2}\notin W_1$. We rename $\mathbf{y}_{f_2}=\mathbf{w}_2$ and also $p_{f_2}=m_2$ so that $\mathbf{w}_2\in a_{m_2}(\mathbf{x}_1)$. Note, in particular, that $\mathbf{w}_2\widetilde{M}_{\mathbf{x}_1}\mathbf{w}_1$. Define

$$W_2:=\{\mathbf{y}_f:\mathbf{y}_f\notin(\{\mathbf{w}_1,\mathbf{w}_2\}\cup W_1),\ \text{and}\ \exists\,\text{a}\ M_{\mathbf{x}_1}\text{-chain from}\ \mathbf{y}_f\ \text{to}\ \mathbf{w}_2\}.$$

Continuing in this way, we define also $f_3<f_4<\cdots$ as long as these exist and so also local roots $\mathbf{w}_i=\mathbf{y}_{f_i}$ and corresponding sets $W_i$, $i=3,4,\dots,$ chained to them by the relation $M_{\mathbf{x}_1}$. We also define indices for the locations of the local roots. Indeed, following the example above for our definitions of $m_1$ and $m_2$, we define $m_i$ such that $\mathbf{w}_i\in a_{m_i}(\mathbf{x}_1)$ for all $i$ such that $\mathbf{w}_i$ exists. Note by definition, since $m_i=p_{f_i}$ and $p_1\geq p_2\geq\cdots$, we have that $m_1\geq m_2\geq\cdots$. In general, for each $i$, we further decompose the set $W_i$ into a disjoint union:

$$(\{\mathbf{w}_{i,1}\}\cup W_{i,1})\cup(\{\mathbf{w}_{i,2}\}\cup W_{i,2})\cup\cdots.$$

Here for each $j=1,2,\dots$, $W_{i,j}$ is a set of vertices chained to the corresponding local root $\mathbf{w}_{i,j}$ by the relation $M_{\mathbf{w}_i}$ as follows. Assume $W_1$ is not empty, else $\mathbf{w}_{1,1}$ and $W_{1,1}$ are undefined. Since $W_1$ is the set of elements chained to $\mathbf{w}_1$, we know there exists $\mathbf{y}\in W_1$ such that $\mathbf{y}M_{\mathbf{x}_1}\mathbf{w}_1$. We take $\mathbf{w}_{1,1}$ as such a vertex $\mathbf{y}$ that minimizes the distance to $\mathbf{w}_1$. We define the index $m_{1,1}$ by the property: $\mathbf{w}_{1,1} \in a_{m_{1,1}}(\mathbf{w}_1)$. Therefore, by definition of $M_{\mathbf{x}_1}$ and the fact that $\mathbf{w}_1 \in a_{m_1}(\mathbf{x}_1)$, we have $m_{1,1} \geq m_1 + 2c$. Note by definition of $\mathbf{w}_{1,1}$, that for any $\mathbf{y} \in W_1$, we have $\mathbf{y} \in a_p(\mathbf{w}_1)$ with some $p \leq m_{1,1}$. We define $W_{1,1}$ as the set of vertices in $W_1 \setminus \{\mathbf{w}_{1,1}\}$ that are chained to the local root $\mathbf{w}_{1,1}$ by the relation $M_{\mathbf{w}_1}$. We perform a similar procedure starting with $W_2$ to define the local root $\mathbf{w}_{2,1}$. In particular, $\mathbf{w}_{2,1} M_{\mathbf{x}_1} \mathbf{w}_2$. Likewise, as long as $W_i$ is not empty, we define $\mathbf{w}_{i,1} \in W_i$ and a set $W_{i,1}$ chained to $\mathbf{w}_{i,1}$ by the relation $M_{\mathbf{w}_i}$. Here the indices $m_{2,1}, m_{3,1}, \ldots$ are defined such that $\mathbf{w}_{i,1} \in a_{m_{i,1}}(\mathbf{w}_i)$, $i = 2, 3, \ldots$. Again, we choose $\mathbf{w}_{i,1}$ such that $m_{i,1}$ is maximal, that is, there does not exist $\mathbf{y} \in a_m(\mathbf{w}_i) \cap W_i$ with $m > m_{i,1}$.

We define $\mathbf{w}_{1,2}$ and $W_{1,2}$ if $W_1 \setminus (W_{1,1} \cup \{\mathbf{w}_{1,1}\})$ is not empty. We do this in the same way that we defined $\mathbf{w}_2$ and $W_2$ from $V_1 \setminus (W_1 \cup \{\mathbf{w}_1\})$. Thus, we order the vertices in $W_1$ as $\mathbf{w}_{1,1}, y_2, y_3, \ldots$, where $y_i \in a_{p_i}(\mathbf{w}_1)$ with $m_{1,1} \geq p_2 \geq p_3 \geq \cdots$. Among all elements of $W_1 \setminus \{\mathbf{w}_{1,1}\}$ that are not chained to the local root $\mathbf{w}_{1,1}$ by the relation $M_{\mathbf{w}_1}$, we choose $\mathbf{w}_{1,2}$ to be the vertex $\mathbf{y}_i$ with least index. Correspondingly, we define $W_{1,2}$ as the elements of $W_1 \setminus (W_{1,1} \cup \{\mathbf{w}_{1,1}, \mathbf{w}_{1,2}\})$ that are chained to $\mathbf{w}_{1,2}$ by the relation $M_{\mathbf{w}_1}$. Note, in particular, that $\mathbf{w}_{1,2} \in a_{m_{1,2}}(\mathbf{w}_1)$ with $m_{1,2} \leq m_{1,1}$. Similarly, we define local roots $\mathbf{w}_{1,j}, j = 3, 4, \ldots$ and for each $i \geq 2$, the roots $\mathbf{w}_{i,j}, j = 2, 3, \ldots$. We define $m_{i_1,i_2}$, for $i_2 = 1, 2, \ldots$, as indices such that $\mathbf{w}_{i_1,i_2} \in a_m(\mathbf{w}_{i_1})$ with $m = m_{i_1,i_2}$. Again by definition $m_{i,1} \geq m_{i,2} \geq \cdots$, and, since $\mathbf{w}_{i,1} \in a_{m_{i,1}}(\mathbf{w}_i)$ and $\mathbf{w}_{i,1} M_{\mathbf{x}_1} \mathbf{w}_i$, we have that $m_{i,1} \geq m_i + 2c$.

We inductively continue this procedure such that for a given local root $\mathbf{w}_{i_1,\ldots,i_k}$ and associated local elements $W_{i_1,\ldots,i_k}$ chained to it, we decompose

$$W_{i_1,\ldots,i_k} = (\{\mathbf{w}_{i_1,\ldots,i_k,1}\} \cup W_{i_1,\ldots,i_k,1}) \cup (\{\mathbf{w}_{i_1,\ldots,i_k,2}\} \cup W_{i_1,\ldots,i_k,2}) \cup \cdots$$

by means of the relation $M_{\mathbf{w}}$ for $\mathbf{w} = \mathbf{w}_{i_1,\ldots,i_k}$. We continue in this way until no further local roots may be defined. In general, for $k \geq 1$, we have $\mathbf{w}_{i_1,\ldots,i_k,1} \in a_m(\mathbf{w}_{i_1,\ldots,i_k})$ for $m = m_{i_1,\ldots,i_k,1}$ and

$$\mathbf{w}_{i_1,\ldots,i_k,1} M_{\mathbf{w}_{i_1,\ldots,i_{k-1}}} \mathbf{w}_{i_1,\ldots,i_k}.$$

Here when $k = 1$, $\mathbf{w}_{i_1,\ldots,i_{k-1}}$ becomes $\mathbf{x}_1$. We also define the index $m_{i_1,\ldots,i_{k+1}}$, in general, by the property that $\mathbf{w}_{i_1,\ldots,i_{k+1}} \in a_m(\mathbf{w}_{i_1,\ldots,i_k})$ with $m = m_{i_1,\ldots,i_{k+1}}$. We have that, for all $k \geq 0$, $m_{i_1,\ldots,i_k,1} \geq m_{i_1,\ldots,i_k,2} \geq \cdots$ and $m_{i_1,\ldots,i_k,1} \geq m_{i_1,\ldots,i_k} + 2c$, where for $k = 0$, $m_{i_1,\ldots,i_k}$ denotes $j_1$.

We now use our parameter $c$ to obtain one further property of the indices not mentioned in the previous paragraph. First, since $W_1$ consists of all vertices $\mathbf{w}$ that may be chained to $\mathbf{w}_1$ by the relation $M_{\mathbf{x}_1}$, we argue that $c$ may be chosen such that

$$\mathbf{w} \in a_p(\mathbf{x}_1) \qquad \text{with } p \geq m_1 - 1, \text{ for all } \mathbf{w} \in W_1. \tag{4.3}$$

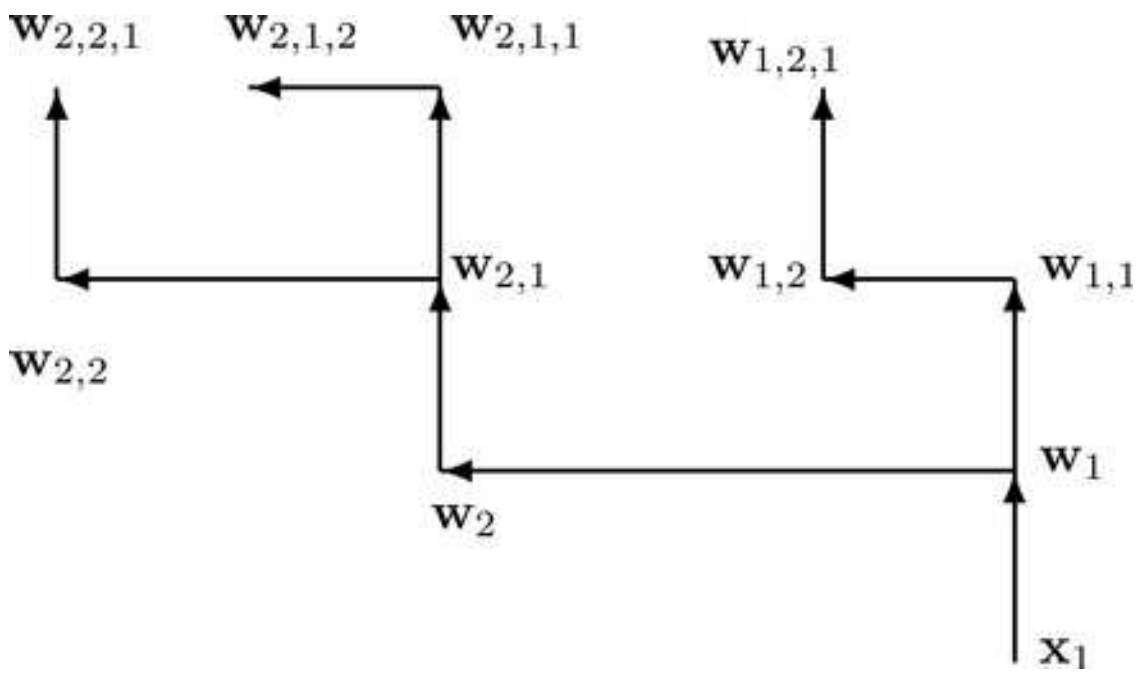

FIG. 1. *The graph* $\mathbf{G}_1$.

Indeed, since there are at most $\tau - 1$ relations with respect to $M_{\mathbf{x}_1}$ that must be satisfied, if $c$ is large enough and if $\mathbf{w} \in a_m(\mathbf{x}_1)$ for some $m \le m_1 - 2$, then the chain will not be able to cross the square annulus $a_{m_1-1}(\mathbf{x}_1)$ to reach $\mathbf{w}_1 \in a_{m_1}(\mathbf{x}_1)$. Therefore, we choose $c$ sufficiently large to guarantee (4.3). Note that the value of $c$ so chosen does not depend on the value $m_1$ or the location $\mathbf{x}_1$. Now, since we have control on the index $p$ for the location of $\mathbf{w}$, it is easy by estimating the sum of distances between successive vertices in a chain of relations leading to $\mathbf{w}_1$ by $(\tau-1)2^{-m_1+1-2c}n$ that again, by choosing $c$ somewhat larger if necessary, we have $m_{1,i_2} \ge m_1 + c$ for all $i_2 = 2, 3 \ldots,$ while, of course, we still have that $m_{1,1} \ge m_1 + 2c$. By the same argument based at any local root, we have, for all $k \ge 0$, that

$$m_{i_1,\ldots,i_k,i_{k+1}} \ge m_{i_1,\ldots,i_k} + c \qquad \text{for } i_{k+1} = 2, 3, \ldots, \tag{4.4}$$
$$m_{i_1,\ldots,i_k,1} \ge m_{i_1,\ldots,i_k} + 2c.$$

We now define the graph $\mathbf{G}_1$ alluded to above. The vertices of $\mathbf{G}_1$ are $\{\mathbf{x}_1\} \cup V_1$. We assume $|V_1| \ge 1$, else the graph is trivial. If $\mathbf{w}_{i_1,\ldots,i_k}$ is a local root, we say that the root is at level $k$. We define a (horizontal) edge at level $k$ between $\mathbf{w}_{i_1,\ldots,i_k}$ and $\mathbf{w}_{i_1,\ldots,i_k+1}$ whenever both these local roots exist. In our diagram below we make the edge go horizontally to the left from $\mathbf{w}_{i_1,\ldots,i_k}$ to $\mathbf{w}_{i_1,\ldots,i_k+1}$ to recall the fact that the associated indices satisfy $m_{i_1,\ldots,i_k+1} \le m_{i_1,\ldots,i_k}$. Next, we call $\mathbf{w}_1$ an immediate successor of $\mathbf{x}_1$. Similarly, if $\mathbf{w}_{i_1,\ldots,i_k}$ is a local root of level $k$ and if the root $\mathbf{w}_{i_1,\ldots,i_k,1}$ exists at level $k+1$, then we call this local root the immediate successor of the former local root. We now define that a (vertical) edge exists between two immediate successors. In our diagram the level increases vertically with $k$. We illustrate $\mathbf{G}_1$ for the following example in Figure 1: $|V_1| = 10$, $|W_1| = 3$, $|W_{1,1}| = 0$, $|W_{1,2}| = 1$ and $|W_{1,2,1}| = 0$; $|W_2| = 5$, $|W_{2,1}| = 2$, $|W_{2,2}| = 1$, $|W_{2,2,1}| = 0$, $|W_{2,1,1}| = 0$ and $|W_{2,1,2}| = 0$.

4.1. *Representation of disjoint boxes.* We fix the graph $\mathbf{G}_1$ and study the problem of verifying that certain boxes centered at its vertices that we now construct are indeed disjoint. Assume $|V_1| \geq 1$. For each vertex $\mathbf{w} = \mathbf{w}_{i_1,\dots,i_k} \in \mathbf{G}_1$, we define

$$m(\mathbf{w}) = \begin{cases} m_{i_1,i_2,\dots,i_k,1}, & \text{if } \mathbf{w}_{i_1,i_2,\dots,i_k,1} \text{ exists,} \\ m_{i_1,i_2,\dots,i_k}, & \text{if } \mathbf{w}_{i_1,i_2,\dots,i_k,1} \text{ does not exist.} \end{cases} \tag{4.5}$$

Note that since we assume that $\mathbf{w}_1$ exists, we also have $m(\mathbf{x}_1) = m_1$. We set the constant value $s := 2c + 4$ where $c$ appears in (4.4).

PROPOSITION 1. *The collection of boxes* $B(\mathbf{w}, 2^{-m(\mathbf{w})-s}n)$, $\mathbf{w} \in \mathbf{G}_1$, *are mutually disjoint.*

PROOF. Let $k \geq 0$ and let $\mathbf{w} = \mathbf{w}_{i_1,\dots,i_k}$ and $\mathbf{w}' = \mathbf{w}_{i'_1,\dots,i'_{k'}}$ be distinct vertices in $\mathbf{G}_1$. Define $l$ as the largest nonnegative integer such that $i'_1 = i_1, \dots, i'_l = i_l$ and set $\mathbf{z} := \mathbf{w}_{i_1,\dots,i_l}$. If one of $\mathbf{w}$ or $\mathbf{w}'$ is $\mathbf{x}_1$, then we set $l = 0$ and put $\mathbf{z} = \mathbf{x}_1$. We consider two cases, namely, (a) one of $\mathbf{w}$ or $\mathbf{w}'$ is equal to $\mathbf{z}$, or (b) neither $\mathbf{w}$ nor $\mathbf{w}'$ is equal to $\mathbf{z}$. In case (a) we assume, without loss of generality, that $\mathbf{w} = \mathbf{z}$. Note therefore that with this choice in case (a), $l = k$, $k' > k$, and $\mathbf{w}' \in W_{i_1,\dots,i_k}$ since $\mathbf{w}_{i_1,\dots,i_k,1}$ exists. In case (b) we must have both $k \geq l+1$ and $k' \geq l+1$, else we are in case (a) again. Thus, in case (b) we may switch the designation of the primed vertex if necessary such that $i'_{l+1} > i_{l+1}$.

We work first with case (a). Since $\mathbf{w}' \in W_{i_1,\dots,i_k}$, we have that

$$\mathbf{w}' \in a_p(\mathbf{w}) \qquad \text{for some } p \leq m_{i_1,\dots,i_k,1}. \tag{4.6}$$

Indeed, (4.6) holds by the definition of $m_{i_1,\dots,i_k,1}$ as the maximal $m$ such that $\mathbf{y} \in a_m(\mathbf{w})$ among all $\mathbf{y} \in W_{i_1,\dots,i_k} \setminus \{\mathbf{w}\}$. Now by (4.5), we have $m(\mathbf{w}) = m_{i_1,\dots,i_k,1}$, so by (4.6),

$$B(\mathbf{w}', 2^{-m(\mathbf{w})-2}n) \quad \text{and} \quad B(\mathbf{w}, 2^{-m(\mathbf{w})-2}n) \qquad \text{are disjoint.} \tag{4.7}$$

Consider first a special case of (a), namely, that

$$\mathbf{w}' \in \{\mathbf{w}_{i_1,\dots,i_k,1}\} \cup W_{i_1,\dots,i_k,1} \tag{4.8}$$

so that $\mathbf{w}'$ is either the immediate successor of $\mathbf{w}$ or is one of the descendants of this immediate successor. It follows by (4.8) that $i'_{k+1} = 1$. Therefore, since $k' \geq k+1$, we have, by (4.4) and (4.8), that

$$m_{i'_1,i'_2,\dots,i'_{k'},1} \geq m_{i'_1,i'_2,\dots,i'_{k'}} + 2c \geq m_{i_1,i_2,\dots,i_k,1} + 2c.$$

Therefore, by (4.5), we have that $m(\mathbf{w}') \geq m(\mathbf{w})$ whether or not the vertex $\mathbf{w}_{i'_1,i'_2,\dots,i'_{k'},1}$ exists. Hence, it follows by (4.7) that

$$B(\mathbf{w}', 2^{-m(\mathbf{w}')-s}n) \quad \text{and} \quad B(\mathbf{w}, 2^{-m(\mathbf{w})-s}n) \qquad \text{are disjoint.}$$

Thus, we have established disjoint boxes under condition (4.8) in case (a).

Suppose next for case (a) that $k' \geq k+1$ with $i'_{k+1} \geq 2$. Thus, we consider the remaining descendants $\mathbf{w}'$ of $\mathbf{w}$ that were not considered in the special case (4.8). Put $\widetilde{\mathbf{w}} := \mathbf{w}_{i_1,\ldots,i_k,i'_{k+1}}$, so either $\mathbf{w}' = \widetilde{\mathbf{w}}$ (when $k' = k+1$) or $\mathbf{w}'$ is a descendant of $\widetilde{\mathbf{w}}$:

$$\mathbf{w}' \in \{\widetilde{\mathbf{w}}\} \cup W_{i_1,\ldots,i_k,i'_{k+1}}.$$

For all such $\mathbf{w}'$, we have that

$$\mathbf{w}' \in a_p(\mathbf{w}) \qquad \text{for some } p \leq m_{i_1,\ldots,i_k,i'_{k+1}}. \tag{4.9}$$

Indeed, by definition, the vertex $\widetilde{\mathbf{w}}$ lies in the annulus $a_m(\mathbf{w})$ where $m$ is maximal: if some $\mathbf{y} \in W_{i_1,\ldots,i_k} \setminus \bigcup_{i=1}^{i'_{k+1}-1} (W_{i_1,\ldots,i_k,i} \cup \{\mathbf{w}_{i_1,\ldots,i_k,i}\})$ lies also in $a_p(\mathbf{w})$, then $m \geq p$. Therefore, since indeed $\mathbf{w}'$ is one such vertex $\mathbf{y}$, the assertion (4.9) is verified. Hence, by (3.13),

$$B(\mathbf{w}', 2^{-p-2}n) \quad \text{and} \quad B(\mathbf{w}, 2^{-p-2}n) \qquad \text{are disjoint.}$$

But, by definition of the indices and (4.9), we have

$$m(\mathbf{w}) = m_{i_1,\ldots,i_k,1} \geq m_{i_1,\ldots,i_k,i'_{k+1}} \geq p.$$

Also, by (4.4) and (4.5), since $l = k$ and $k' \geq k+1$, we have that

$$m(\mathbf{w}') \geq m_{i'_1,\ldots,i'_{k'}} \geq m_{i_1,\ldots,i_k,i'_{k+1}} + c \geq p + c \geq p.$$

Thus, since $s \geq 2$, we obtain the desired conclusion. This completes the proof of disjoint boxes for case (a).

We now proceed to study case (b). We first note that since $k, k' \geq l+1$,

$$\begin{aligned} &\text{(i)} \;\; \mathbf{w} \in \{\mathbf{w}_{i_1,\ldots,i_l,i_{l+1}}\} \cup W_{i_1,\ldots,i_l,i_{l+1}}, \\ &\text{(ii)} \;\; \mathbf{w}' \in \{\mathbf{w}_{i_1,\ldots,i_l,i'_{l+1}}\} \cup W_{i_1,\ldots,i_l,i'_{l+1}}. \end{aligned} \tag{4.10}$$

Therefore, just as in (4.9), we find by (4.10) that

$$\begin{aligned} &\text{(i)} \;\; \mathbf{w} \in a_p(\mathbf{z}) &&\text{for some } p \leq m_{i_1,\ldots,i_l,i_{l+1}}, \\ &\text{(ii)} \;\; \mathbf{w}' \in a_{p'}(\mathbf{z}) &&\text{for some } p' \leq m_{i_1,\ldots,i_l,i'_{l+1}}. \end{aligned} \tag{4.11}$$

Now we claim that for $p$ given in (4.11), we have

$$\mathbf{w}' \notin B(\mathbf{w}, 2^{-p-2c}n). \tag{4.12}$$

Indeed, on the contrary, we would have $\mathbf{w}' M_{\mathbf{z}} \mathbf{w}$. Therefore, we could chain $\mathbf{w}'$ to $\mathbf{w}_{i_1,\ldots,i_l,i_{l+1}}$ by the relation $M_{\mathbf{z}}$. Indeed, if $l \leq k-2$, then $\mathbf{w}$ is already chained in this way to $\mathbf{w}_{i_1,\ldots,i_l,i_{l+1}}$, while if $l = k-1$, then $\mathbf{w} = \mathbf{w}_{i_1,\ldots,i_l,i_{l+1}}$ so we would have directly that $\mathbf{w}' M_{\mathbf{z}} \mathbf{w}_{i_1,\ldots,i_l,i_{l+1}}$. Therefore, on the one hand, we have the inclusion (ii) of (4.10) and, on the other hand, we would have that $\mathbf{w}' \in \bigcup_{i=1}^{i_{l+1}} W_{i_1,\dots,i_l,i}$ since $\mathbf{w}'$ is chained to $\mathbf{w}_{i_1,\dots,i_l,i_{l+1}}$. But these two inclusions are in contradiction since $i'_{l+1} > i_{l+1}$. Hence, we must not have that this chain relation exists and, therefore, (4.12) holds.

To finish the argument for case (b), suppose first that $p' < p-2$. Then by (4.11) alone and (3.11), we have that $B(\mathbf{w}', 2^{-p'-2}n)$ and $B(\mathbf{w}, 2^{-p-2}n)$ are disjoint. But by (4.4) and (4.5),

$$m(\mathbf{w}') \geq m_{i_1,\dots,i_l,i'_{l+1}} \geq p' \quad \text{and} \quad m(\mathbf{w}) \geq m_{i_1,\dots,i_l,i_{l+1}} \geq p.$$

Thus, since $s \geq 2$, we obtain the desired disjoint boxes condition. Suppose finally that $p' \geq p-2$. We have by (4.12) and (3.13) that $B(\mathbf{w}', 2^{-p-2c-2}n)$ and $B(\mathbf{w}, 2^{-p-2c-2}n)$ are disjoint. Therefore, we obtain the disjoint boxes condition by using (4.4), (4.5) and (4.11) to obtain the following two strings of inequalities:

$$m(\mathbf{w}') + s \geq m_{i'_1,\dots,i'_{k'}} + s \geq m_{i_1,\dots,i_l,i'_{l+1}} + s \geq p' + s \geq p + 2c + 2$$

and

$$m(\mathbf{w}) + s \geq m_{i_1,\dots,i_k} + s \geq m_{i_1,\dots,i_l,i_{l+1}} + s \geq p + s \geq p + 2c + 2.$$

This completes the proof of case (b). Therefore, the proof of Proposition 1 is complete. □

**5. Upper bounds in the pivotal case.** In this section we will prove in detail upper bounds for the first and second moments of $|Q_n|$. To do this, we will recall the approach of Kesten [4] to lay the groundwork that allows us to establish certain "horseshoe" estimates that we describe below. Let $B_1 \subset B(n)$ be a box centered at $\mathbf{x}$ near the right boundary of $B(n)$ such that the right boundary of $B_1$ lies on the right boundary of $B(n)$, and let $B_2 \subset B(n)$ be a box containing $B_1$ such that the right edge of $B_1$ is centered in the right edge of $B_2$. Thus, $B_2 \setminus B_1$ is a semi-annular region that we call a horseshoe. To estimate $E(|Q_n|)$, we bound the $P(\mathcal{Q}(\mathbf{x}, n))$ by the product of probabilities of two subevents of $\mathcal{Q}(\mathbf{x}, n)$, namely, (i) there exists a four-arm path from $\mathbf{x}$ to $\partial B_1$, and (ii) there exists a three-arm crossing of the horseshoe. The probability of the latter event will be handled by Lemma 5. To organize the sizes of the larger boxes $B_2$ that fit inside $B(n)$, we introduce a partition of the box $B(n)$ that is dual to the original partition of concentric annuli introduced in Section 3. For the second moment, we must estimate $P(\mathcal{Q}(\mathbf{x}, n) \cap \mathcal{Q}(\mathbf{y}, n))$. We employ the same "near to" definition employed in Section 4. When $\mathbf{y}\widetilde{N}\mathbf{x}$, so that $\mathbf{x}$ and $\mathbf{y}$ are isolated root vertices, we determine first whether these vertices are separated sufficiently to give rise to one or two horseshoes. The boxes and horseshoes we construct for our probability estimates will remain disjoint. We then utilize independence of events and Lemma 5 applied to each horseshoe that appears in our construction.

From this point of view, our method for the pivotal case may be termed the method of disjoint horseshoes. However, if $\mathbf{y}N\mathbf{x}$, then it does not suffice to simply apply a disjoint boxes argument combined with Lemma 2, because this leads to a divergent sum in our dyadic summation method. Thus we need another result, namely, Lemma 7, that is proved in the Appendix.

Let $B_1 = B_1(2^\rho) \subset B(n)$ be a fixed box of radius $2^\rho$ and for each $\nu \geq \rho$ such that $\nu - \rho$ is an integer, let $B_2 = B_2(2^\nu) \subset B(n)$ be a box of radius $2^\nu$ containing $B_1$ such that the right edge of $\partial B_1$ is centered in the right edge of $\partial B_2$. Denote by $H := H(\rho,\nu) := B_2(2^\nu) \setminus B_1(2^\rho)$ the corresponding horseshoe. Consider $\partial H$ with the right edges in common with the right edge of $\partial B_2$ removed. The resulting set of vertices consists of two concentric semi-rings of vertices in $\partial H$. The smaller semi-ring we denote by $\partial_1 H$ and call the inner horseshoe boundary and the larger semi-ring we denote by $\partial_2 H$ and call the outer horseshoe boundary. Define the event

$$\begin{aligned}\mathcal{J}(\rho,\nu) := {} & \text{there exists an open path } r_1 \text{ in } H = H(\rho,\nu) \text{ that}\\ & \text{connects } \partial_2 H \text{ to } \partial_1 H \text{ and there exist two disjoint}\\ & \text{closed paths } r_2 \text{ and } r_4 \text{ in } H(\rho,\nu) \text{ that connect}\\ & \partial_2 H \text{ to } \partial_1 H;\ r_4 \text{ is oriented counterclockwise and}\\ & r_2 \text{ clockwise from } r_1 \text{ as viewed from } \partial_2 H.\end{aligned} \tag{5.1}$$

LEMMA 5. *Define the event that there is a three-arm crossing of the horseshoe $H(\rho,\nu)$ in $B(n)$ with inner radius $2^\rho$ and outer radius $2^\nu$ by (5.1). Then there is a function $\varepsilon(u) \to 0$ as $u \to \infty$ and constant $C$ such that $P(\mathcal{J}(\rho,\nu)) \leq C 2^{\rho(2+\varepsilon(\rho))}/2^{\nu(2+\varepsilon(\nu))}$.*

PROOF. The first main step is to establish (5.2). To do this, we have to recall the proof of Kesten's [4] Lemma 4. Since Kesten's connection arguments will continue to play a role in our proof of Lemma 7, we repeat the main outlines of these arguments here for the sake of completeness. For any box $B = B(\mathbf{x}, r)$, we define the $i$th side, $i = 1,2,3,4$, as the part of the boundary of $B$, that is, respectively, on the left, bottom, right or top of $B$. Define disjoint filled squares $\beta_i = \beta_i(\rho)$, $i = 1,2,4$, that lie outside but adjacent to the sides of $B_1(2^\rho)$, where the index $i$ refers to the $i$th side, so that the squares are listed in counterclockwise order around the boundary of $B(\mathbf{y}, 2^\rho)$. Here and in the sequel a square will be synonymous with a box $B(\mathbf{x}, r)$ for some center $\mathbf{x}$ and radius $r$. We assume that the squares $\beta_i$ are of radius $2^{\rho-3}$ with spacing $7(2^{\rho-3})$ on either side. See Figure 2. Define the event $\mathcal{H}(\rho,\nu)$ as the event $\mathcal{J}(\rho,\nu)$ with the additional requirements that the path $r_1$ $h$-tunnels through $\beta_1$ and the paths $r_2$ and $r_4$ $v$-tunnel through $\beta_2$ and $\beta_4$, respectively, and further, there is a vertical open crossing of $\beta_1$ and there are horizontal closed crossings of $\beta_2$ and $\beta_4$. We will show

$$P(\mathcal{J}(\rho,\nu)) \leq C P(\mathcal{H}(\rho,\nu)). \tag{5.2}$$

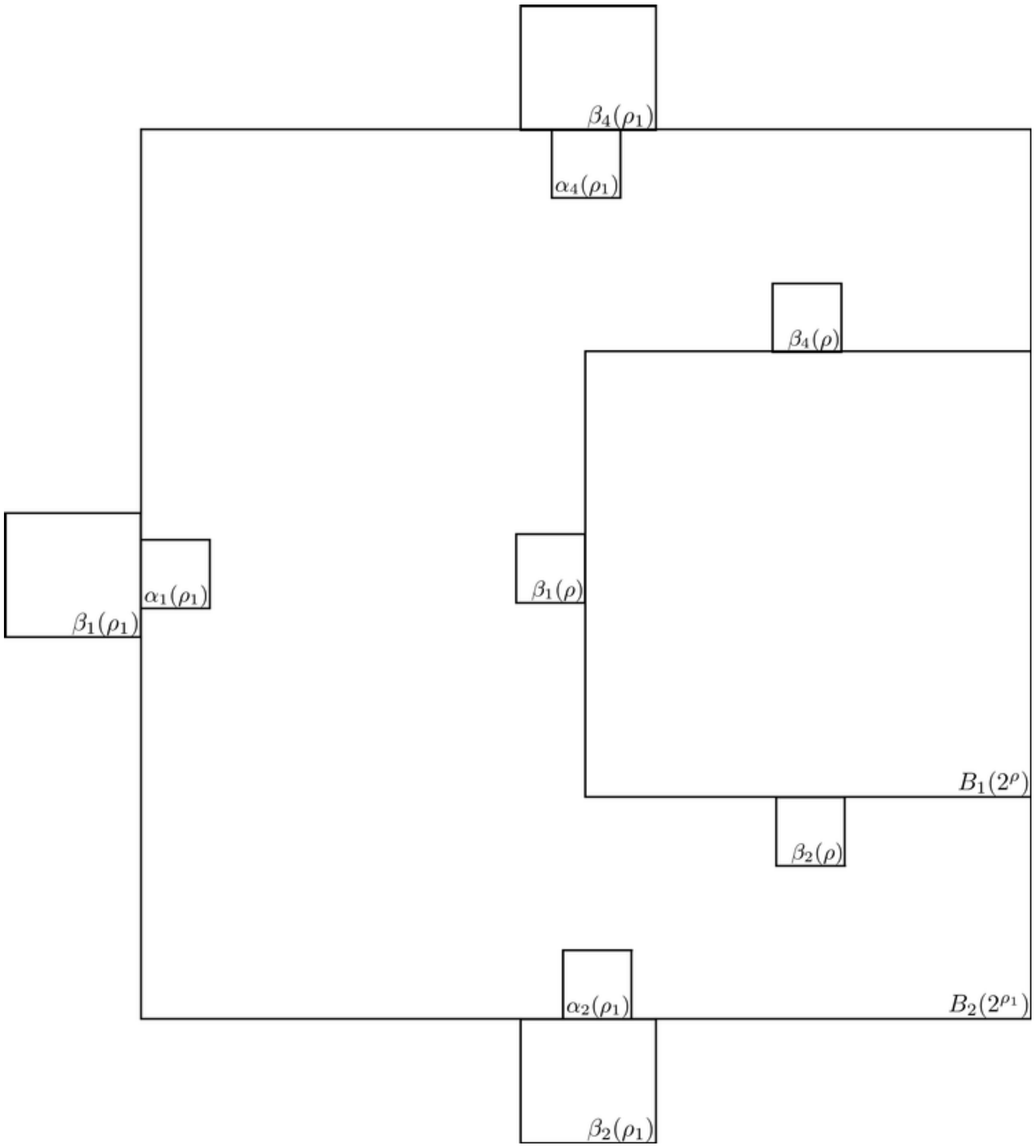

FIG. 2. *Arrangement of the connection boxes in the proof of Lemma 5. Here $\rho_1 = \rho + 1$.*

Define $H(\rho_1, \nu) := B_2(2^\nu) \setminus B_2(2^{\rho_1})$ for any $\rho \le \rho_1 \le \nu$, where by our definition above, $B_2(2^{\rho_1}) = B_1(2^\rho)$ for $\rho_1 = \rho$. We now take $\rho_1 = \nu - k$ and so view a nested sequence of boxes $B_2(2^{\nu-k})$, $k \ge 1$, each in a similar relationship to the box $B_1(2^\rho)$ as the original box $B_2(2^\nu)$. Introduce disjoint squares $\alpha_i = \alpha_i(\nu - k)$, $i = 1, 2, 4$, of radii $2^{\nu-k-3}$ that lie now inside but adjacent and centered to the $i$th sides of $B_2(2^{\nu-k})$, $k = 0, \dots, \nu - \rho - 1$. Likewise, by similarity to the squares $\beta_i(\rho)$ on the outside of $B_1(2^\rho)$, introduce corresponding squares $\beta_i(\nu - k - 1)$ of radii $2^{\nu-k-4}$ on the outside of $B_2(2^{\nu-k-1})$. First note for the case $k = 0$, that, by the existence of vertical open crossings of the squares $\alpha_1(\nu)$ and $\beta_1(\nu - 1)$ and horizontal closed crossings of the squares $\alpha_i(\nu)$ and $\beta_i(\nu-1)$, $i=2,4$, and by the existence of appropriate connecting paths that $h$-tunnel through both $\alpha_1(\nu)$ and $\beta_1(\nu-1)$ and that $v$-tunnel through $\alpha_i(\nu)$ and $\beta_i(\nu-1)$ for each $i=2,4$, and by FKG, there exists a constant $c_1$ such that

$$P(\mathcal{H}(\nu-1,\nu)) \geq c_1.$$

Now we iterate this argument with $k\geq 1$, while keeping track of the probability of connecting paths from one step to the next. Indeed, we replace in the above argument the squares $\alpha_i(\nu)$ and $\beta_i(\nu-1)$ by the squares $\alpha_1(\nu-k)$ and $\beta_1(\nu-k-1)$, and only require, besides the horizontal closed crossings and vertical open crossings, the existence of connecting paths that, as appropriate, either $h$-tunnel or $v$-tunnel through all three of $\beta_i(\nu-k)$ and $\alpha_i(\nu-k)$, and $\beta_i(\nu-k-1)$, to show by induction that there exists a constant $c_2$ such that

$$P(\mathcal{H}(\nu-k-1,\nu)) \geq c_1 c_2^{-k} \qquad \text{all } k=0,\dots,\nu-\rho-1. \tag{5.3}$$

We may assume that $\nu>\rho+2$, so we now do so. Define the event $\mathcal{J}(\rho+2,\nu)$ by replacing the horseshoe $H(\rho,\nu)$ in (5.1) by $H(\rho+2,\nu)$, so that obviously $\mathcal{J}(\rho,\nu)\subset\mathcal{J}(\rho+2,\nu)$. Consider the event $\mathcal{K}(\rho+2,\nu)$ that the paths $r_i$, $i=1,2,4$, defining $\mathcal{J}(\rho+2,\nu)$ can be chosen such that each has a certain fence around it at the location that it meets the inner horseshoe boundary $\partial_1 H(\rho+2,\nu)$; see [4], page 134, for the precise definition of the fence. Kesten shows, by adroit application of the FKG inequality (see [4], Lemma 3), that each fence, in turn, will allow an extension of the chosen path $r_i$ into $H(\rho,\rho+2)$ by means of a certain corridor it will travel through, with the result that there is only a multiplicative constant cost in probability that the path will $h$-tunnel or $v$-tunnel, as appropriate, through the corresponding square $\beta_i$:

$$P(\mathcal{K}(\rho+2,\nu)) \leq C_f P(\mathcal{H}(\rho,\nu)), \tag{5.4}$$

where $C_f$ depends on the parameter of the fence.

On the exceptional set, where one of the paths $r_i$ cannot be chosen to have such a fence, one obtains, following Kesten [4], page 131, a bound

$$P(\mathcal{J}(\rho+2,\nu)\setminus\mathcal{K}(\rho+2,\nu)) \leq \delta P(\mathcal{J}(\rho+3,\nu)). \tag{5.5}$$

The parameter $\delta$ can be made as small as desired by adjusting the parameter of the fence (see [4], Lemma 2). Therefore, by (5.1) and (5.5), one obtains

$$P(\mathcal{J}(\rho,\nu)) \leq P(\mathcal{J}(\rho+2,\nu)) \leq P(\mathcal{K}(\rho+2,\nu)) + \delta P(\mathcal{J}(\rho+3,\nu)). \tag{5.6}$$

By iteration of (5.6) and by applying (5.4) and (5.3) at the end, one obtains, just as in [4], page 131, that

$$P(\mathcal{J}(\rho,\nu)) \leq \sum_{t\geq 0}\delta^t P(\mathcal{K}(\rho+3t+2,\nu)) + C\delta^{(\nu-\rho)/3}$$

$$(5.7)\qquad \begin{aligned} &\le \sum_{t\ge 0} C_f \delta^t P(\mathcal{H}(\rho+3t,\nu)) + C(\delta c_2^3)^{(\nu-\rho)/3} P(\mathcal{H}(\rho,\nu)) \\ &\le P(\mathcal{H}(\rho,\nu))\left(\sum_{t\ge 0} C_f c_1^{-1}(\delta c_2^3)^t + C(\delta c_2^3)^{(\nu-\rho)/3}\right). \end{aligned}$$

Since $\delta$ is arbitrary, by (5.7), the desired estimate (5.2) follows.

We continue the proof of the lemma. Let $\mathbf{y}'$ be the center vertex of the right side of $B_1 := B_1(2^\rho)$. Recall that $\partial_1 H$ denotes the inner horseshoe boundary of the horseshoe $H(\rho,\nu)$. Let the squares $\alpha_i(\rho)$, $i=1,2,4$, as defined above lie inside the boundary of $B_2(2^\rho) = B_1(2^\rho)$. Define the events

$$(5.8)\qquad \begin{aligned} \mathcal{E}(\rho) := &\text{ there exists an open path } r_1 \text{ in } B_1(2^\rho) \text{ from } \mathbf{y}' \text{ to} \\ &\partial_1 H \text{ and there exist two disjoint closed paths} \\ &r_2 \text{ and } r_4 \text{ in } B_1(2^\rho) \text{ from } \mathbf{y}' \text{ to } \partial_1 H; \\ &r_2 \text{ is oriented counterclockwise and } r_4 \\ &\text{clockwise from } r_1 \text{ as viewed from the vertex } \mathbf{y}' \end{aligned}$$

and

$$(5.9)\qquad \begin{aligned} \mathcal{D}(\rho) := &\ \mathcal{E}(\rho) \text{ occurs, the path } r_1\ h\text{-tunnels through } \alpha_1(\rho), \\ &\text{and, for each } i=2,4, \text{ the paths } r_i\ v\text{-tunnel through } \alpha_i(\rho). \\ &\text{Further, there exists a vertical open crossing of } \alpha_1(\rho) \\ &\text{and horizontal closed crossings of } \alpha_2(\rho) \text{ and } \alpha_4(\rho). \end{aligned}$$

By Kesten's arguments again, $P(\mathcal{E}(\rho)) \le CP(\mathcal{D}(\rho))$ and $P(\mathcal{H}(\rho,\nu))P(\mathcal{D}(\rho)) \le CP(\mathcal{E}(\nu))$. Therefore, by (5.2) and these two inequalities, $P(\mathcal{J}(\rho,\nu))P(\mathcal{E}(\rho)) \le CP(\mathcal{H}(\rho,\nu))P(\mathcal{D}(\rho)) \le CP(\mathcal{E}(\nu))$. Therefore,

$$(5.10)\qquad P(\mathcal{J}(\rho,\nu)) \le CP(\mathcal{E}(\nu))/P(\mathcal{E}(\rho)).$$

Finally, to complete the proof of the lemma, we recall Smirnov and Werner's semi-annulus version of Lemma 1 (Theorem 3 of [12]) as follows. Let $\mathcal{G}_\kappa(r_0,r)$ denote the event that there exist $\kappa$ disjoint crossings of the semi-annulus $A_+(r_0,r) = \{z\in\mathbf{C} : r_0 < |z| < r, \Im z > 0\}$ for the hexagonal tiling of fixed mesh 1 in $\mathbf{C}$. Then for all $\kappa \ge 1$,

$$(5.11)\qquad P(\mathcal{G}_\kappa(r_0,r)) = r^{-\kappa(\kappa+1)/6+o(1)} \qquad \text{as } r\to\infty.$$

Therefore, by (5.11) with $\kappa=3$, we have

$$(5.12)\qquad P(\mathcal{E}(\rho)) = 2^{-\rho(2+o(1))} \qquad \text{as } \rho\to\infty.$$

Hence, by (5.10) and (5.12), the proof of the lemma is complete. □

For our proof of Lemma 7, we will also need the following result that is a restatement of Kesten's [4] Lemma 5. Let $B(l)$ be a box centered at the origin with radius $l\ge 2$, and let $B(\mathbf{x},m) \subset B(l/2)$. Define disjoint filled squares $\beta_i$, $i = 1, 2, 3, 4$, that lie outside but adjacent to the sides of $B(\mathbf{x}, m)$, where the index $i$ refers to the $i$th side. We take the squares to have radii $m/8$ and to be centered in the sides of $B(\mathbf{x}, m)$. Let $\mathcal{U}_4(\mathbf{x}, m; l)$ be as defined in (2.5). Let $\mathcal{V}_4(\mathbf{x}, m; l)$ be defined by (2.5) with the following additional requirements: the open paths $r_1$ and $r_3$ that exist from $\partial B(l)$ to $\partial B(\mathbf{x}, m)$ will $h$-tunnel through $\beta_1$ and $\beta_3$, respectively, on their ways to $\partial B(\mathbf{x}, m)$, and, likewise, the closed paths $r_2$ and $r_4$ will $v$-tunnel through $\beta_2$ and $\beta_4$, respectively, on their ways to $\partial B(\mathbf{x}, m)$, and, further, there exist vertical open crossings of $\beta_1$ and $\beta_3$ and horizontal closed crossings of $\beta_2$ and $\beta_4$.

LEMMA 6. *There is a constant $C$ such that*

$$P(\mathcal{U}_4(\mathbf{x}, m; l)) \leq CP(\mathcal{V}_4(\mathbf{x}, m; l)).$$

5.1. *Expectation bound for pivotal sites.* We are now ready to estimate $E(|Q_n|)$. We will refine the partition of the box $B(n)$ defined by the concentric annuli $A_j$ of (3.1) by cutting these annuli transversally. Define an increasing sequence of regions $B^*(j^*)$, $j^* \geq 0$, each lying inside $B(n)$ by

$$B^*(j^*) := \{(x_1, x_2) \in B(n) : \min\{|x_1|, |x_2|\} \leq (1 - 2^{-j^*-1})n\}. \tag{5.13}$$

The set $B^*(j^*)$ is the box $B(n)$ with squares of diameter $2^{-j^*-1}n$ removed from each of its corners. We define the dual sets to the annuli $A_j$ by taking the successive differences of the sets $B^*(j^*)$:

$$\begin{aligned} &A_0^* := B^*(0), \\ &A_{j^*}^* := B^*(j^*) \setminus B^*(j^* - 1), \qquad j^* \geq 1. \end{aligned} \tag{5.14}$$

Thus, for $j^* > 0$, $A_{j^*}^*$ consists of four "L"-shaped regions. For each such region, the "L" cuts off a square in the corresponding corner of $B(n)$. The collection $\{A_{j^*}^*,\ j^* \geq 0\}$ is a partition of $B(n)$. Moreover, the following properties hold:

$$A_j \cap A_{j^*}^* \begin{cases} = \varnothing, & j^* > j, \\ \text{is a union of eight rectangles}, & \text{if } 0 < j^* < j, \\ \text{is a union of four corner squares}, & \text{if } 0 < j^* = j, \\ \text{is a union of four rectangles}, & \text{if } j^* = 0 \text{ and } j > 0, \\ \text{one central square}, & \text{if } j^* = j = 0. \end{cases} \tag{5.15}$$

Note, by (3.1), (5.13) and (5.14), that we have the estimate

$$|A_j \cap A_{j^*}^*| \leq C2^{-j-j^*}n^2 \qquad \text{all } 0 \leq j^* \leq j,\ n \geq 1 \tag{5.16}$$

for some constant $C$. Let $\{A_j,\ 0 \leq j \leq j_0\}$ be the partition of $B(n)$ defined by (3.1). Thus, by (5.15), the collection $\{A_j \cap A_{j^*}^*,\ 0 \leq j^* \leq j \leq j_0\}$ comprises a joint partition of $B(n)$. Hence, we can write

$$E(|Q_n|) = \sum_{\mathbf{x}\in B(n)} P(\mathcal{Q}(\mathbf{x},n)) = \sum_{j=0}^{j_0}\sum_{j^*=0}^{j}\sum_{\mathbf{x}\in A_j\cap A^*_{j^*}} P(\mathcal{Q}(\mathbf{x},n)). \tag{5.17}$$

Now for any $0 \le j^* \le j$, we consider $\mathbf{x} \in A_j \cap A^*_{j^*}$. Choose real numbers $\rho = \rho(j,n)$ and $\nu = \nu(j^*,n)$ such that $2^\rho \asymp 2^{-j}n$ and $2^\nu \asymp 2^{-j^*}n$ and $\nu - \rho$ is integer. Here $f \asymp g$ over a range of arguments for the functions $f$ and $g$ means that there exists a constant $C > 0$ such that $(1/C)g \le f \le Cg$ over this range. We choose $B_1(2^\rho)$ to have center $\mathbf{x}$ and make the definition of $\rho$ such that $B_1(2^\rho) \subset B(n)$, but also such that $\partial B_1(2^\rho) \subset \partial B(n)$. This is possible since the box $B(\mathbf{x}, 2^{-j-2}n)$ lies interior to $B(n)$ by construction, so now we expand the radius of this box such that its boundary just meets that of $B(n)$. Notice therefore that while $\rho$ is not independent of $\mathbf{x}$, the value of $2^\rho$ only varies by a constant factor with $\mathbf{x}$. We also construct a box $B_2(2^\nu) \subset B(n)$ containing $B_1(2^\rho)$, as in the context of Lemma 5, such that $B_2(2^\nu)$ and $B_1(2^\rho)$ share boundary points along the side of $B(n)$ corresponding to the side of the annulus $A_j$ that $\mathbf{x}$ belongs to. This is possible by our construction of the dual $A^*_{j^*}$. Thus, by the definition (5.1) of $J(\rho,\nu)$ and the definition of the four-arm path (2.5), and by independence, we have that

$$P(\mathcal{Q}(\mathbf{x},n)) \le P(\mathcal{U}_4(\mathbf{0},2^\rho))P(\mathcal{J}(\rho,\nu)). \tag{5.18}$$

Let $\varepsilon > 0$. By Lemma 2, there exists a constant $C_{\varepsilon,1}$ such that

$$P(\mathcal{U}_4(\mathbf{0},r)) \le C_{\varepsilon,1} r^{-5/4+\varepsilon} \qquad \text{all } r \ge 1.$$

Similarly, by Lemma 5, there exists a constant $C_{\varepsilon,2}$ such that $P(\mathcal{J}(\rho,\nu)) \le C_{\varepsilon,2}2^{(\rho-\nu)(2-\varepsilon)}$, all $\rho \le \nu$. Therefore, by these considerations with $r = C2^{-j}n \ge C/2$ and with $2^{-j}n$ and $2^{-j^*}n$ in place of $2^\rho$ and $2^\nu$, respectively, we have, by (5.16) and (5.18), that there exists a constant $C_\varepsilon$ such that

$$\begin{aligned} E(|Q_n|) &\le C_\varepsilon \sum_{j=0}^{j_0}\sum_{j^*=0}^{j} 2^{-j-j^*}n^2(2^{-j}n)^{-5/4+\varepsilon}2^{(-j+j^*)(2-\varepsilon)} \\ &\le C_\varepsilon n^{3/4+\varepsilon}\sum_{j=0}^{\infty} 2^{-3j/4} \le C_\varepsilon n^{3/4+\varepsilon}. \end{aligned} \tag{5.19}$$

This concludes the case $\tau = 1$ of item 3 of Theorem 1.

5.2. *Second moment for pivotal sites.* In the case of second and higher moments we will have to consider the condition that a given root vertex $\mathbf{x}_{e_i}$ is not isolated (so that $|V_i| \ge 1$; see Section 4). To handle the need for an extra convergence factor in our dyadic summation method, in this case we introduce the following lemma. We will use this lemma, in particular, to estimate the probability of the event $\mathcal{Q}(\mathbf{x},n)\cap\mathcal{Q}(\mathbf{y},n)$ in case $\mathbf{y}N\mathbf{x}$ for the second moment estimate below. Let $R$ be a filled-in rectangle of vertices with sides parallel to the coordinate axes. For any vertex $\mathbf{w}$ contained in the interior of $R$, denote the event

$$\mathcal{U}_4(\mathbf{w};R) := \exists \text{ a four-arm path in } R \text{ from } \mathbf{w} \text{ to } \partial R. \tag{5.20}$$

This is simply an extension of the definition $\mathcal{U}_4(\mathbf{w},n)$ in (2.5) with $R$ in place of $B(n)$. Recall also definition (2.8).

LEMMA 7. *Let $R = R(\mathbf{x}')$ be a rectangle centered at $\mathbf{x}'$ with its shortest half-side of length $l \geq 1$ and longest half-side of length $L \geq 1$ such that $1 \leq L/l \leq 2$. Let $R$ contain a vertex $\mathbf{x}$ such that $\|\mathbf{x}-\mathbf{x}'\| \leq l/2$. Suppose further that $R$ contains a collection of disjoint boxes $B_i := B(\mathbf{y}_i, 2^{\lambda_i})$, $i=1,\ldots,v$, which also have the property that for each $i=1,\ldots,v$, $\mathbf{x}\notin B_i$. Then there exist constants $C$, $d$ and $c_1$, depending only on $v$, such that*

$$P\left(\mathcal{U}_4(\mathbf{x};R)\cap\left(\bigcap_{i=1}^{v}\mathcal{U}_4(\mathbf{y}_i;R)\right)\right) \leq CP(\mathcal{T}_4(\mathbf{0},l/d))\prod_{i=1}^{v}P(\mathcal{U}_4(\mathbf{0},2^{\lambda_i-c_1})).$$

We prove Lemma 7 in the Appendix.

We are now ready to estimate the second moment of $|Q_n|$. As in the estimation of the first moment, we use the partition $\{A_j\cap A^*_{j^*},\ 0\leq j^*\leq j\leq j_0\}$ of $B(n)$. Denote $p_n(\mathbf{x},\mathbf{y}) := P(\mathcal{Q}(\mathbf{x},n)\cap\mathcal{Q}(\mathbf{y},n))$. Hence, as in (4.1)–(4.2), it suffices to estimate

$$\Sigma_0 := \sum_{j=0}^{j_0}\sum_{j^*=0}^{j}\sum_{\mathbf{x}\in A_j\cap A^*_{j^*}}\sum_{k=j}^{j_0}\sum_{k^*=0}^{k}\sum_{\mathbf{y}\in A_k\cap A^*_{k^*}} p_n(\mathbf{x},\mathbf{y}). \tag{5.21}$$

Recall that, as mentioned at the beginning of Section 4, we may assume that the vertices $\mathbf{x}$ and $\mathbf{y}$ are distinct in (5.21). We now define a diagonal subsum of the sum (5.21) according to the condition $\mathbf{y}\in B(\mathbf{x},2^{-j-2c}n)$ that is $\mathbf{y}N\mathbf{x}$. For this case, we recall by (3.11) that $\bigcup_{m=j+2c}^{j_0} a_m(\mathbf{x}) = B(\mathbf{x},2^{-j-2c}n)$, where $c$ is defined in Section 4 by (4.4). Therefore, we write this diagonal sum as

$$I := \sum_{j=0}^{j_0}\sum_{j^*=0}^{j}\sum_{\mathbf{x}\in A_j\cap A^*_{j^*}}\sum_{m=j+2c}^{j_0}\sum_{\mathbf{y}\in a_m(\mathbf{x})} p_n(\mathbf{x},\mathbf{y}). \tag{5.22}$$

Let $\varepsilon>0$. We estimate $I$. We use that if $\mathbf{y}\in a_m(\mathbf{x})$, then (3.13) holds. If $m\geq j+2c+4$, we apply Lemma 7 with $v=1$ and $\mathbf{x}'=\mathbf{x}$ for a square $R\subset B(n)$ with half-side $l\asymp 2^{-j}n$ such that $B(\mathbf{y},2^{-m-2}n)\subset R$. Note that, indeed, $\mathbf{x}\notin B(\mathbf{y},2^{-m-2}n)$ for $\mathbf{y}\in a_m$ so that the hypothesis of Lemma 7 is satisfied with $2^{\lambda_1} = 2^{-m-2}n$. If instead $j + 2c \le m < j + 2c + 4$, then we can still define the square $R$ for the same asymptotic size of $l$ but such that now the box $B(\mathbf{y}, 2^{-m-2}n) \subset B(n)$ is disjoint from $R$. In this latter case we simply apply independence of events. Finally, we find a box $B_1$ in $B(n)$ of radius $2^\rho \asymp 2^{-j}n$, one of whose edges lies in $\partial B(n)$ and that contains both $R$ and the box $B(\mathbf{y}, 2^{\lambda_1})$. For $\mathbf{x} \in A_j \cap A^*_{j^*}$, we construct a box $B_2 \subset B(n)$ of radius $2^\nu \asymp 2^{-j^*}n$ such that the horseshoe pair $(B_1, B_2)$ conforms to the context of Lemma 5. Therefore, since on the event $\mathcal{Q}(\mathbf{x}, n) \cap \mathcal{Q}(\mathbf{y}, n)$ there must be a three-arm crossing of the horseshoe, by Lemma 5, independence and Lemma 7, and by Lemma 2 applied to both $P(\mathcal{U}_4(\mathbf{0}, 2^{\lambda_1 - c_1}))$ and $P(\mathcal{U}_4(\mathbf{0}, l/d))$ for $l \asymp 2^\rho$, we have that

$$p_n(\mathbf{x}, \mathbf{y}) \le C_\varepsilon 2^{(\rho+\lambda_1)(-5/4+\varepsilon)+(\rho-\nu)(2-\varepsilon)} \le C_\varepsilon n^{-5/2+2\varepsilon} 2^{2j^*+(-3j+5m)/4}. \tag{5.23}$$

Therefore, by (3.12), (5.16), (5.22) and (5.23), we have that

$$\begin{aligned} I &\le C_\varepsilon n^{-5/2+2\varepsilon} \sum_{j=0}^{\infty} \sum_{j^*=0}^{j} |A_j \cap A^*_{j^*}| \sum_{m=j}^{\infty} |a_m| 2^{2j^*+(-3j+5m)/4} \\ &\le C_\varepsilon n^{3/2+2\varepsilon} \sum_{j=0}^{\infty} \sum_{j^*=0}^{j} \sum_{m=j}^{\infty} 2^{j^*+(-7j-3m)/4} \le C_\varepsilon n^{3/2+2\varepsilon} \sum_{j=0}^{\infty} 2^{-3j/2}. \end{aligned} \tag{5.24}$$

We now turn to the remaining sum $II := \Sigma_0 - I$. The vertices $\mathbf{x}$ and $\mathbf{y}$ left to consider in this sum are isolated root vertices, so that $\mathbf{y}\widetilde{N}\mathbf{x}$. We would like to construct a horseshoe along the side of $B(n)$ for each vertex in this pair of vertices. But boxes centered at these vertices, defined by the condition that each box just comes to the side of $B(n)$, may overlap. To treat this case, we define an (asymmetric) horseshoe relationship as follows. We write $\mathbf{y}K\mathbf{x}$ if $\mathbf{y} \in B(\mathbf{x}, 2^{-j+6}n)$. If $\mathbf{y}K\mathbf{x}$, then we define one root horseshoe vertex, that is, $\mathbf{x}$, while if $\mathbf{y}\widetilde{K}\mathbf{x}$, then both $\mathbf{x}$ and $\mathbf{y}$ are defined as root horseshoe vertices. We shall refer to these cases, respectively, by the number $h$ of root horseshoe vertices, namely, $h = 1$ or $h = 2$. The horseshoe relationship provides a useful way to organize our construction of estimates.

Consider first that $h = 1$ and write $II_1$ for the sum over pairs of vertices in $II$ corresponding to this condition. Though there are few vertices $\mathbf{y}$ near $\mathbf{x}$ that satisfy both the conditions $\|\mathbf{x} - \mathbf{y}\| > 2^{-j-2c}n$ and $h = 1$, the fact that $\mathbf{y}$ may lie near the boundary of $B(n)$ requires us to construct a horseshoe at $\mathbf{y}$ if $k$ is much larger than $j$. So we write $II_1 = II_1a + II_1b$, where the sums $II_1a$ and $II_1b$ correspond, respectively, to the cases (a) $j \le k \le j + 2c + 2$, and (b) $k > j + 2c + 2$. In case (a), because $h = 1$, we can fit two disjoint boxes centered at our vertices each with radius asymptotic to $2^{-j}n$ inside a box $B_1$ that has a radius $2^\rho \asymp 2^{-j}n$ of the same asymptotic order. Yet $B_1$ also has one edge in $\partial B(n)$. Again, we find a box $B_2 \subset B(n)$ of radius $2^\nu \asymp 2^{-j^*}n$ such that the pair $(B_1, B_2)$ conforms to the context of Lemma 5. Note that the size of the set of vertices that $\mathbf{y}$ is confined to by the conditions $h=1$ and (a) is bounded by $C2^{-2j}n$. Therefore, in a similar fashion as the estimation of $I$ but now without the use of Lemma 7, we have that

$$
\begin{aligned}
II_1a &\le C_\varepsilon n^{-5/2+2\varepsilon} \sum_{j=0}^{\infty}\sum_{j^*=0}^{j} |A_j \cap A^*_{j^*}| 2^{-2j} n^2 2^{2\rho(-5/4+\varepsilon)+(\rho-\nu)(2-\varepsilon)} \\
&\le C_\varepsilon n^{3/2+2\varepsilon} \sum_{j=0}^{\infty}\sum_{j^*=0}^{j} 2^{j^*-5j/2}.
\end{aligned} \tag{5.25}
$$

To estimate $II_1b$, we first construct a pair of boxes $(B_1, B_2)$ as in the context of Lemma 5 with the parameter $\sigma$ playing the role of $\rho$ as follows. We find $B_1 := B(\mathbf{y}, 2^\sigma) \subset B(n)$ with $2^\sigma \asymp 2^{-k}n$ such that $B_1$ has one side in the boundary of $B(n)$. We take $B_2$ accordingly by defining its radius $2^\nu \asymp 2^{-j}n$ such that $B_2$ is disjoint from $B(\mathbf{x}, 2^\rho)$ for $2^\rho = 2^{-j-2c-1}n$. That $B_1$ and $B_2$ will exist follows by (b) and the assumption that the vertices are isolated roots. Now since $h=1$, we can also find another inner horseshoe box $\widetilde{B}_1$ with radius $2^{\rho_1} \asymp 2^{-j}n$ that now contains both $B(\mathbf{x}, 2^{-j-2c-1}n)$ and $B_2$. We pair the box $\widetilde{B}_1$ with an outer horseshoe box $\widetilde{B}_2$ with radius $2^{\nu_1} \asymp 2^{-j^*}n$. Thus, we have the horseshoe formed by the pair $(B_1, B_2)$ nested inside the horseshoe formed by $(\widetilde{B}_1, \widetilde{B}_2)$. Hence, by independence and Lemmas 2 and 5, we find that

$$
\begin{aligned}
p_n(\mathbf{x}, \mathbf{y}) &\le P(\mathcal{U}_4(\mathbf{0}, 2^\rho)) P(\mathcal{U}_4(\mathbf{0}, 2^\sigma)) P(\mathcal{J}(\sigma, \nu)) P(\mathcal{J}(\rho_1, \nu_1)) \\
&\le C_\varepsilon n^{-5/2+2\varepsilon} 2^{2j^*+(5j-3k)/4}.
\end{aligned} \tag{5.26}
$$

Therefore, since there are only on the order of $2^{-j-k}n^2$ vertices $\mathbf{y}$ accounted for when $\mathbf{y} \in A_k$ in the sum $II_1b$, we find by (5.26) that

$$
II_1b \le C_\varepsilon n^{3/2+2\varepsilon} \sum_{j=0}^{\infty}\sum_{j^*=0}^{j}\sum_{k=j}^{\infty} 2^{j^*+(-3j-7k)/4}. \tag{5.27}
$$

Consider next that $h=2$. Write $II_2$ for the sum over pairs of vertices falling under $II$ that correspond to this condition. In this way $II = II_1 + II_2$. As before, we define $j$, $j^*$, $k$ and $k^*$ by the inclusions $\mathbf{x} \in A_j \cap A^*_{j^*}$ and $\mathbf{y} \in A_k \cap A^*_{k^*}$. We consider two cases: either (a) $|k^* - j^*| \le 1$ or (b) $|k^* - j^*| > 1$. Accordingly, we will break up the sum $II_2$ into the sum $II_2 = II_2a + II_2b$ with summands corresponding, respectively, to these cases.

We first work with the case (a). Since $k^*$ is almost equal to $j^*$, we shall in effect lengthen the set $A_k \cap A^*_{j^*}$ in the long directions of $A_k$, and denote this lengthening by

$$
A_{k,j^*} := A_k \cap \left( \bigcup_{k^*=j^*-1}^{j^*+1} A^*_{k^*} \right). \tag{5.28}
$$

Consider one of the eight connected components $A'_{k,j^*} \subset A_{k,j^*}$ that is a rectangular section of $A_k$ on the same side of $A_k$ that $\mathbf{x}$ belongs to in $A_j$. The other components of $A_{k,j^*}$ can be handled similarly. Note that by (5.28) and (5.16), $A'_{k,j^*}$ has dimensions on order of $2^{-k}n$ by $2^{-j^*}n$. Assume, without loss of generality, that $\mathbf{x} = (x_1, x_2)$ and $\mathbf{y} = (y_1, y_2) \in A'_{k,j^*}$ both belong to the right side $A_j$ and $A_k$, respectively.

We introduce bands of vertices $b_v := b_v(\mathbf{x}, j, k, j^*)$ in $A'_{k,j^*}$, for $v \geq 1$, by

$$b_v := \{(y_1, y_2) \in A'_{k,j^*} : v2^{-j+5}n < |x_2 - y_2| \leq (v+1)2^{-j+5}n\}. \tag{5.29}$$

Here $v$ ranges up to order $v_{\max} \asymp 2^{j-j^*}$. By our construction, these bands of vertices cross $A'_{k,j^*}$ transversally. Note that by the assumption $h = 2$ that $\mathbf{y} \in b_v$ only for some $v \geq 2$ so that $A'_{k,j^*} = \bigcup_{v=2}^{v_{\max}} b_v$. As for the sizes of the bands $b_v$, we have by the definition of $A_k$ and (5.29) that, independent of $j^*$ and $v$,

$$|b_v| \leq C2^{-j-k}n^2. \tag{5.30}$$

Now choose $2^\rho \asymp 2^{-j}n$ such that the right edge of $B_1(\mathbf{x}) := B(\mathbf{x}, 2^\rho)$ just meets $\partial B(n)$. Also, for $\mathbf{y} \in b_v$, define $2^\sigma \asymp 2^{-k}n$ such that the right edge of $B_1(\mathbf{y}) := B(\mathbf{y}, 2^\sigma)$ just meets $\partial B(n)$. These are the inner boxes of horseshoes we will construct at each of $\mathbf{x}$ and $\mathbf{y}$. For each $v = 2, \ldots, v_{\max}$, we define an exponent $\nu$ by $2^\nu \asymp v2^{-j}n$, uniformly in $v \geq 2$, so that boxes $B_2(\mathbf{x}) \subset B(n)$ and $B_2(\mathbf{y}) \subset B(n)$, each with radius $2^\nu$, exist and are disjoint such that $(B_1(\mathbf{x}), B_2(\mathbf{x}))$ and $(B_1(\mathbf{y}), B_2(\mathbf{y}))$ each form a horseshoe pair as in the context of Lemma 5. The outer boxes remain in $B(n)$ by (5.28) and the expression for $v_{\max}$. Moreover, the outer boxes $B_2(\mathbf{x})$ and $B_2(\mathbf{y})$, while disjoint, are nested inside another box $\widetilde{B}_1$ of radius $C2^\nu$ whose right edge lies in $\partial B(n)$. Since we are in case (a), we may again pair $\widetilde{B}_1$ with an outer horseshoe box $\widetilde{B}_2$ of radius $2^{\nu_1} \asymp 2^{-j^*}n$. Therefore, by independence and by application of Lemma 5 to the horseshoe pairs, and by Lemma 2, we obtain

$$\begin{aligned} p_n(\mathbf{x}, \mathbf{y}) &\leq C_\varepsilon 2^{(\rho+\sigma)(-5/4+\varepsilon)+(\sigma+\rho-\nu-\nu_1)(2-\varepsilon)} \\ &\leq C_\varepsilon n^{-5/2+2\varepsilon} v^{-2+\varepsilon} 2^{2j^*+(5j-3k)/4}. \end{aligned} \tag{5.31}$$

Hence, by (5.16), (5.30) and (5.31), we have

$$II_2a \leq C_\varepsilon n^{3/2+2\varepsilon} \sum_{j=0}^{\infty} \sum_{j^*=0}^{j} \sum_{k=j}^{\infty} \sum_{v=1}^{\infty} v^{-2+\varepsilon} 2^{j^*+(-3j-7k)/4}. \tag{5.32}$$

Finally, we turn to the sum $II_2b$. By (b) we have $|x_2 - y_2| \geq 2^{-j^*-1}n \geq 2^{-k^*}n$ if $k^* > j^* + 1$, while $|x_2 - y_2| \geq 2^{-k^*-1}n \geq 2^{-j^*}n$ if $k^* < j^* - 1$. We can therefore define two different values of $\nu$, namely, $\nu_1$ and $\nu_2$, by $2^{\nu_1} \asymp 2^{-j^*}n$ and $2^{\nu_2} \asymp 2^{-k^*}n$ to obtain two disjoint horseshoes with outer radii $2^{\nu_1}$ and $2^{\nu_2}$. In detail, we have two pairs of boxes $(B_1, B_2)$, where each pair of boxes conforms to the context of Lemma 5, and where both larger boxes $B_2$ are disjoint and belong to $B(n)$. In one pair $B_1 = B(\mathbf{x}, 2^\rho)$ and $B_2$ has radius $2^{\nu_1}$ and in the other pair $B_1 = B(\mathbf{y}, 2^\sigma)$ and $B_2$ has radius $2^{\nu_2}$. Here $2^\rho \asymp 2^{-j}n$ and $2^\sigma \asymp 2^{-k}n$ are chosen such that the inner boxes $B_1$ lie along the boundary of $B(n)$ and are disjoint by $h = 2$. See Figure 3. Therefore, by independence and Lemmas 2 and 5, we estimate that

$$\begin{aligned} p_n(\mathbf{x}, \mathbf{y}) &\le C_\varepsilon 2^{(\rho+\sigma)(-5/4+\varepsilon)+(\rho+\sigma-\nu_1-\nu_2)(2-\varepsilon)} \\ &\le C_\varepsilon n^{-5/2+2\varepsilon} 2^{(2j^*+2k^*)} 2^{(-3j-3k)/4}. \end{aligned} \tag{5.33}$$

Hence, by (5.16) and (5.33), we have

$$II_2 b \le C_\varepsilon n^{3/2+2\varepsilon} \sum_{j=0}^{\infty} \sum_{j^*=0}^{j} \sum_{k=j}^{\infty} \sum_{k^*=0}^{k} 2^{j^*+k^*+(-7j-7k)/4}. \tag{5.34}$$

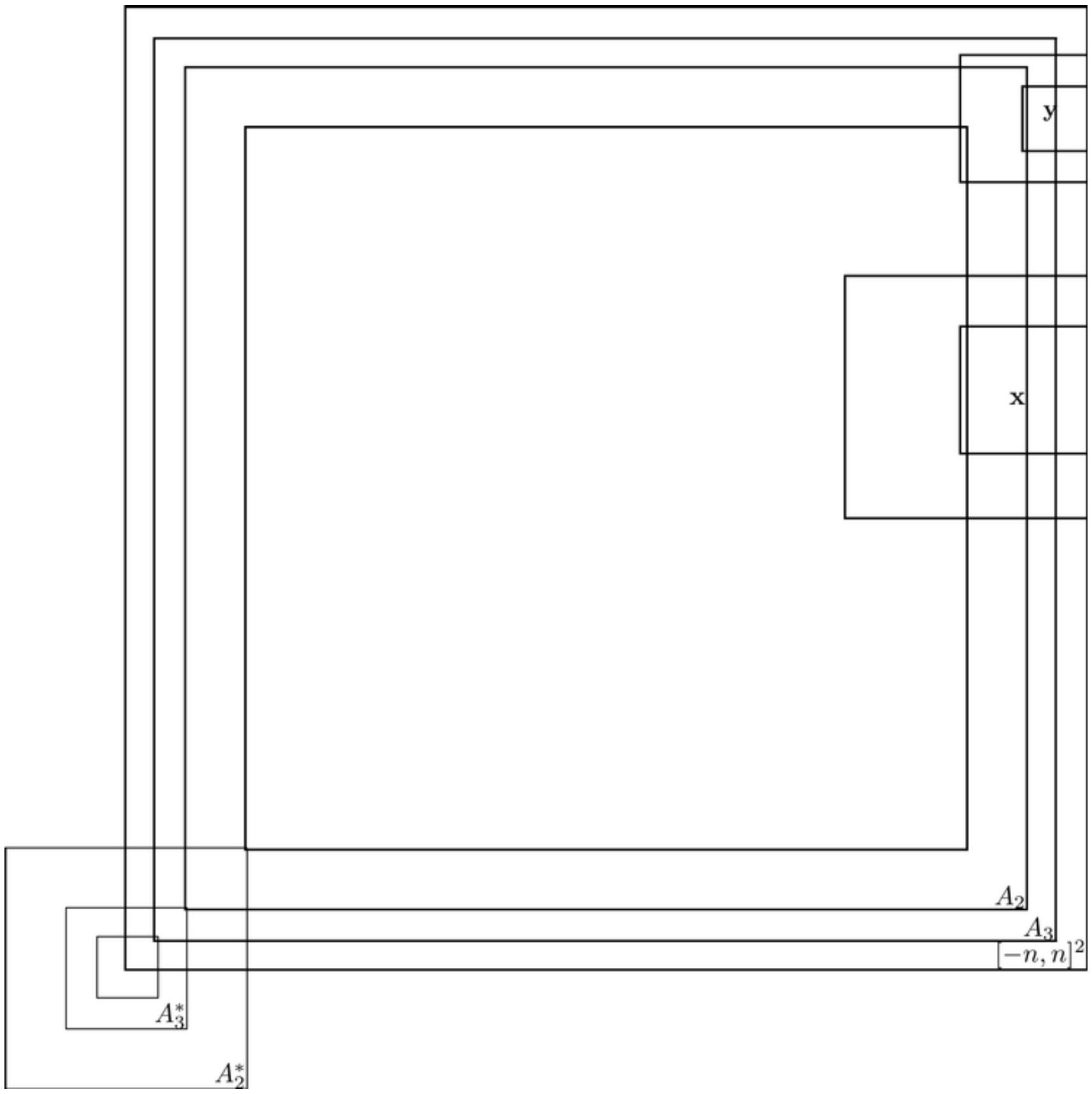

FIG. 3. *Horseshoes for $\tau = 2$ in case $|j^* - k^*| > 1$; $\mathbf{x} \in A_2 \cap A_0^*$, $\mathbf{y} \in A_3 \cap A_2^*$.*

Thus, by (5.17), (5.24), (5.25), (5.27), (5.32) and (5.34), we have proved item 3 of Theorem 1 for $\tau = 2$.

**6. Higher moments for pivotal sites.** In this section we show how to generalize the first and second moments for the number of pivotal sites shown in the previous section. We will outline the main ingredients for establishing a general $\tau$th moment by considering in some detail the case $\tau = 3$. The main issues not covered so far will be to determine (a) the sizes and numbers of horseshoes to construct, (b) the manner in which Proposition 1 is applied, and (c) the way that Lemma 7 is applied. We recall the definition $p_{n,\tau}(\mathbf{x}_1, \ldots, \mathbf{x}_\tau) := P(\bigcap_{i=1}^{\tau} \mathcal{Q}(\mathbf{x}_i, n))$ and the sum $\Sigma_0$ that we must estimate in (4.1)–(4.2). As in the previous section, we assume that each of the vertices in (4.1) belongs to the right side of its respective annulus $A_{j_i}$.

To organize our construction of estimates, we generalize the horseshoe relation $K$ on the set of root vertices $\mathbf{x}_{e_i}$, $i = 1, \ldots, r$, defined in Section 4. Write

$$\mathbf{x}_{e_k} K \mathbf{x}_{e_i} \qquad \text{if } \mathbf{x}_{e_k} \in B(\mathbf{x}_{e_i}, 2^{-j_{e_i}+6} n).$$

The constant in the exponent allows some breathing room so that, in particular, if $\mathbf{x}_{e_k} \widetilde{K} \mathbf{x}_{e_i}$, then there exists $l \asymp 2^{-j_{e_k}} n$ such that the right edge of the box $B(\mathbf{x}_{e_k}, l)$ lies on $\partial B(n)$ and such that $\|\mathbf{x}_{e_i} - \mathbf{x}_{e_k}\| \geq 4l$.

We define root horseshoe vertices among the set of root vertices by analogy with the definition of root vertices in Section 4 but now for the horseshoe relationship. We denote these root horseshoe vertices by $\mathbf{x}_{f_a}$, $a = 1, \ldots, h$, for some $h \leq r$ where $f_a = e_{i_a}$ and, in particular, $f_1 = 1$. We shall define the sets $U_a$ of root vertices chained to the root horseshoe vertices $\mathbf{x}_{f_a}$ in a way that is different from the default definition given by the method of Section 4. The reason for this is that the organization of certain probability estimates we make below is sensitive to the order of the indices in the roots that are not root horseshoe vertices. We proceed inductively as follows. First, if $\mathbf{x}_{e_2} K \mathbf{x}_{e_1}$, then $\mathbf{x}_{e_2} \in U_1$, else by our definition of root horseshoe vertices $\mathbf{x}_{e_2}$ is the root horseshoe vertex $\mathbf{x}_{f_2}$. Suppose now that $f_2 = e_2$ so indeed $\mathbf{x}_{e_2}$ is the second root horseshoe vertex. Then if $\mathbf{x}_{e_3} K \mathbf{x}_{e_1}$ but $\mathbf{x}_{e_3} \widetilde{K} \mathbf{x}_{f_2}$, we put $\mathbf{x}_{e_3} \in U_1$. This is the default arrangement that we spoke of. However, if instead $\mathbf{x}_{e_3} K \mathbf{x}_{f_2}$, then we put instead $\mathbf{x}_{e_3} \in U_2$. In general, if each $\mathbf{x}_{e_k}$, $1 \leq k \leq i$, has been designated as some root horseshoe vertex $\mathbf{x}_{f_a}$ with some $1 \leq a \leq b$ or placed as an element of $U_a$ for some $1 \leq a \leq b$, then we determine the designation or placement of $\mathbf{x}_{e_{i+1}}$ by the following rule. If $\mathbf{x}_{e_{i+1}} \widetilde{K} \mathbf{x}_{e_k}$ for all $1 \leq k \leq i$, then $e_{i+1} = f_{b+1}$, that is, we have a new root horseshoe vertex. Else we place $\mathbf{x}_{e_{i+1}}$ in the set $U_a$ of highest index $a$ such that $\mathbf{x}_{e_{i+1}}$ is related to some vertex in the current set $\{\mathbf{x}_{f_a}\} \cup U_a$. Thus, the sets $U_a$ are continuously updated, but elements may only be added and not subtracted, and they are only added at the highest possible level subject to a chain condition available at the current step. As a consequence, we obtain after construction that if $i < j$ and if $\mathbf{x}_{e_i} \in \{\mathbf{x}_{f_a}\} \cup U_a$, and if $\mathbf{x}_{e_j} \in \{\mathbf{x}_{f_b}\} \cup U_b$ for some $b < a$, then $\mathbf{x}_{e_j} \widetilde{K} \mathbf{x}_{e_i}$ since indeed $\mathbf{x}_{e_j}$ is not related to any element of $\{\mathbf{x}_{f_a}\} \cup U_a$. To see how this property is used, see the comments of Section 6.1.5 following (6.32).

6.1. *The third moment.* We organize our discussion of the case $\tau = 3$ at first according to the value of $r$. Subsequent levels of organization derive from the values of $h$ and a further parameter $t \le h$ that we shall define below.

6.1.1. $\tau = 3$, $r = 1$. We assume $r = 1$ so that $|V_1| = 2$. Thus, $\mathbf{G} = \mathbf{G}_1$ has vertices given by either (i) $\{\mathbf{x}_1, \mathbf{w}_1, \mathbf{w}_{1,1}\}$ or (ii) $\{\mathbf{x}_1, \mathbf{w}_1, \mathbf{w}_2\}$. We write $I$ to denote the sub-sum of $\Sigma_0$ that corresponds to $r = 1$. We also write $I = Ii + Iii$, where the sums $Ii$ and $Iii$ correspond, respectively, to the cases (i) and (ii).

Assume first that (i) holds. Consider now the following subcases under (i):

(a) $m_{1,1} \ge m_1 + s + 2$,
(b) $m_{1,1} < m_1 + s + 2$,

where $s$ is the constant $2c + 4$. Partition the sum $Ii = Iia + Iib$ accordingly. Consider first subcase (a). Put $R := B(\mathbf{w}_1, l)$, for $l := 2^{-m_1 - s}$. Since $\mathbf{w}_{1,1} \in a_{m_{1,1}}(\mathbf{w}_1)$, we have by (a) that $B(\mathbf{w}_{1,1}, 2^{-m_{1,1}-2}n) \subset R$, while also $\mathbf{w}_1 \notin B(\mathbf{w}_{1,1}, 2^{-m_{1,1}-2}n)$. Therefore, we may apply Lemma 7 with $v = 1$, $\mathbf{x}' = \mathbf{x} = \mathbf{w}_1$, $\mathbf{y}_1 = \mathbf{w}_{1,1}$, and $2^{\lambda_1} = 2^{-m_{1,1}-2}n$. We also apply Proposition 1 to the subgraph of $G_1$ with vertex set $\{\mathbf{x}_1, \mathbf{w}_1\}$ only. Therefore, $B(\mathbf{x}_1, l)$ and $R$ are disjoint. Hence, by these results and independence, we find that

$$p_{n,3}(\mathbf{x}_1, \mathbf{w}_1, \mathbf{w}_{1,1}) \le P(\mathcal{U}_4(\mathbf{0}, l))P(\mathcal{U}_4(\mathbf{0}, 2^{\lambda_1 - c_1}))P(\mathcal{U}_4(\mathbf{0}, l/d)). \tag{6.1}$$

Note that automatically, because $r = 1$, we have that $h = 1$ in any case. Therefore, we may construct a horseshoe with inner radius $2^{\rho} \asymp 2^{-j_1}n$ and outer radius $2^{\nu} \asymp 2^{-j_1^*}n$ for $\mathbf{x} \in A_{j_1} \cap A^*_{j_1^*}$ whose inner box contains all the boxes discussed above. Hence, by Lemma 5, we will be able to improve the estimate (6.1) by a factor $2^{(\rho-\nu)(2-\varepsilon)}$. Thus, by (6.1), Lemma 2 and this last observation, we estimate that

$$p_{n,3}(\mathbf{x}_1, \mathbf{w}_1, \mathbf{w}_{1,1}) \le C_\varepsilon n^{-15/4+3\varepsilon} 2^{2j_1^* - 2j_1 + 5(2m_1 + m_{1,1})/4}. \tag{6.2}$$

Since we apply the size estimate (3.12) for each of the indices $m = m_1$ and $m = m_{1,1}$, we obtain by (6.2) that

$$Iia \le C_\varepsilon n^{9/4+3\varepsilon} \sum_{j_1=0}^{\infty} \sum_{j_1^*=0}^{j_1} \sum_{m_1=j_1}^{\infty} \sum_{m_{1,1}=m_1}^{\infty} 2^{j_1^* - 3j_1 + (2m_1 - 3m_{1,1})/4}. \tag{6.3}$$

Now consider subcase (b). By Proposition 1 applied directly to $G_1$, we have that the boxes $B(\mathbf{w}_{1,1}, 2^{-m_{1,1}-s}n)$, $B(\mathbf{w}_1, 2^{-m_1-2s-2}n)$ and $B(\mathbf{x}_1, 2^{-m_1-s}n)$ are mutually disjoint under (i). Then since $m_{1,1}$ is at most a constant different than $m_1$ by (b) and (4.4), by independence and Lemma 2 alone, we obtain that (6.2) continues to hold with $m_1$ in place of $m_{1,1}$. Therefore, by substitution of $m_1$ for $m_{1,1}$ also in the size estimate of $a_{m_{1,1}}$ and by eliminating the sum on $m_{1,1}$, we obtain

$$Iib \le C_\varepsilon n^{9/4+3\varepsilon} \sum_{j_1=0}^{\infty} \sum_{j_1^*=0}^{j_1} \sum_{m_1=j_1}^{\infty} 2^{j_1^*-3j_1-m_1/4}. \tag{6.4}$$

To help with case (ii), as well as further cases arising in higher moment calculations, we first state a general consequence of Proposition 1 that we will use to set up our application of Lemma 7.

PROPOSITION 2. *Consider the graph $G_1 = \{\mathbf{x}_1, \mathbf{w}_1, \ldots\}$. Set $D := 2^{-j_1-2c}n$. There exists a constant $c_0$ and a rectangle $R$ centered at $\tilde{\mathbf{x}}_1$ with smallest half-side of length $l$ satisfying $D/10 \le l \le D/5$ and largest half-side of length $L$ satisfying $L/l \le 2$ and with center satisfying $\|\tilde{\mathbf{x}}_1 - \mathbf{x}_1\| \le l/2$ such that each box $B'(\mathbf{w}) := B(\mathbf{w}, 2^{-m(\mathbf{w})-s-c_0}n)$, $\mathbf{w} \in G_1$, $\mathbf{w} \ne \mathbf{x}_1$, lies either entirely inside $R$ or entirely outside $R$.*

PROOF. First, by Proposition 1, the boxes $B(\mathbf{w}) = B(\mathbf{w}, 2^{-m(\mathbf{w})-s}n)$, $\mathbf{w} \in G_1$, are mutually disjoint. Since the sum of the radii of the boxes $B(\mathbf{w})$, $\mathbf{w} \ne \mathbf{x}_1$, is bounded by $(\tau - 1)D$, we may choose $c_0$ so large that the sum of the diameters of the corresponding shrunken boxes $B'(\mathbf{w})$ is at most $D/16$. Therefore, by the same argument as given in the proof of Lemma 7 in the Appendix for the construction of $\widetilde{R}$, with $D$ here playing the role there of the distance $D_v$, the proof is complete. □

We now continue our discussion of case (ii). We apply Proposition 2 directly to the graph $G_1$ to obtain the rectangle $R$ having the properties stated there so that, in particular, $B(\mathbf{x}_1, l/2) \subset R \subset B(\mathbf{x}_1, D/5)$ for $D := 2^{-j_1-2c}n$. We apply Lemma 7 with $v = 2$, $\mathbf{x} = \mathbf{x}_1$, and $2^{\lambda_i} = 2^{-m(\mathbf{w}_i)-s-c_0}n$. Hence, by Lemma 2, and by independence applied to any shrunken box lying outside $R$, we have that

$$\begin{aligned} p_{n,3}(\mathbf{x}_1, \mathbf{w}_1, \mathbf{w}_2) &\le CP(\mathcal{U}_4(\mathbf{0}, l/d)) \prod_{i=1}^{2} P(\mathcal{U}_4(\mathbf{0}, 2^{\lambda_i - c_1}n)) \\ &\le C_\varepsilon n^{-15/4+3\varepsilon} 2^{5(j_1+m_1+m_2)/4}. \end{aligned} \tag{6.5}$$

As in case (i), we may construct a horseshoe with inner radius $2^\rho$ and outer radius $2^\nu$ as chosen above, whose inner box contains all the boxes implied by the estimate (6.5). So we may improve this estimate by the same factor coming from Lemma 5 as before. Therefore, by (6.5), this repeated observation, and the size estimate (3.12) for each index $m = m_1$ and $m = m_2$, we obtain

$$Iii \le C_\varepsilon n^{9/4+3\varepsilon} \sum_{j_1=0}^{\infty} \sum_{j_1^*=0}^{j_1} \sum_{m_2=j_1}^{\infty} \sum_{m_1=m_2}^{\infty} 2^{j_1^*+(-7j_1-3m_2-3m_1)/4}. \tag{6.6}$$

6.1.2. $\tau = 3$, $r = 2$. We first assume that $|V_1| = 1$ and $|V_2| = 0$. For simplicity, we assume that $\mathbf{x}_3$ is the second root. Thus, the graph $\mathbf{G}_1$ has vertices $\{\mathbf{x}_1, \mathbf{w}_1\}$ and the graph $\mathbf{G}_2$ is trivial over the (isolated) root vertex $\{\mathbf{x}_3\}$. Next we determine whether we have a horseshoe relationship between the two root vertices or not. If $\mathbf{x}_3 K \mathbf{x}_1$, then we have $h = 1$, else we have $h = 2$.

We work first with the case $h = 1$. Similar to our analysis of the corresponding case of the second moment estimation, we have two possibilities under $h = 1$: either (a), $j_1 \le j_3 \le j_1 + 2c + 2$, or (b), $j_3 > j_1 + 2c + 2$. We write $II_1$ to denote the sub-sum of $\Sigma_0$ that corresponds to $r = 2$ and $h = 1$. We also write $II_1 = II_1a + II_1b$, where the sums $II_1a$ and $II_1b$ correspond, respectively, to the cases (a) and (b). We study first case (a). We apply Proposition 2 to the graph $G_1$ to obtain a rectangle $R$ such that $B(\mathbf{x}_1, l/2) \subset R \subset B(\mathbf{x}_1, D/5)$ for $D := 2^{-j_1-2c}n$ and $l \ge D/10$, so that $B'(\mathbf{w}_1)$ lies either inside or outside $R$. By construction of the original roots and by the cases $m_1 = j_1 + 2c$ or $m_1 \ge j_1 + 2c + 1$, we find that the box $B(\mathbf{x}_3, 2^{-j_3-4c-2}n)$ is disjoint from both $R$ and the box $B(\mathbf{w}_1, 2^{-m_1-2c-2}n)$ [see (7.4)]. Hence, we can apply independence and Lemmas 2 and 7 to estimate

$$\begin{aligned} p_{n,3}(\mathbf{x}_1, \mathbf{w}_1, \mathbf{x}_3) &\le P(\mathcal{U}_4(\mathbf{0}, l/d))P(\mathcal{U}_4(\mathbf{0}, C2^{-m_1}n))P(\mathcal{U}(\mathbf{0}, C2^{-j_3}n)) \\ &\le C_\varepsilon n^{-15/4+3\varepsilon} 2^{5(2j_1+m_1)/4}, \end{aligned} \tag{6.7}$$

where the last inequality holds because under (a), $j_3$ is within a constant of $j_1$. Since $h = 1$, it is again an easy matter to construct a horseshoe as in each case of Section 6.1.1 with inner radius $2^\rho \asymp 2^{-j_1}n$ and outer radius $2^\nu \asymp 2^{-j_1^*}n$ whose inner box contains all the boxes implied by the estimate (6.7). So we may improve the estimate (6.7) by an application of Lemma 5. Hence, because $|a_{m_1}| \le C2^{-2m_1}n^2$ and since there are only $C2^{-2j_1}n^2$ vertices $\mathbf{x}_3$ to account for when $h = 1$ and (a) holds, we have that

$$II_1a \le C_\varepsilon n^{9/4+3\varepsilon} \sum_{j_1=0}^{\infty} \sum_{j_1^*=0}^{j_1} \sum_{m_1=j_1}^{\infty} 2^{j_1^*+(-10j_1-3m_1)/4}. \tag{6.8}$$

We next study case (b). Since $h=1$, we have that the rectangle $R$ exists as constructed above. But now because $j_3$ is sufficiently larger than $j_1$ and since $\mathbf{x}_3$ is isolated, there is room to construct a pair of boxes $B_1$ and $B_2$ as in the context of Lemma 5 with $B_1$ centered at $\mathbf{x}_3$ as follows. We find $B_1 = B(\mathbf{x}_3, 2^{\rho_3}) \subset B(n)$, with $2^{\rho_3} \asymp 2^{-j_3}n$ such that the right-hand side of $B_1$ lies in $\partial B(n)$. We take $B_2$ accordingly by defining its radius as $2^{\nu_3} \asymp 2^{-j_1}n$ such that $B_2$ is disjoint from both $B'(\mathbf{w}_1)$ and $R$, where $B'(\mathbf{w}_1)$ is the shrunken box given by Proposition 2. Hence, by independence and Lemmas 2, 5 and 7, we find

$$\begin{aligned} p_{n,3}(\mathbf{x}_1, \mathbf{w}_1, \mathbf{x}_3) &\le P(\mathcal{U}_4(\mathbf{0}, l/d)) \\ &\quad \times P(\mathcal{U}_4(\mathbf{0}, C2^{-m_1}n))P(\mathcal{U}_4(\mathbf{0}, 2^{\rho_3}))P(\mathcal{J}(\rho_3, \nu_3)) \\ &\le C_\varepsilon n^{-5/2+2\varepsilon} 2^{(j_1+m_1+j_3)(5/4-\varepsilon)} 2^{(j_1-j_3)(2-\varepsilon)}. \end{aligned} \tag{6.9}$$

Again we may improve this estimate by introducing a horseshoe whose inner box contains $R$, $B'(\mathbf{w}_1)$ and the horseshoe pair $(B_1, B_2)$. So we multiply the right-hand side of (6.9) by the factor $2^{(\rho-\nu)(2-\varepsilon)}$, where the radii $2^\rho$ and $2^\nu$ are defined up to multiplicative constants in the previous case (a). Since by $h=1$, there are only $C2^{-j_1-j_3}n^2$ vertices $\mathbf{x}_3$ accounted for with $\mathbf{x}_1 \in A_{j_1}$ and $\mathbf{x}_3 \in A_{j_3}$, we find by these observations that

$$II_1 b \le C_\varepsilon n^{9/4+3\varepsilon} \sum_{j_1=0}^{\infty} \sum_{j_1^*=0}^{j_1} \sum_{m_1=j_1}^{\infty} \sum_{j_3=j_1}^{\infty} 2^{j_1^*+(-3j_1-3m_1-7j_3)/4}. \tag{6.10}$$

We now pass to the case $h=2$. Note that the sum over $\mathbf{x}_3$ in $II_2$ is no longer localized strictly nearby $\mathbf{x}_1$ via the horseshoe relation, so we will be able to construct a larger horseshoe at $\mathbf{x}_3$, but how large now depends on the relation between the dual indices of the original two roots. Let $\mathbf{x}_i \in A^*_{j_i^*}$, $i = 1, 3$. There are two cases to consider regarding the dual indices:

$$\text{either (a) } |j_1^* - j_3^*| \le 1 \quad \text{or} \quad \text{(b) } |j_1^* - j_3^*| \ge 2.$$

Define sums $II_2 a$ and $II_2 b$ by partitioning the sum $II_2$ according to these cases.

We work first with case (a). As in Section 5, we have $\mathbf{x}_3 \in A_{j_3, j_1^*}$ where the latter set is defined by (5.28). For convenience, we rewrite the bands $b_v = b_v(\mathbf{x}_1, j_1, j_3, j_1^*)$ of (5.29) by

$$b_v := \{(y_1, y_2) \in A'_{j_3, j_1^*} : v2^{-j_1+5}n < |(\mathbf{x}_1)_2 - y_2| \le (v+1)2^{-j_1+5}n\}$$

for each $v = 1, 2, \dots, v_{\max}$, with $v_{\max} \asymp 2^{j_1 - j_1^*}$. Here we have simply substituted $\mathbf{x}_1$, $j_1$, $j_3$ and $j_1^*$ for $\mathbf{x}$, $j$, $k$ and $j^*$, respectively, in the original definition. Notice by the definition of the horseshoe relation that $\mathbf{x}_3 \in b_v$ only for $v \ge 2$. We utilize the rectangle $R$ and the box $B'(\mathbf{w}_1)$ we have constructed for $h=1$. First choose the radius $2^{\rho_3} \asymp 2^{-j_3}n$ such that the right edge of $B_1(\mathbf{x}_3) := B(\mathbf{x}_3, 2^{\rho_3})$ just meets $\partial B(n)$. Also choose $2^{\rho_1} \asymp 2^{-j_1}n$ such that the right edge of $B_1(\mathbf{x}_1) := B(\mathbf{x}_1, 2^{\rho_1})$ just meets $\partial B(n)$. These are the inner boxes of horseshoes we will construct at each of $\mathbf{x}_3$ and $\mathbf{x}_1$, respectively. Note that $B_1(\mathbf{x}_1)$, in fact, contains both $R$ and $B'(\mathbf{w}_1)$ because these latter sets are chosen via the constant $c$ of Section 4 to both lie within a box of radius $2^{-j_1-c}n$ centered at $\mathbf{x}_1$. For each $v = 2, \dots, v_{\max}$, we define $2^{\nu} \asymp v2^{-j_1}n$, uniformly in $v \geq 2$, so that boxes $B_2(\mathbf{x}_1) \subset B(n)$ and $B_2(\mathbf{x}_3) \subset B(n)$, each with radius $2^{\nu}$, exist and are disjoint such that $(B_1(\mathbf{x}_1), B_2(\mathbf{x}_1))$ and $(B_1(\mathbf{x}_3), B_2(\mathbf{x}_3))$ each form a horseshoe pair as in the context of Lemma 5. The outer boxes remain in $B(n)$ by (5.28) and the expression for $v_{\max}$. Moreover, the outer boxes $B_2(\mathbf{x}_1)$ and $B_2(\mathbf{x}_3)$, while disjoint, are nested inside another box $\widetilde{B}_1$ of radius $C2^{\nu}$ whose right edge also lies in $\partial B(n)$. Since we are in case (a), we may again pair $\widetilde{B}_1$ with an outer horseshoe box $\widetilde{B}_2$ of radius $2^{\nu_1} \asymp 2^{-j_1^*}n$. Therefore, by independence and by application of Lemma 5 to the horseshoe pairs, and by an application of Lemma 7 as in the case $h = 1$, for $\mathbf{x}_3 \in b_v$, we have that, for all $v \geq 2$, $p_n := p_{n,3}(\mathbf{x}_1, \mathbf{w}_1, \mathbf{x}_3)$ satisfies

$$\begin{aligned} p_n &\leq C_\varepsilon n^{-5/2+2\varepsilon} 2^{(j_1+m_1+j_3)(5/4-\varepsilon)} 2^{(\rho_1+\rho_3-\nu-\nu_1)(2-\varepsilon)} \\ &\leq C_\varepsilon n^{-15/4+3\varepsilon} v^{-2+\varepsilon} 2^{2j_1^*+(5j_1+5m_1-3j_3)/4}. \end{aligned} \tag{6.11}$$

Hence, by (6.11), (5.16), (3.12) and the estimate $|b_v| \leq C2^{-j_1-j_3}n^2$ [cf. (5.30)], we have

$$II_2 a \leq C_\varepsilon n^{9/4+3\varepsilon} \sum_{j_1=0}^{\infty} \sum_{j_1^*=0}^{j_1} \sum_{m_1=j_1}^{\infty} \sum_{j_3=j_1}^{\infty} \sum_{v=1}^{\infty} v^{-2+\varepsilon} 2^{j_1^*+(-3j_1-3m_1-7j_3)/4}. \tag{6.12}$$

Consider next case (b) under $h = 2$. The difference with case (a) is that now the outer boxes $B_2(\mathbf{x}_1)$ and $B_2(\mathbf{x}_3)$ found there may be chosen with larger radii while still remaining disjoint. Thus, the horseshoe pair $(\widetilde{B}_1, \widetilde{B}_2)$ is no longer needed in this case. Indeed, we may now take the radii of these outer boxes as $2^{\nu_1} \asymp 2^{-j_1^*}n$ and $2^{\nu_3} \asymp 2^{-j_3^*}n$, respectively. By definition of the dual partition, these larger boxes still remain in $B(n)$ and are disjoint. Therefore, we obtain in place of (6.11) the bound

$$\begin{aligned} p_n &\leq C_\varepsilon n^{-5/2+2\varepsilon} 2^{(j_1+m_1+j_3)(5/4-\varepsilon)} 2^{(\rho_1+\rho_3-\nu_1-\nu_3)(2-\varepsilon)} \\ &\leq C_\varepsilon n^{-15/4+3\varepsilon} 2^{2j_1^*+2j_3^*+(-3j_1+5m_1-3j_3)/4}. \end{aligned} \tag{6.13}$$

Hence, by (6.13), we have

$$II_2 b \leq C_\varepsilon n^{9/4+3\varepsilon} \sum_{j_1=0}^{\infty} \sum_{m_1=j_1}^{\infty} \sum_{j_3=j_1}^{\infty} \sum_{j_1^*=0}^{j_1} \sum_{j_3^*=0}^{j_3} 2^{j_1^*+j_3^*+(-7j_1-3m_1-7j_3)/4}. \tag{6.14}$$

Finally, in the case that the isolated root vertex is instead $\mathbf{x}_1$ (so $V_1 = \varnothing$) and the second root vertex is $\mathbf{x}_2$ (so $V_2 = \{\mathbf{x}_3\}$), we obtain a wholly analogous estimation by writing out the cases (a) $j_2 \le j_1 + 2c + 2$, and (b) $j_2 > j_1 + 2c + 2$ under $h = 1$ and by writing out the cases (a) $|j_2^* - j_1^*| \le 1$, and (b) $|j_2^* - j_1^*| > 1$ under $h = 2$.

6.1.3. $\tau = 3$, $r = 3$. In the case $r = 3$ we must have three isolated roots relative to the original partition $\{A_j\}$ that are simply $\mathbf{x}_1$, $\mathbf{x}_2$ and $\mathbf{x}_3$. The basic plan in all that follows is that a horseshoe must be constructed whenever there is room at a given level of algebraic or dyadic division to do so. Algebraic levels of division arise according to placement of a root vertex in a band $b_v$ or in another closely related band $b'_u$ that we define in Section 6.1.4. The horseshoe structure depends on the room that exists between root vertices. This spacing will be accounted for by various joint inequalities in the dyadic indices $j_i$, or in another spacing relation that we introduce for the dual indices $j_i^*$ in Section 6.1.5.

6.1.4. $\tau = 3$, $r = 3$, $h = 1$. We assume first that $h = 1$. The set of root vertices chained to $\mathbf{x}_1$ by the horseshoe relationship is therefore $U_1 = \{\mathbf{x}_2, \mathbf{x}_3\}$. We write $III_1$ to denote the sub-sum of $\Sigma_0$ that corresponds to $r = 3$ and $h = 1$. We write four conditions:

$$\begin{aligned} &\text{(a1) } j_1 \le j_2 \le j_1 + 2c + 8, \quad \text{(b1) } j_2 > j_1 + 2c + 8, \\ &\text{(a2) } j_2 \le j_3 \le j_2 + 2c + 8, \quad \text{(b2) } j_3 > j_2 + 2c + 8. \end{aligned} \tag{6.15}$$

We also write $III_1 = III_1a1a2 + III_1b1a2 + III_1a1b2 + III_1b1b2$, where the summands correspond, respectively, to these four joint cases.

For any vertex $\mathbf{x}_i \in U_1$, we define bands of vertices $b'_u = b'_u(\mathbf{x}_i, j_i)$, for all $u = 0, 1, \ldots$ and $i \ge 2$, that divide $R_0 := B(\mathbf{x}_1, 2^{-j_1+6+\tau} n) \cap B(n)$ into horizontal sections by

$$b'_u := \{(y_1, y_2) \in R_0 : u 2^{-j_i+5} n < |(\mathbf{x}_i)_2 - y_2| \le (u+1) 2^{-j_i+5} n\}, \tag{6.16}$$

where $u$ ranges up to $u_{\max} \asymp 2^{j_i - j_1}$ for $i \ge 2$. Here the exponent in the definition of $R_0$ is chosen such that, by the definition of the horseshoe relation, any vertex in $U_1$ lies in $R_0$. The main difference between the bands $b'_u$ and the bands $b_v$ that we defined in (5.29) is that, contrary to that definition, here we place no restriction that $b'_u$ lie in some single annulus $A_{j_k}$. Although these new bands play a similar role as the original ones, we apply them at a different level of the construction of estimates. We apply, in general, the $b'_u$ within a horseshoe set $U_a$ with $U_a \ne \varnothing$. We apply the $b_v$ instead in a region between such horseshoe sets. Now, by (6.16), for any annulus $A_{j_k}$, we have that

$$|b'_u(\mathbf{x}_i, j_i) \cap A_{j_k}| \le C 2^{-j_i - j_k} n^2 \qquad \text{for all } u \ge 0. \tag{6.17}$$

Study first the joint case (a1)–(a2). In this case all three vertices are located either in $A_{j_1}$ or a nearby annulus. Therefore, since $h=1$, $\mathbf{x}_2$ and $\mathbf{x}_3$ are each confined to a set of vertices of size at most $C2^{-2j_1}n^2$. We call such a size estimate a confinement factor. Since the roots are isolated, we may construct disjoint boxes with centers at the vertices $\mathbf{x}_i$, $i=1,2,3$, such that each has a radius $2^{\rho_1} \asymp 2^{-j_1}n$ and each lies in $R_0$. Construct a box $B_1$ of radius $2^{\rho} \asymp 2^{-j_1}n$ whose right edge meets $\partial B(n)$ and that contains $R_0$. Pair $B_1$ with an associated outer horseshoe box $B_2$ of radius $2^{\nu} \asymp 2^{-j_1^*}n$. Hence, by independence and Lemmas 2 and 5, and the confinement factors, we easily have

$$III_1 a1a2 \le C_\varepsilon n^{9/4+3\varepsilon} \sum_{j_1=0}^{\infty} \sum_{j_1^*=0}^{j_1} 2^{j_1^* - 13j_1/4}. \tag{6.18}$$

Under the joint case (b1)–(a2) we study two subcases,

$$\text{either (i) } \mathbf{x}_3 \in b_0'(\mathbf{x}_2, j_2) \quad \text{or} \quad \text{(ii) } \mathbf{x}_3 \in b_u'(\mathbf{x}_2, j_2) \qquad \text{for some } u \ge 1. \tag{6.19}$$

Partition the sum $III_1 b1a2$ accordingly: $III_1 b1a2 = III_1 b1a2i + III_1 b1a2ii$. Study first subcase (i) under (b1)–(a2). Since the vertices are isolated roots, by (a2), we can choose radii $2^{\sigma} \asymp 2^{-j_2}n$ and $2^{\rho} \asymp 2^{-j_1}n$ such that the boxes $B(\mathbf{x}_1, 2^{\rho})$, $B(\mathbf{x}_2, 2^{\sigma})$ and $B(\mathbf{x}_3, 2^{\sigma})$ lie in $B(n)$ and are mutually disjoint. Moreover, by (b1) and (i), the boxes centered at $\mathbf{x}_2$ and $\mathbf{x}_3$ are both contained in a box $B_1$ of radius $C2^{\sigma} \le 2^{-j_1-2}n$ whose right edge lies on $\partial B(n)$, where we choose $2^{\rho}$ such that $B_1$ is also disjoint from the box centered at $\mathbf{x}_1$. We construct a second box $B_2$ so that the pair $(B_1, B_2)$ conforms to the context of Lemma 5 where the outer box has radius $C2^{\rho}$ and is disjoint from the box centered at $\mathbf{x}_1$. We also construct an inner horseshoe box $\widetilde{B}_1$ of radius $2^{\rho_1} \asymp 2^{-j_1}n$ that contains all the boxes constructed so far and pair it with an outer horseshoe box $\widetilde{B}_2$ of radius $2^{\nu_1} \asymp 2^{-j_1^*}n$. Thus, by independence and Lemma 2 and two applications of Lemma 5, we have

$$\begin{aligned} p_{n,3}(\mathbf{x}_1, \mathbf{x}_2, \mathbf{x}_3) &\le C_\varepsilon 2^{(\rho+2\sigma)(-5/4+\varepsilon)+(\sigma+\rho_1-\rho-\nu_1)(2-\varepsilon)} \\ &\le C_\varepsilon n^{-15/4+3\varepsilon} 2^{2j_1^*+(5j_1+2j_2)/4}. \end{aligned} \tag{6.20}$$

Since under $h=1$ and (a2) the confinement factor for $\mathbf{x}_2$ is $C2^{-j_1-j_2}n^2$, and since under the added condition (i) the confinement factor for $\mathbf{x}_3$ is $C2^{-2j_2}n^2$, we obtain by (6.20) that

$$III_1 b1a2i \le C_\varepsilon n^{9/4+3\varepsilon} \sum_{j_1=0}^{\infty} \sum_{j_1^*=0}^{j_1} \sum_{j_2=j_1}^{\infty} 2^{j_1^*+(-3j_1-10j_2)/4}. \tag{6.21}$$

Consider next subcase (ii) under (b1)–(a2). The difference with case (i) is that we now create two disjoint inner horseshoes instead of just one.

We take inner boxes $B_1(\mathbf{x}_i)$ centered at $\mathbf{x}_i$, $i=2,3$, with radii $2^{\rho_i}\asymp 2^{-j_2}n$, $i=2,3$, such that these inner boxes meet the boundary of $B(n)$ and are disjoint and are, moreover, disjoint from a box centered at $\mathbf{x}_1$ with radius $2^{\rho}\asymp 2^{-j_1}n$. We take the associated outer boxes $B_2(\mathbf{x}_i)$, $i=2,3$, each with a radius $2^{\nu}\asymp u2^{-j_2}n$, for $u$ ranging up to order $2^{j_2-j_1}$. Both outer boxes of these horseshoes are disjoint from $B(\mathbf{x}_1,2^{\rho})$ by (b1). Further, we construct a third horseshoe by taking an inner box $B_1$ of radius $C2^{\nu}$ that contains both the outer boxes $B_2(\mathbf{x}_i)$, $i=2,3$, and that admits an outer box $B_2$ of radius $2^{\rho_1}\asymp 2^{-j_1}n$ that is still disjoint from the box $B(\mathbf{x}_1,2^{\rho})$. Finally, we construct a fourth horseshoe pair $(\widetilde{B}_1,\widetilde{B}_2)$ such that $\widetilde{B}_1$ contains all the previous outer boxes, as well as the box $B(\mathbf{x}_1,2^{\rho})$. We take the inner and outer radii of this last pair to be, respectively, $C2^{\rho_1}$ and $2^{\nu_1}\asymp 2^{-j_1^*}n$. Therefore, we obtain

$$\begin{aligned}p_{n,3}(\mathbf{x}_1,\mathbf{x}_2,\mathbf{x}_3)&\le C_\varepsilon 2^{(\rho+\rho_2+\rho_3)(-5/4+\varepsilon)+(\rho_2+\rho_3-\nu-\nu_1)(2-\varepsilon)}\\&\le C_\varepsilon n^{-15/4+3\varepsilon}u^{-2+\varepsilon}2^{2j_1^*+(5j_1+2j_2)/4}.\end{aligned}\tag{6.22}$$

Now by (a2) and (6.17), we have that $|b_u'\cap A_{j_3}|\le C2^{-2j_2}n^2$. Also, $\mathbf{x}_2$ is confined to a region of size $C2^{-j_1-j_2}n^2$. Therefore, by independence and Lemmas 2 and 5, we obtain by (6.22) the estimate

$$\begin{aligned}III_1b1a2ii&\le C_\varepsilon n^{9/4+3\varepsilon}\\&\quad\times\sum_{j_1=0}^{\infty}\sum_{j_1^*=0}^{j_1}\sum_{j_2=j_1}^{\infty}\sum_{u=1}^{\infty}u^{-2+\varepsilon}2^{j_1^*+(-3j_1-10j_2)/4}.\end{aligned}\tag{6.23}$$

We now consider the joint case (a1)–(b2). Since $j_3$ is sufficiently larger than $j_2$ and the roots are isolated, we can construct a horseshoe at $\mathbf{x}_3$ with inner radius $2^{\rho_3}\asymp 2^{-j_3}n$ and outer radius $2^{\nu}\asymp 2^{-j_2}n\asymp 2^{-j_1}n$ and choose the radius $2^{\rho}\asymp 2^{-j_1}n$ so that the outer box of this horseshoe will be disjoint from both the boxes $B(\mathbf{x}_1,2^{\rho})$ and $B(\mathbf{x}_2,2^{\rho})$ that lie in $B(n)$ and are themselves constructed to be disjoint. Again, we construct a large horseshoe pair $(\widetilde{B}_1,\widetilde{B}_2)$ with inner and outer radii $2^{\rho_1}$ and $2^{\nu_1}$, respectively, as in the previous cases such that $\widetilde{B}_1$ contains the smaller horseshoe, as well as the boxes $B(\mathbf{x}_1,2^{\rho})$ and $B(\mathbf{x}_2,2^{\rho})$. Therefore, by independence and Lemmas 2 and 5, and by (a1), we have that

$$\begin{aligned}p_{n,3}(\mathbf{x}_1,\mathbf{x}_2,\mathbf{x}_3)&\le C_\varepsilon 2^{(2\rho+\rho_3)(-5/4+\varepsilon)+(\rho_1+\rho_3-\nu-\nu_1)(2-\varepsilon)}\\&\le C_\varepsilon n^{-15/4+3\varepsilon}2^{2j_1^*+(10j_1-3j_3)/4}.\end{aligned}\tag{6.24}$$

By $h=1$ and (a1), we have that $\mathbf{x}_2$ is confined by the factor $2^{-2j_1}n^2$, while by $h=1$ alone, $\mathbf{x}_3$ is confined by the factor $2^{-j_1-j_3}n^2$. Therefore, since by (a1) we eliminate the sum over $j_2$, we obtain by (6.24) that

$$III_1a1b2\le C_\varepsilon n^{9/4+3\varepsilon}\sum_{j_1=0}^{\infty}\sum_{j_1^*=0}^{j_1}\sum_{j_3=j_1}^{\infty}2^{j_1^*+(-6j_1-7j_3)/4}.\tag{6.25}$$

Consider finally the joint case (b1)–(b2). Again we apply the dichotomy (6.19). We partition the sum $III_1 b1b2 = III_1 b1b2i + III_1 b1b2ii$ accordingly. In subcase (i), under (b1)–(b2), we take a horseshoe at $\mathbf{x}_3$ with inner radius $2^{\rho_3} \asymp 2^{-j_3} n$ and outer radius $2^{\nu_3} \asymp 2^{-j_2} n$, where the outer box is disjoint from the boxes centered at $\mathbf{x}_2$ and $\mathbf{x}_1$ of radii $2^{\sigma} \asymp 2^{-j_2} n$ and $2^{\rho} \asymp 2^{-j_1} n$, respectively. These last two boxes are chosen to be small enough that, even doubling their radii, they would not meet $\partial B(n)$. We next take a horseshoe $(B_1, B_2)$ such that $B_1$ has a radius $2^{\rho_2} \asymp 2^{-j_2} n$ and contains both the box centered at $\mathbf{x}_2$ and the horseshoe at $\mathbf{x}_3$. We take $B_2$ to have a radius $2^{\nu_2} \asymp 2^{-j_1} n$ that is disjoint from the box centered at $\mathbf{x}_1$. Again we construct a large horseshoe pair $(\widetilde{B}_1, \widetilde{B}_2)$ with inner and outer radii $2^{\rho_1}$ and $2^{\nu_1}$, respectively, as in the previous cases such that $\widetilde{B}_1$ contains both the smaller nested horseshoes, as well as the box $B(\mathbf{x}_1, 2^{\rho})$. Therefore, we obtain by independence and Lemmas 2 and 5 that

$$\begin{aligned} p_{n,3}(\mathbf{x}_1, \mathbf{x}_2, \mathbf{x}_3) &\le C_\varepsilon 2^{(\rho+\sigma+\rho_3)(-5/4+\varepsilon)+(\rho_1+\rho_2+\rho_3-\nu_1-\nu_2-\nu_3)(2-\varepsilon)} \\ &\le C_\varepsilon n^{-15/4+3\varepsilon} 2^{2j_1^* + (5j_1+5j_2-3j_3)/4}. \end{aligned} \tag{6.26}$$

By $h = 1$, we have that $\mathbf{x}_2$ is confined by the factor $2^{-j_1-j_2} n^2$, while in addition by (i), $\mathbf{x}_3$ is confined by the factor $2^{-j_2-j_3} n^2$. Therefore, we obtain by (6.26) that

$$III_1 b1b2i \le C_\varepsilon n^{9/4+3\varepsilon} \sum_{j_1=0}^{\infty} \sum_{j_1^*=0}^{j_1} \sum_{j_2=j_1}^{\infty} \sum_{j_3=j_1}^{\infty} 2^{j_1^* + (-3j_1-3j_2-7j_3)/4}. \tag{6.27}$$

Finally, in subcase (ii), under (b1)–(b2), we take inner boxes $B_1(\mathbf{x}_i)$ centered at $\mathbf{x}_i$, $i = 2, 3$, with radii $2^{\rho_i} \asymp 2^{-j_i} n$, $i = 2, 3$, such that these inner boxes meet the boundary of $B(n)$ and are disjoint and are, moreover, disjoint from a box centered at $\mathbf{x}_1$ with radius $2^{\rho} \asymp 2^{-j_1} n$. We take the associated outer boxes $B_2(\mathbf{x}_i)$, $i = 1, 2$, each with a radius $2^{\nu} \asymp u 2^{-j_2} n$, for $u$ ranging up to order $2^{j_2 - j_1}$. Both outer boxes of these horseshoes are disjoint from $B(\mathbf{x}_1, 2^{\rho})$ by (b1). This almost looks like subcase (ii) under (b1)–(a2), except notice that here the box $B_1(\mathbf{x}_3)$, while still having a radius distinct from the box $B_1(\mathbf{x}_2)$, has now an asymptotically smaller radius since we are in case (b2). All the remaining arrangements of boxes and horseshoes are exactly as in subcase (ii) of (b1)–(a2), with the same formulae for asymptotic radii. Thus, we have a total of four horseshoes. Therefore, we obtain

$$\begin{aligned} p_{n,3}(\mathbf{x}_1, \mathbf{x}_2, \mathbf{x}_3) &\le C_\varepsilon 2^{(\rho+\rho_2+\rho_3)(-5/4+\varepsilon)+(\rho_2+\rho_3-\nu-\nu_1)(2-\varepsilon)} \\ &\le C_\varepsilon n^{-15/4+3\varepsilon} u^{-2+\varepsilon} 2^{2j_1^* + (5j_1+5j_2-3j_3)/4}. \end{aligned} \tag{6.28}$$

Now by (6.17), we have that $|b_u' \cap A_{j_3}| \le C 2^{-j_2-j_3} n^2$. Also, $\mathbf{x}_2$ is confined to a region of size $C 2^{-j_1-j_2} n^2$ by $h = 1$. Therefore, by independence and

Lemmas 2 and 5, we obtain by (6.22) the estimate

(6.29)
$$III_1b1b2ii \le C_\varepsilon n^{9/4+3\varepsilon} \times \sum_{j_1=0}^{\infty}\sum_{j_1^*=0}^{j_1}\sum_{j_2=j_1}^{\infty}\sum_{j_3=j_2}^{\infty}\sum_{u=1}^{\infty} u^{-2+\varepsilon}2^{j_1^*+(-3j_1-3j_2-7j_3)/4}.$$

This concludes our analysis of the case $h=1$.

6.1.5. *$\tau=3$, $r=3$, $h=2$.* We consider next that $h=2$. We assume first that $U_1=\{\mathbf{x}_2\}$ so that we have root horseshoe vertices $\mathbf{x}_1$ and $\mathbf{x}_3$. Again, because $r=3$, all the roots are isolated and, moreover, the root $\mathbf{x}_3$ is an isolated root horseshoe vertex. We write $III_2$ for the part of the sum $\Sigma_0$ corresponding to this arrangement. Recall that we locate the vertices in the dual partition, in general, by the dual indices $j_i^*$ such that $\mathbf{x}_i \in A^*_{j_i^*}$, $i=1,2,\dots$.

Perhaps by now part of the outline is clear. Initially, we consider two possibilities:

$$\text{either (a) } j_2 \le j_1+2c+8 \quad \text{or} \quad \text{(b) } j_2 > j_1+2c+8.$$

However, while we will construct a horseshoe at $\mathbf{x}_3$ in either case, we must delineate its size. Its outer radius may be as small as $Cv2^{-j_2}n$ when $\mathbf{x}_3 \in b_v(\mathbf{x}_2,j_2,j_3,j_1^*)$, and it may be as large as $C2^{-j_1^*}$ when $\mathbf{x}_3 \in A^*_{j_3^*}$ for $|j_3^*-j_1^*| \ge 2$. So, in general, we need to know the manner in which the root horseshoe vertices $\mathbf{x}_{f_a}$ are separated in the dual partition. We have already seen such an analysis in case $r=2$ and $h=2$. We generalize the approach shown there. For any dual indices $k^*$ and $j^*$ of root horseshoe vertices (we call such indices also as dual horseshoe indices), write that $k^*J^*j^*$ if $|k^*-j^*| \le 1$. We say that a dual index $j^*_{f'}$ of a root horseshoe vertex $\mathbf{x}_{f'}$ is chained to the dual index $j^*_f$ of another root horseshoe vertex $\mathbf{x}_f$ if there exists a sequence of $J^*$ relations from $j^*_{f'}$ to $j^*_f$. By the method of Section 4, we define root dual horseshoe indices $j^*_{f_1^*} < j^*_{f_2^*} < \cdots < j^*_{f_t^*}$, with $f_1^*=f_1$ and some $t \le h$, where $\mathbf{x}_{f_1}=\mathbf{x}_1$ is the first root horseshoe vertex.

In the current case we have either $t=1$ or $t=2$, where $t=1$ means that $|j_3^*-j_1^*| \le 1$ and $t=2$ means that the opposite inequality holds. We write $III_{h,t}$ for the part of the sum $\Sigma_0$ corresponding to $r=3$ and the given values of $h$ and $t$. Since here $h=2$ and $t=1$ or 2, we have $III_2 = III_{2,1}+III_{2,2}$. We further partition $III_{2,1}=III_{2,1}a+III_{2,1}b$ and $III_{2,2}=III_{2,2}a+III_{2,2}b$ for the arrangements of vertices corresponding, respectively, to cases (a) and (b).

We consider first an estimate of $III_{2,1}a$ so that, in particular, $t=1$. Put $b_v := b_v(\mathbf{x}_2,j_2,j_3,j_1^*)$. Since $h=2$ and $t=1$, we have $\mathbf{x}_3 \in b_v$ for some $v \ge 2$ (cf. the case $h=2$ of Section 6.1.2). Here $v$ ranges up to $v_{\max} \asymp 2^{j_2-j_1^*} \le$ $C2^{j_1-j_1^*}$ by (a). We choose a radius $2^{\rho_3} \asymp 2^{-j_3}n$, so that the right boundary of the box $B_1(\mathbf{x}_3) := B(\mathbf{x}_3, 2^{\rho_3})$ meets $\partial B(n)$. This is the inner box of a horseshoe at $\mathbf{x}_3$. We construct boxes $B(\mathbf{x}_i, 2^\rho)$, $i = 1, 2$, lying inside $B(n)$ with a common radius $2^\rho \asymp 2^{-j_1}n$ that are themselves disjoint and also disjoint from $B_1(\mathbf{x}_3)$. This is possible since $r = 3$ and $h = 2$ and with the given radius for $i = 2$ by (a). Yet by (a) again, the vertex $\mathbf{x}_2$, although in a horseshoe relationship to $\mathbf{x}_1$, may stray as far away as $v_0 2^{-j_2}n$ from $\mathbf{x}_1$ for some constant integer $v_0 \geq 2$ since we are here measuring the distance in terms of the exponent $j_2$. This can easily be dealt with by breaking up the analysis into the subcases $v \leq 4v_0$ and $v > 4v_0$. For $v \leq 4v_0$ we take the outer box $B_2(\mathbf{x}_3)$ of a horseshoe at $\mathbf{x}_3$ to have radius $2^{\nu_3} \asymp 2^{-j_2}n \asymp 2^{-j_1}n$. For $v \leq 4v_0$, we do not yet construct a second horseshoe of outer radius $2^{\nu_3}$. For $v > 4v_0$, we do construct another such horseshoe as follows. We construct an inner horseshoe box $B_1$ that contains $B(\mathbf{x}_i, 2^\rho)$ for each $i = 1, 2$ that has radius $C2^\rho$ and is disjoint also from $B_1(\mathbf{x}_3)$ by our choice of large enough $v$. Accordingly, we adjust the radius $B_2(\mathbf{x}_3)$ upward to $2^{\nu_3} \asymp v2^{-j_2}n \asymp v2^{-j_1}n$ and also define an outer box $B_2$ paired with $B_1$ in a horseshoe formation, by taking the radius of $B_2$ as also $2^{\nu_3} \asymp v2^{-j_1}n$. We choose this radius such that $B_2(\mathbf{x}_3)$ and $B_2$ are disjoint. Hence, if $v > 4v_0$ then we have two horseshoes of equal outer radii. The outer boxes remain in $B(n)$ by (5.28) for $v \leq v_{\max}$. Moreover, these outer boxes are nested inside another box $\widetilde{B}_1$ of radius $C2^{\nu_3}$ whose right edge also lies in $\partial B(n)$. Since we are in case $t = 1$, we may again pair $\widetilde{B}_1$ with an outer horseshoe box $\widetilde{B}_2$ of radius $2^{\nu_1} \asymp 2^{-j^*}n$. By (a), we have that the confinement factor for $\mathbf{x}_2$ is $C2^{-2j_1}n^2$. The confinement factor for $\mathbf{x}_3$ at level $v$ is by (a) and (5.30), $|b_v| \leq C2^{-j_1-j_3}n^2$, independent of $v$. Therefore, by independence and Lemmas 2 and 5, and by using (a) to eliminate the sum on $j_2$, we obtain the following estimation:

$$III_{2,1}a \leq C_\varepsilon n^{9/4+3\varepsilon} \sum_{j_1=0}^{\infty} \sum_{j_1^*=0}^{j_1} \sum_{j_3=j_1}^{\infty} \sum_{v\geq 1} v^{-2+\varepsilon} 2^{j_1^*+(-6j_1-7j_3)/4}. \tag{6.30}$$

Consider next an estimation of $III_{2,1}b$. We have two subcases as follows. Subcase (i): $\mathbf{x}_3$ belongs to the band $b_{1,v_1} := b_{v_1}(\mathbf{x}_2, j_1, j_2, j_1^*)$ for some $v_1 \geq 2$, where $v_1$ ranges up to $v_{1,\max}$ with $v_{1,\max} \asymp 2^{j_1-j_1^*}$. But since we may have $\mathbf{x}_3 \in \bigcup_{v_1=0}^{1} b_{1,v_1}$, we have also subcase (ii): $\mathbf{x}_3$ belongs to the band $b_{2,v_2} := b_{v_2}(\mathbf{x}_2, j_2, j_3, j_1^*)$ for some $v_2 \geq 2$, where now $v_2$ only ranges up to $v_{2,\max}$ with $v_{2,\max} \asymp 2^{j_2-j_1}$. These subcases comprise a dichotomy since $\mathbf{x}_3$ is not in a horseshoe relation to either $\mathbf{x}_2$ or $\mathbf{x}_1$. Partition the sum $III_{2,1}b = III_{2,1}bi + III_{2,1}bii$ accordingly. We study first subcase (ii) under (b). Since $j_2$ is sufficiently larger than $j_1$, we will now be able to construct horseshoes at both $\mathbf{x}_2$ and $\mathbf{x}_3$ with outer radii of each given as $2^{\nu_3} \asymp v_2 2^{-j_2}n$. We also construct a horseshoe with inner box of radius $Cv_2 2^{-j_2}n$ containing both the horseshoes at $\mathbf{x}_2$ and $\mathbf{x}_3$, and with outer radius $2^\rho \asymp 2^{-j_1}n$. We choose $2^\rho$ small enough subject to this asymptotic relation so that the box $B(\mathbf{x}_1, 2^\rho)$ is outside this last horseshoe. We also construct a large horseshoe pair $(\widetilde{B}_1, \widetilde{B}_2)$ with inner and outer radii $2^{\rho_1} \asymp 2^{-j_1} n$ and $2^{\nu_1} \asymp 2^{-j_1^*} n$, respectively, as in the previous cases such that $\widetilde{B}_1$ contains both the smaller nested horseshoes, as well as the box $B(\mathbf{x}_1, 2^\rho)$. The confinement factor for $\mathbf{x}_2$ is $C 2^{-j_1-j_2} n^2$, while that for $\mathbf{x}_3$ under subcase (ii) is $C 2^{-j_2-j_3} n^2$. Therefore, we obtain

$$III_{2,1}bii \le C_\varepsilon n^{9/4+3\varepsilon} \times \sum_{j_1=0}^{\infty} \sum_{j_1^*=0}^{j_1} \sum_{j_2=j_1}^{\infty} \sum_{j_3=j_2}^{\infty} \sum_{v_2 \ge 1} v_2^{-2+\varepsilon} 2^{j_1^* + (-3j_1 - 3j_2 - 7j_3)/4}. \tag{6.31}$$

We turn to subcase (i) under (b). Now $\|\mathbf{x}_3 - \mathbf{x}_i\| \ge C 2^{-j_1} n$, $i = 1, 2$, for all $v_1 \ge 2$. Similar as in case (a), we delay the construction of a horseshoe near $\mathbf{x}_1$ until $v_1 > 4$ since we may have $\|\mathbf{x}_2 - \mathbf{x}_1\| \ge v_1 2^{-j_1+5}$ for some $v_1 \le 2$. But since this will not affect the estimate, we assume $v_1 > 4$. We still construct an inner horseshoe at $\mathbf{x}_2$ with inner horseshoe box $B_1(\mathbf{x}_2)$ of radius $2^\sigma \asymp 2^{-j_2} n$, and outer horseshoe box $B_2(\mathbf{x}_2)$ of radius $2^{\nu_2} \asymp 2^{-j_1} n$. We take $B_2(\mathbf{x}_2)$ to be disjoint from a box $B(\mathbf{x}_1, 2^\rho)$ of radius $2^{\rho_1} \asymp 2^{-j_1} n$. We construct an inner horseshoe box $B_1(\mathbf{x}_3)$ of radius $2^{\rho_3} \asymp 2^{-j_3} n$ and an associated outer horseshoe box $B_2(\mathbf{x}_3)$ of radius $2^{\nu_3} \asymp v_1 2^{-j_1} n$. We also construct a box $B_1$ whose right edge meets the boundary of $B(n)$ that also contains the outer horseshoe box $B_2(\mathbf{x}_2)$ and the box $B(\mathbf{x}_1, 2^\rho)$ that were already constructed to be disjoint. We pair $B_1$ with an outer horseshoe box $B_2$ of radius $2^{\nu_3}$ so that $B_2(\mathbf{x}_3)$ and $B_2$ are disjoint. We finally construct a large horseshoe pair $(\widetilde{B}_1, \widetilde{B}_2)$ but now with a new inner radius $C 2^{\nu_3}$, while the outer radius remains $2^{\nu_1} \asymp 2^{-j_1^*} n$. So we have four horseshoes in all. The confinement factor for $\mathbf{x}_i$ is $C 2^{-j_1-j_i} n^2$, $i = 2, 3$. Therefore, we obtain

$$III_{2,1}bi \le C_\varepsilon n^{9/4+3\varepsilon} \times \sum_{j_1=0}^{\infty} \sum_{j_1^*=0}^{j_1} \sum_{j_2=j_1}^{\infty} \sum_{j_3=j_2}^{\infty} \sum_{v_1 \ge 1} v_1^{-2+\varepsilon} 2^{j_1^* + (j_1 - 7j_2 - 7j_3)/4}. \tag{6.32}$$

We comment on the situation that instead $U_1 = \{\mathbf{x}_3\}$ and $\mathbf{x}_{f_2} = \mathbf{x}_2$ with $U_2 = \varnothing$ when $h = 2$ and $t = 1$. In this case we do not have to consider the possibility that $\mathbf{x}_3 K \mathbf{x}_2$ since, by our definition of the sets $U_a$, if this relation did hold, then we would instead have the case $U_1 = \varnothing$ and $U_2 = \{\mathbf{x}_3\}$ that is analogous to the one we have just considered. We consider now the cases (a) $j_3 \le j_2 + 2c + 8$, and (b) $j_3 > j_2 + 2c + 8$. In case (a) we have an analogous situation as in the previous case (a) except now we have $\mathbf{x}_3 \in b_v(\mathbf{x}_2, j_2, j_3, j_1^*)$, so in the generic case that $v$ is sufficiently large, we construct horseshoes at each of $\mathbf{x}_2$ and $\mathbf{x}_3$ of asymptotically equal inner radii by condition (a) of order $2^{-j_2}n$ and equal outer radii of order $v2^{-j_2}n$. We obtain an estimate for the corresponding sum $III'_{2,1}a$ as

$$III'_{2,1}a \le C_\varepsilon n^{9/4+3\varepsilon} \sum_{j_1=0}^{\infty} \sum_{j_1^*=0}^{j_1} \sum_{j_2=j_1}^{\infty} 2^{j_1^*+(-3j_1-10j_2)/4}.$$

In case (b) we break up the analysis by the dichotomy (i) $\mathbf{x}_3 \in b_{v_2}(\mathbf{x}_2, j_2, j_3, j_1^*)$ for some $v_2 \ge 2$, or (ii) $\mathbf{x}_2 \in b_{v_1}(\mathbf{x}_1, j_1, j_2, j_1^*)$ for some $v_1 \ge 2$. With this minor change in notation, we obtain estimates for the sum corresponding to these cases with the same forms of estimation as shown for $III_{2,1}bi$ and $III_{2,1}bii$. This concludes our discussion of the case $t=1$ under $h=2$.

We next discuss the case $t=2$. We assume, as in the original discussion of $h=2$ and $t=1$, that $U_1=\{\mathbf{x}_2\}$. Now, however, the second root horseshoe index is far from both $\mathbf{x}_2$ and $\mathbf{x}_1$ by assumption, that is, the length scale of this distance is $\max\{2^{-j_3^*}n, 2^{-j_1^*}n\}$ as compared to $2^{-j_1}n$ for the length scale between $\mathbf{x}_2$ and $\mathbf{x}_1$. In the generic case that $j_2$ is sufficiently larger than $j_1$, we construct a horseshoe at $\mathbf{x}_2$ of inner box of radius $2^{\rho_2} \asymp 2^{-j_2}n$ and outer box of radius $2^{\nu_2} \asymp 2^{-j_1}n$ that remains disjoint from a box $B(\mathbf{x}_1, 2^\rho)$ with radius $2^\rho \asymp 2^{-j_1}n$. This horseshoe is nested in a large horseshoe $(\widetilde{B}_1(\mathbf{x}_1), \widetilde{B}_2(\mathbf{x}_1))$ with inner radius $2^{\rho_1} \asymp 2^{-j_1}n$ and outer radius $2^{\nu_1} \asymp 2^{-j_1^*}n$. We also construct a horseshoe $(\widetilde{B}_1(\mathbf{x}_3), \widetilde{B}_2(\mathbf{x}_3))$ with inner box of radius $2^{\rho_3} \asymp 2^{-j_3}n$ and outer box of radius $2^{\nu_3} \asymp 2^{-j_3^*}n$ such that the outer boxes $\widetilde{B}_2(\mathbf{x}_1)$ and $\widetilde{B}_2(\mathbf{x}_3)$ are disjoint. Thus, in this generic case (b) we obtain an estimate

$$III_{2,2}b \le C_\varepsilon n^{9/4+3\varepsilon} \sum_{j_1=0}^{\infty} \sum_{j_1^*=0}^{j_1} \sum_{j_3=j_1}^{\infty} \sum_{j_3^*=0}^{j_3} 2^{j_1^*+j_3^*+(-3j_1-7j_2-7j_3)/4}. \tag{6.33}$$

This concludes our discussion of the case $h=2$.

6.1.6. $\tau=3,\ r=3,\ h=3.$ We finally consider the case $h=3$. Write $III_3$ to denote the sub-sum of $\Sigma_0$ that corresponds to $h=3$. We partition $III_3 = III_{3,1} + III_{3,2} + III_{3,3}$ according to the cases $t=1,2,3$, respectively. Consider first that $h=3$ and $t=1$. We take $\mathbf{x}_3 \in b_v(\mathbf{x}_1, j_1, j_3, j_1^*)$ with some $v \ge 2$. We take $\mathbf{x}_2 \in b_u(\mathbf{x}_1, j_1, j_2, j_1^*)$ with some $u \ge 2$. We have three cases for $v \ge u$:

$$\text{(a) } 2u \ge v \ge u+2, \text{ (b) } v > 2u \quad \text{and} \quad \text{(c) } 0 \le v-u \le 1.$$

Partition $III_{3,1} = III_{3,1}a + III_{3,1}b + III_{3,1}c$ accordingly.

We study first case (a). We define inner horseshoe boxes $B_1(\mathbf{x}_i)$ of radii $2^{\rho_i} \asymp 2^{-j_i}n$, $i=1,2,3$, for the three respective vertices. We define the radii of the associated outer horseshoe boxes $B_2(\mathbf{x}_i)$ for $i=2,3$ to both equal $2^{\nu_2} \asymp (v-u)2^{-j_1}n$. We define the radius of the associated outer horseshoe box $B_2(\mathbf{x}_1)$ to be $2^{\nu_1} \asymp u2^{-j_1}n$. These asymptotic relations are chosen such that all three of $B_2(\mathbf{x}_i)$, $i = 1, 2, 3$, are disjoint. We also construct a horseshoe $(B_1, B_2)$ with inner box $B_1$ of radius $C2^{\nu_2}$ containing both outer boxes $B_2(\mathbf{x}_i)$, $i = 2, 3$, and with outer box $B_2$ of radius $2^{\nu_1}$ and disjoint from $B_2(\mathbf{x}_1)$. We finally construct a large horseshoe pair $(\widetilde{B}_1, \widetilde{B}_2)$ with an inner radius $Cu2^{-j_1}n$ and outer radius $C2^{-j_1^*}n$ such that $\widetilde{B}_1$ contains all the previous outer horseshoe boxes. We have three inner horseshoes, two of which are nested in a fourth horseshoe, and all four of these are nested in a fifth horseshoe. See Figure 4. The confinement factors for $\mathbf{x}_i$ are $2^{-j_1-j_i}$, $i=2,3$. Therefore, after the substitution $\Delta := v-u$, we obtain an estimate

$$III_{3,1}a \le C_\varepsilon n^{9/4+3\varepsilon} \times \sum_{j_1=0}^{\infty}\sum_{j_1^*=0}^{j_1}\sum_{j_2=j_1}^{\infty}\sum_{j_3=j_2}^{\infty}\sum_{u,\Delta\ge 2}(u\Delta)^{-2+\varepsilon}2^{j_1^*+(j_1-7j_2-7j_3)/4}. \tag{6.34}$$

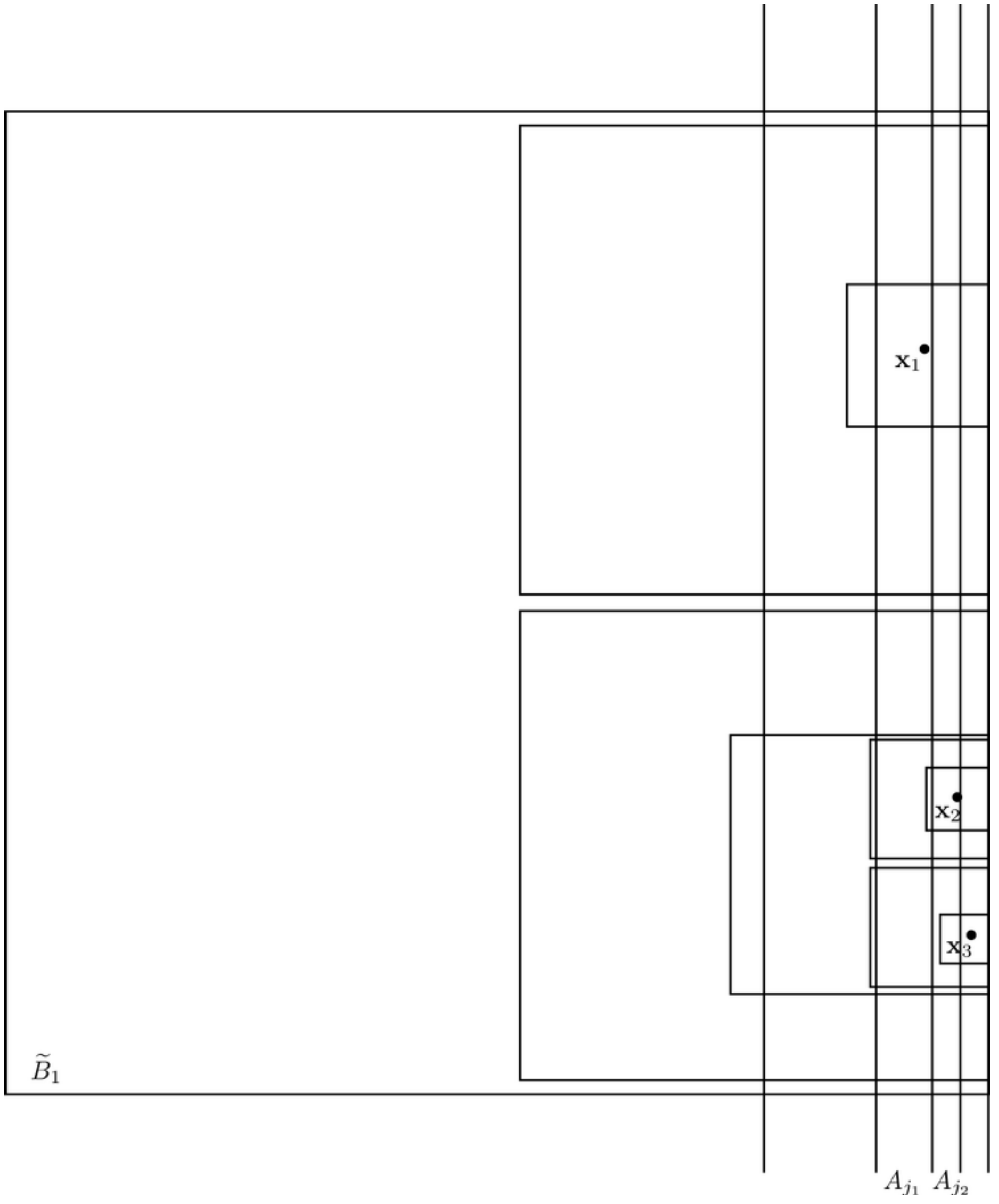

FIG. 4. *Horseshoe construction for the estimate of $III_{3,1}a$.*

In case (b) we change the three inner horseshoe pairs $B_1(\mathbf{x}_i), B_2(\mathbf{x}_i)$, $i=1,2,3$, as follows. We use the same asymptotic formulae $2^{\rho_i} \asymp 2^{-j_i}n$ as in case (a) for the inner radii, and we still write $2^{\nu_1} \asymp u2^{-j_1}n$ and $2^{\nu_2} \asymp (v-u)2^{-j_1}n$, but now define the outer horseshoe radii by changing the vertex at which the larger outer box sits from $\mathbf{x}_1$ to $\mathbf{x}_3$. Indeed, we now take the outer boxes $B_2(\mathbf{x}_i)$, $i=1,2$, to have radii $2^{\nu_1}$, and the outer box $B_2(\mathbf{x}_3)$ to have radius $2^{\nu_2}$ such that all three of these outer boxes are disjoint. Therefore, much as in case (a), by the substitution $\Delta := v-u$, we obtain

$$III_{3,1}b \le C_\varepsilon n^{9/4+3\varepsilon} \times \sum_{j_1=0}^{\infty}\sum_{j_1^*=0}^{j_1}\sum_{j_2=j_1}^{\infty}\sum_{j_3=j_2}^{\infty}\sum_{u\ge 2,\Delta\ge u}(u\Delta)^{-2+\varepsilon}2^{j_1^*+(j_1-7j_2-7j_3)/4}. \tag{6.35}$$

Finally, in case (c) we still take $\mathbf{x}_2 \in b_u(\mathbf{x}_1, j_1, j_2, j_1^*)$ with some $u \ge 2$ but now consider $\mathbf{x}_3 \in b_w(\mathbf{x}_2, j_2, j_3, j_1^*)$ with some $w \ge 2$, where $w$ ranges up to order $2^{j_2-j_1}$. We again construct three inner horseshoe pairs $B_1(\mathbf{x}_i), B_2(\mathbf{x}_i)$, $i=1,2,3$. We again take the corresponding inner radii to be $2^{\rho_i} \asymp 2^{-j_i}n$. We take the radii of the outer boxes $B_2(\mathbf{x}_i)$, $i=2,3$, to be $2^{\nu_3} \asymp w2^{-j_2}n$, so that, by the range of $w$, this is asymptotically no larger than $2^{-j_1}n$. We take the radius of the outer box $B_2(\mathbf{x}_1)$ to be $2^{\nu_1} \asymp u2^{-j_1}n$. We choose the boxes subject to these asymptotic formulae such that all three outer boxes are disjoint. We also construct a fourth horseshoe $(B_1, B_2)$ with inner box $B_1$ of radius $C2^{\nu_3}$ containing both outer boxes $B_2(\mathbf{x}_i)$, $i=2,3$, and with outer box $B_2$ of radius $2^{\nu_1}$ that is disjoint from $B_2(\mathbf{x}_1)$. We finally construct a large horseshoe pair $(\widetilde{B}_1, \widetilde{B}_2)$ with an inner radius $Cu2^{-j_1}n$ and outer radius $C2^{-j_1^*}n$ such that $\widetilde{B}_1$ contains all the previous outer horseshoe boxes. Therefore, since the confinement factor for $\mathbf{x}_3$ is now instead $2^{-j_2-j_3}$,

$$III_{3,1}c \le C_\varepsilon n^{9/4+3\varepsilon} \times \sum_{j_1=0}^{\infty}\sum_{j_1^*=0}^{j_1}\sum_{j_2=j_1}^{\infty}\sum_{j_3=j_2}^{\infty}\sum_{u\ge 2,w\ge 2}(uw)^{-2+\varepsilon}2^{j_1^*+(-3j_1-3j_2-7j_3)/4}. \tag{6.36}$$

If instead we consider $v<u$ in the original setting, we apply the same method but with the roles of $u$ and $v$ switched. This completes our discussion of the case $t=1$.

Consider next that $t=2$. Say that $j_3^*$ is the second root dual horseshoe index. We let $\mathbf{x}_2 \in b_u(\mathbf{x}_1, j_1, j_2, j_1^*)$ for some $u \geq 2$. We construct three horseshoe pairs $(B_1(\mathbf{x}_i), B_2(\mathbf{x}_i))$, $i=1,2,3$, with inner radii $2^{\rho_i} \asymp 2^{-j_i}n$. We take the radii of $B_2(\mathbf{x}_i)$, $i=1,2$, to be $2^{\nu_1} \asymp u2^{-j_1}n$, but because $\mathbf{x}_3$ is now far from both $\mathbf{x}_1$ and $\mathbf{x}_2$, we construct the radius of the outer box $B_2(\mathbf{x}_3)$ to be $2^{\nu_3} \asymp 2^{-j_3^*}n$. As before, we construct these three outer boxes to be disjoint. Finally, we construct a fourth horseshoe pair $(B_1, B_2)$ with inner box $B_1$ of radius $C2^{\nu_1}$ containing both outer boxes $B_2(\mathbf{x}_i)$, $i=1,2$, and with outer box $B_2$ of radius $2^{-j_1^*}n$ that is disjoint from $B_2(\mathbf{x}_3)$. This is possible due to the separation of the dual indices, where $u$ ranges up to order $2^{j_1^*-j_1}$. Therefore, since $|b_u| \leq 2^{-j_1-j_2}$, we have by (5.15) that

$$
\begin{aligned}
&III_{3,2} \leq C_\varepsilon n^{9/4+3\varepsilon} \\
&\quad \times \sum_{j_1=0}^{\infty} \sum_{j_1^*=0}^{j_1} \sum_{j_2=j_1}^{\infty} \sum_{j_3=j_2}^{\infty} \sum_{j_3^*=0}^{j_3} \sum_{u\geq 2} (u)^{-2+\varepsilon} 2^{j_1^*+j_3^*+(-3j_1-7j_2-7j_3)/4}.
\end{aligned}
\tag{6.37}
$$

We consider finally the case $t=3$. This is the easiest case. It refers to widely separated vertices. We construct three horseshoes $B_1(\mathbf{x}_i), B_2(\mathbf{x}_i)$, $i=1,2,3$, with inner radii $2^{\rho_i} \asymp 2^{-j_i}n$ and outer radii $2^{\nu_i} \asymp 2^{-j_i^*}n$, $i=1,2,3$, respectively. Therefore, by (5.15), we obtain

$$
\begin{aligned}
&III_{3,3} \leq C_\varepsilon n^{9/4+3\varepsilon} \\
&\quad \times \sum_{j_1=0}^{\infty} \sum_{j_1^*=0}^{j_1} \sum_{j_2=j_1}^{\infty} \sum_{j_2^*=0}^{j_2} \sum_{j_3=j_2}^{\infty} \sum_{j_3^*=0}^{j_3} 2^{j_1^*+j_2^*+j_3^*+(-7j_1-7j_2-7j_3)/4}.
\end{aligned}
\tag{6.38}
$$

This completes our discussion of the case $h=3$ when $r=3$.

**7. The general case.** We show in this section a two-fold argument for establshing the general $\tau$th moment for the number of pivotal sites. The first part of the argument is to isolate discussion of the (nonroot) vertices that are chained to a given root. This is accomplished by utilizing Lemma 7 and its generalization in Lemma 8, together with Proposition 1. Lemma 8 is required to handle the case $\tau \geq 5$. This part of the analysis does not require any horseshoe estimates. The second part of the argument is to explain the general strategy for the construction of horseshoes at root vertices, as well as the construction of nested horseshoes.

7.1. *Nonroot vertices.* We begin with an example to understand how to generalize the argument of Section 6.1.1. Let $\tau=5$ and $r=1$. Assume $\mathbf{G}_1 = \{\mathbf{x}_1, \mathbf{w}_1, \mathbf{w}_{1,1}, \mathbf{w}_{1,2}, \mathbf{w}_{1,2,1}\}$. Consider the following dichotomy:

$$\text{either (i) } m_{1,2} \geq m_1 + s_1 + 2 \quad \text{or} \quad \text{(ii) } m_{1,2} < m_1 + s_1 + 2.$$

Define boxes $R_1 := B(\mathbf{w}_1, l_1)$, for $l_1 := 2^{-m_1-s_1}n$ and $R_{1,2} := B(\mathbf{w}_{1,2}, l_{1,2})$, for $l_{1,2} := 2^{-m_{1,2}-s_1}n$, and more generally, $R_{i_1,\ldots,i_k} := B(\mathbf{w}_{i_1,\ldots,i_k}, l_{i_1,\ldots,i_k})$, for $l_{i_1,\ldots,i_k} := 2^{-m_{i_1,\ldots,i_k}-s_1}n$, where $s_1$ is the constant $s_1 := s + s_0$ for $s = 2c + 4$ (see Section 4) and some $s_0 > 0$ to be determined below.

We study first the generic case (i). Since $\mathbf{w}_{1,2} \in a_{m_{1,2}}(\mathbf{w}_1)$, we have that $R_{1,2} \subset B(\mathbf{w}_{1,2}, 2^{-m_{1,2}-2}n) \subset R_1$ while also $\mathbf{w}_1 \notin B(\mathbf{w}_{1,2}, 2^{-m_{1,2}-2}n)$, so, in particular, $\mathbf{w}_1 \notin R_{1,2}$. Consider now the following subcases:

(a1) $m_{1,1} \geq m_1 + s_1 + 2$ or (b1) $m_{1,1} < m_1 + s_1 + 2$,
(a2) $m_{1,2,1} \geq m_{1,2} + s_1 + 2$ or (b2) $m_{1,2,1} < m_{1,2} + s_1 + 2$.

Consider first the generic joint subcase (a1)–(a2) under case (i). By the very same reasoning as given for (i), by (a1), we have that $B(\mathbf{w}_{1,1}, 2^{-m_{1,1}-2}n) \subset R_1$ while also $\mathbf{w}_1 \notin B(\mathbf{w}_{1,1}, 2^{-m_{1,1}-2}n)$, so, in particular, $\mathbf{w}_1 \notin R_{1,1}$. By the same reasoning again under (a2), $R_{1,2,1} \subset B(\mathbf{w}_{1,2,1}, 2^{-m_{1,2,1}-2}n) \subset R_{1,2}$ while also $\mathbf{w}_{1,2} \notin B(\mathbf{w}_{1,2,1}, 2^{-m_{1,2,1}-2}n)$, so, in particular, $\mathbf{w}_{1,2} \notin R_{1,2,1}$. Furthermore, by Proposition 1, we have that $R_{1,1}$ and $R_{1,2}$ are disjoint. Hence, we have in all the following picture: $R_{1,2,1} \subset R_{1,2} \subset R_1$ and $R_{1,1} \subset R_1$ with $R_{1,1} \cap R_{1,2} = \varnothing$. See Figure 5. We note that the context of Lemma 7 may be generalized to the present circumstance as follows.

LEMMA 8. *Let $R = R(\mathbf{x}')$ be a rectangle centered at $\mathbf{x}'$ with its shortest half-side of length $l \geq 1$ and longest half-side of length $L \geq 1$ such that $1 \leq L/l \leq 2$. Let $R$ contain a vertex $\mathbf{x}$ such that $\|\mathbf{x} - \mathbf{x}'\| \leq l/2$. Suppose that $R$ contains a collection of boxes $B_i = B(\mathbf{y}_i, 2^{\lambda_i})$, $i = 1, \ldots, v$, such that for every $i$, $\mathbf{x} \notin B_i$. Assume that any two of the boxes $B_i$ are either disjoint or one of them is contained entirely within the other and that whenever box $B_i \subset B_j$, the smaller box $B_i$ does not contain the center $\mathbf{y}_j$ of the larger box $B_j$. Denote $V := \{\mathbf{y}_i : i = 1, \ldots, v\}$ and $D_V = D_V(\mathbf{x}) := \max_{\mathbf{y}_i \in V} \|\mathbf{y}_i - \mathbf{x}\|$ and assume that $\sum_{\mathbf{y}_i \in V} 2^{\lambda_i} \leq \frac{1}{64} D_V(\mathbf{x})$. For any $W \subset V$ such that $B_i \subset B_{i_0}$ for all $\mathbf{y}_i \in W$ and some $\mathbf{y}_{i_0} \in V \setminus W$, denote also the maximal distance $D_W(\mathbf{y}_{i_0}) := \max_{\mathbf{y}_i \in W} \|\mathbf{y}_i - \mathbf{y}_{i_0}\|$, and assume that $\sum_{\mathbf{y}_i \in W} 2^{\lambda_i} \leq \frac{1}{64} D_W(\mathbf{y}_{i_0})$. Then the conclusion of Lemma 7 continues to hold.*

The proof of Lemma 8 follows by induction on $v$ and is an extension of the proof of Lemma 7. To show first how to apply Lemma 8 for our construction, we work on the current example. We begin with the assumption that we are in the joint subcase (a1)–(a2) of (i). Take $R = R_1$, $\mathbf{x} = \mathbf{w}_1$ and $v = 3$, where the centers $\mathbf{y}_i$, $i = 1, 2, 3$, are the vertices $V := \{\mathbf{w}_{1,1}, \mathbf{w}_{1,2}, \mathbf{w}_{1,2,1}\}$. The boxes $B_i$ are the corresponding boxes $R_{i_1,\ldots,i_k}$ with radii $2^{\lambda_i} = l_{i_1,\ldots,i_k}$ as above. We assume that $s_0$ is so large in the definition of the radii of these boxes through the parameter $s_1 = s + s_0$ that we do not need to shrink these boxes to establish the sum of diameters conditions. This is possible since we do not change the positions of the vertices of $\mathbf{G}_1$ (we do not change $s = 2c + 4$), but only adjust the increment $s_0$ so that the radii of the boxes centered at these vertices change. In particular, for $W = V$, we have that the sum of all the diameters of the boxes is at most $2\tau 2^{-m_{1,2}-s-s_0}n$ since both $m_{1,2,1} \geq m_{1,2}$ and $m_{1,1} \geq m_{1,2}$ while $D_V(\mathbf{w}_1) \geq 2^{-m_{1,2}-2}n$, so it suffices for this case to find $s_0$ such that $8\tau 2^{-s-s_0} < 1$. If instead $W = \{\mathbf{w}_{1,1}\}$, then $D_W(\mathbf{w}_1) \geq 2^{-m_{1,1}-2}n$ while the radius of the box $R_{1,1}$ is $2^{-m_{1,1}-s-s_0}n$. For the case $W = \{\mathbf{w}_{1,2,1}\}$, we subtract to find

$$\|\mathbf{w}_{1,2,1} - \mathbf{w}_1\| \geq \|\mathbf{w}_{1,2} - \mathbf{w}_1\| - \|\mathbf{w}_{1,2} - \mathbf{w}_{1,2,1}\| \geq 2^{-m_{1,2}-2} - 2^{-m_{1,2}-c} \geq 2^{-m_{1,2}-3}$$

by our choice of $c$ in Section 4. Hence, again $D_W(\mathbf{w}_1) \geq C2^{-m_{1,2}} \geq 642^{-m_{1,2,1}-s_1}$. In summary, we can see that the sum of diameters condition for the maximal distance $D_W(\mathbf{w}_1)$ is satisfied because, first, the quantity $D_W(\mathbf{w}_1)$ is estimated below by a constant times the distance $\|\mathbf{w}_1 - \mathbf{w}_{1,j}\|$ (and this distance is at least $2^{-m_{1,j}-2}$) for the maximal $j \geq 1$ such that some vertex $\mathbf{w}_{1,j,i_3,\ldots,i_k} \in W$. Second, the radius of each box whose center is in $W$ is at most $2^{-m_{1,j}-s-s_0}$. We shall construct below a nested or disjoint boxes condition for each joint subcase referred to above such that the boxes $B_i$ that will be contained in a given box $B_{i_0}$ with center $\mathbf{y}_{i_0} = \mathbf{w}_0$ (or a rectangular variant of $B_{i_0}$) will correspond to centers that are children of $\mathbf{w}_0$. Therefore, we see by the graphical structure of Section 4 that the argument above for verifying the sum of diameters condition is independent of the joint subcase. We fix the value of $s_1$ in the definitions of the radii $l_{i_1,\ldots,i_k}$ and in the conditions (i), (ii), (a1), (a2), (b1) and (b2). Thus, the disjoint and nested relations of our joint subcase (a1)–(a2) of (i) are preserved as in Figure 5 and so we have verified the hypothesis of Lemma 8 for this case. We next verify the conclusion of Lemma 8 for this case as well.

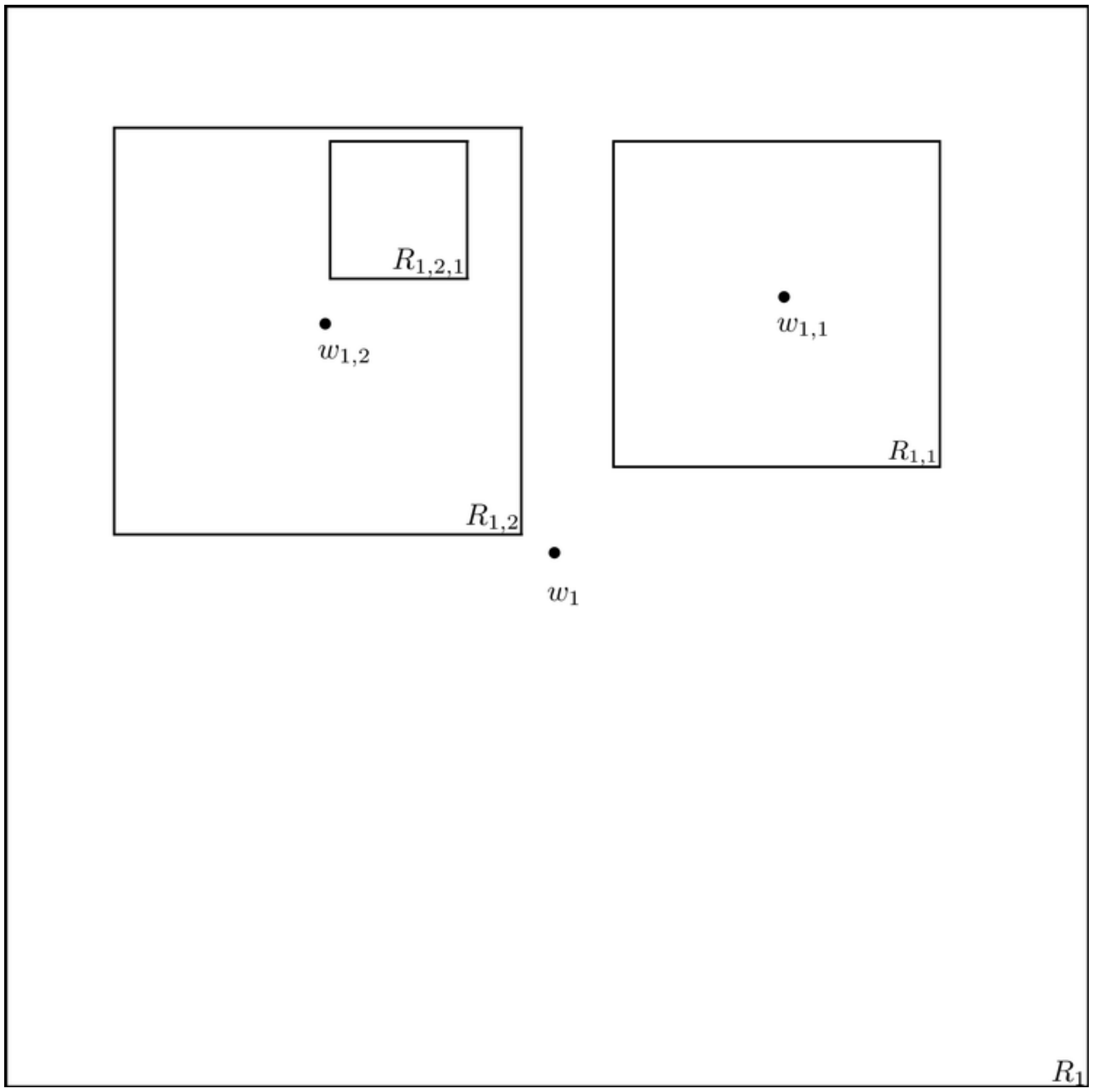

FIG. 5. *Configuration of boxes $R_{i_1,\dots,i_k}$ in joint subcase* (a1)–(a2) *under* (i).

The induction hypothesis is the conclusion of Lemma 8 for some number of boxes $B_i$, $i = 1, \ldots, u$, with $u$ in place of $v$. The proof of the induction hypothesis for $v = 1$ in Lemma 8 is the same as the proof for $v = 1$ in the original statement of Lemma 7 with $2^{\lambda_1}$ in place of $2^{\lambda_1 - 1}$ in (A.1). We establish the inductive step of Lemma 8 for the current example. We construct a rectangle $\widetilde{R}_1$, as in the proof of Lemma 7, that has a center $\widetilde{\mathbf{x}}_1'$ such that $\|\widetilde{\mathbf{x}}_1' - \mathbf{w}_1\| \leq \frac{1}{20} D_V \leq \widetilde{l}_1/2$ for $\frac{1}{5} D_V \leq \widetilde{l}_1 \leq \frac{1}{10} D_V$ and $D_V \leq l_1$. Let us assign indices by $B_1 = R_{1,1}$, $B_2 = R_{1,2}$ and $B_3 = R_{1,2,1}$, so that $B_3$ does not contain the center of $B_2$ and $B_3 \subset B_2$. By the proof of Lemma 7, all the boxes $B_i$ are either entirely inside or outside $\widetilde{R}_1$. Therefore, if one of $B_1$ or $B_2$ does lie inside $\widetilde{R}_1$, then because at least one other does lie outside $\widetilde{R}_1$, we can apply the induction hypothesis applied to the rectangle $\widetilde{R}_1$ in place of $R$. So let us assume that there are no boxes $B_i$ inside $\widetilde{R}_1$. In this case we proceed to construct a second rectangle $\widetilde{R}_2 = B(\mathbf{w}_{1,2}, r)$ inside $R_{1,2}$, as in the proof of the initial case of Lemma 7 applied with $R = R_{1,2}$ and with $r = \|\mathbf{w}_{1,2,1} - \mathbf{w}_{1,2}\|/4$, so that $B_3$ is disjoint from $\widetilde{R}_2$. That is, we apply the induction step only to the context of a single box $B_3$ inside the box $R = R_{1,2}$. Explicitly, we have that

$$\begin{aligned} &P\left(\mathcal{U}_4(\mathbf{w}_1; R) \cap \left(\bigcap_{i=1}^{3} \mathcal{U}_4(\mathbf{y}_i; R)\right)\right) \\ &\qquad \leq (\mathcal{U}_4(\mathbf{w}_1; \widetilde{R}_1)) P(\mathcal{U}_4(\mathbf{w}_{1,1}; R_{1,1})) \\ &\qquad\quad \times P(\mathcal{U}_4(\mathbf{0}, 2^{\lambda_3})) P(\mathcal{U}_4(\mathbf{w}_{1,2}; \widetilde{R}_2)) P(\mathcal{U}_4(\mathbf{0}, 5r; l_{1,2}/2)), \end{aligned} \tag{7.1}$$

where $2^{\lambda_3} = 2^{m_{1,2,1} - s_1} n$ and where if $24r \geq l_{1,2}$, we omit the last factor in this inequality. If $24r < l_{1,2}$, then the product of the last two factors in this inequality is bounded by $CP(\mathcal{T}_4(\mathbf{0}; l_{1,2}/d))$, as in the proof of the case $v = 1$ of Lemma 7. Therefore, we obtain the desired conclusion of Lemma 8 for the particular example.

Consider next subcase (a1)–(b2) under (i). We still have the boxes $B_i$, $i = 1,2,3$, as defined above, but now they may not satisfy the nested or disjoint condition of Lemma 8. But, by (b2), $m_{1,2}+s_1+s+2 > m_{1,2,1}+s$. Therefore, the boxes $B(\mathbf{w}_{1,2}, 2^{-m_{1,2}-s_1-s-2}n)$ and $B(\mathbf{w}_{1,2,1}, 2^{-m_{1,2,1}-s}n)$ are disjoint by Proposition 1. Therefore, since this last box contains $B_3$ (because $s_1 > s$), in fact, we may shrink the box $B_2$ by the constant factor $2^{-s}$ to obtain a box $B_2'$ such that now the boxes $B_1$, $B_2'$, $B_3$ are mutually disjoint. Also, due to the geometric series estimate $\|\mathbf{w}_{1,2,1}-\mathbf{w}_1\| \le 2^{-m_{1,2}-1}n$, we have by (i) that $B_3 \subset R_1$. Hence, in fact, we may apply Lemma 7 with the $B_2'$ in place of the $B_2$ and with $B_1$ and $B_3$ as before and still with $R = R_1$ to again reach the desired conclusion. This trick must be modified in general. In subcase (b1)–(a2) under (i) we have that the boxes $B(\mathbf{w}_1, 2^{-m_1-s_1-s-2}n)$ and $B_1 = B(\mathbf{w}_{1,1}, 2^{-m_{1,1}-s_1}n)$ are disjoint by (b1) and Proposition 1, and, moreover, $B_3$ is nested in $B_2$ by (a2). Further, by Proposition 1 applied to the trimmed graph $\mathbf{G}_1$ without the vertex $\mathbf{w}_{1,2,1}$, we have that $B_1$ and $B_2$ are disjoint. To obtain a nested or disjoint condition, we apply the method of Proposition 2 with $\mathbf{w}_1$ playing the role of the root and $B_2 = B(\mathbf{w}_{1,2}, 2^{-m_{1,2}-s_1})$ playing the role of a shrunken box, and with $D = 2^{-s-2}l_1 = 2^{-m_1-s_1-s-2}$. Since $m_{1,2} \ge m_1$, the diameter of the box $B_2$ is small compared with $D$, depending on the parameter $s_0$, so we may adjust this parameter upward if necessary to construct a rectangle $R_1'$ with center $\mathbf{x}'$ that satisfies $\|\mathbf{x}'-\mathbf{w}_1\| \le l_1'/2$, where $R_1'$ has half sides $l_1'$ and $L_1'$ with $1 \le L_1'/l_1' \le 2$ such that both $B_2$ and $B_3$ lie either nested together inside $R_1'$ or nested together outside $R_1'$. Here in the method of proof of Lemma 7 we take $D/10 \le l_1' \le D/5$ so that $R_1' \subset B(\mathbf{w}_1, l') \subset B(\mathbf{w}_1, 2^{-s-2}l_1)$. Thus, we obtain that $B_1$ lies outside $R_1'$. By working again with the method of Proposition 2, this time with $D = C2^{-j_1}n$, for any of the joint subcases, we may enclose $\mathbf{x}_1$ by a suitable rectangle $R_0'$ with shortest half-side $l \asymp 2^{-j_1}n$ such that $R_0'$ satisfies a nested or disjoint condition with the other rectangles constructed thus far. Now Lemma 8 is applied with $\mathbf{x} = \mathbf{x}_1$ and the rectangles contained by $R_0'$ if any such exist. If $R_1'$ is disjoint from $R_0'$, then we apply Lemma 8 separately with $\mathbf{x} = \mathbf{w}_1$ and the rectangles contained by $R_1'$. Thus, in each joint subcase organized as above, there will exist a specific arrangement of boxes and rectangles with radii given asymptotically by $2^{-m_{i_1,\dots,i_k}}n$ such that Lemma 8 will apply to disjoint pieces of the arrangement. As shown in Section 6.1.1, we may alternatively choose to apply Lemma 8 first with $R = R_1'$ and use that $B(\mathbf{x}_1, 2^{-m_1-s}n)$ is disjoint from $R_1$ and all other boxes by Proposition 1. However, it is convenient to use the more general format with $\mathbf{x} = \mathbf{x}_1$ and $R = R_0'$ for one application of Lemma 8 to obtain the same result in all joint subcases, namely, by Lemma 2, that

$$p_{n,3}(\mathbf{x}_1,\dots,\mathbf{w}_{1,2,1}) \le C_\varepsilon n^{-25/4+5\varepsilon} 2^{5(j_1+m_1+m_{1,1}+m_{1,2}+m_{1,2,1})/4}. \tag{7.2}$$

Hence, by using the confinement factors (3.12), by (7.2) and by constructing one horseshoe at $\mathbf{x}_1$, we have that any sub-sum of $\Sigma_0$ corresponding to $\tau = 5$ and $r = 1$ under the eight different joint subcase combinations is bounded by

$$(7.3)\quad C_\varepsilon n^{15/4+5\varepsilon} \sum_{j_1=0}^{\infty} \cdots \sum_{m_{1,2,1}=m_{1,2}}^{\infty} 2^{j_1^* + (-7j_1 - 3m_1 - 3m_{1,2} - 3m_{1,1} - 3m_{1,2,1})/4},$$

where the intervening summations are indexed by the conditions $0 \le j_1^* \le j_1$, $m_1 \ge j_1$, $m_{1,2} \ge m_1$ and $m_{1,1} \ge m_{1,2}$. In conclusion, since the confinement factor for each vertex $\mathbf{w}_{i_1,\ldots,i_k} \in V_1$ is of the form $2^{-2m_{i_1,\ldots,i_k}} n^2$, we are able to establish a convergent sum analogous to (7.3) for any graph $\mathbf{G}_1$ since we have shown that we can take $2^{\lambda_i} \asymp 2^{-m_{i_1,\ldots,i_k}} n$ corresponding to the vertex $\mathbf{w}_{i_1,\ldots,i_k}$ in Lemma 8.

It remains to make some comments about the case when there is a second graph $\mathbf{G}_2$ with root $\mathbf{x}_e$. We argued briefly in the lines preceding (6.7) of Section 6.1.2 that, even with $h = 1$, there would exist a box $B = B(\mathbf{x}_e, l)$ with radius $l \asymp 2^{-j_e} n$ such that $B$ is disjoint from all the rectangles constructed as above in Section 7.1 with centers or mock centers at vertices in $\mathbf{G}_1$. The argument is based on the fact that if, in the construction of the vertices $V_1$ for $\mathbf{G}_1$, we obtain $m_1 = j_1 + 2c$, then of course since $j_e + 4c + 2 \ge j_1 + 4c + 2 \ge m_1 + 2c + 2$, the boxes mentioned in the line preceding (6.7) with $\mathbf{x}_e$ in place of $\mathbf{x}_3$ are disjoint. If, on the other hand, $m_1 = j_1 + 2c + 1$, then, for all $\mathbf{w}_{i_1,\ldots,i_k} \in W_1$ with $k \ge 2$, we have

$$(7.4)\qquad \|\mathbf{w} - \mathbf{x}_1\| \le 2^{-j_1 - 2c - 1} + 2^{-m_1 - c} \le 3 \cdot 2^{-j_1 - 2c - 2}.$$

Hence, again $B(\mathbf{x}_e, 2^{-j_e - 4c - 2} n)$ is disjoint from $B(\mathbf{w}_{i_1,\ldots,i_k}, 2^{-m_{i_1,\ldots,i_k} - 2c - 2} n)$. This last statement continues to hold with $k \ge 1$ and $i_1 \ge 2$ since we may similarly argue for $i \ge 2$ that $j_1 + c \le m_i \le j_1 + 2c$ or $m_i \ge j_1 + 2c + 1$. Therefore, we may extend the analysis of (7.3) to the case of a second graph $\mathbf{G}_2$ with an appropriate horseshoe at the second root vertex $\mathbf{x}_e$. In conclusion of this section, we have shown a convergent dyadic sum process up to the construction of horseshoes at the root vertices.

7.2. *Horseshoes.* The remainder of the general argument is based on a pattern of nonoverlapping horseshoes and the corresponding probability estimates that provide for convergence factors in the pivotal case. We have shown in Section 6 a parametrization of spacing between the root vertices provided by the definitions of the root vertices, the root horseshoe vertices, certain bands dividing space between nearby roots and the separation of dual horseshoe indices. These levels of organization determine confinement factors associated to each root vertex in terms of the dyadic indices $j_i$ and dual indices $j_i^*$. In fact, there exists a sufficiently large constant $C_0$ in this parametrization such that as long as a root vertex $\mathbf{x}_k$ is not confined to belong to a set of vertices of size at most $C_0 2^{-2j_k} n^2$, then a horseshoe is constructed at this vertex. We have, in particular, that a horseshoe is constructed at each root vertex $\mathbf{x}_k$ such that its dyadic index $j_k$ does not satisfy $j_i \le j_k \le j_i + c_0$ for some root vertex $\mathbf{x}_i$ with $i < k$ and dyadic index $j_i$, where $c_0$ is a constant positive integer. If indeed the condition $j_i \le j_k \le j_i + c_0$ does hold and if in addition the vertex $\mathbf{x}_k$ is in a horseshoe relation to $\mathbf{x}_i$, then the vertex $\mathbf{x}_k$ is confined to belong to a set of vertices of size at most $C 2^{-2j_i} n^2$. In this case a convergence factor for this vertex is accounted for by its confinement factor alone. This follows because the probability that a four-arm path issues from the center of a box of radius $2^\rho \asymp 2^{-j_k} n$ and then exits this box is at most $C_\varepsilon 2^{5j_k/4} n^{-5/4+\varepsilon}$. Thus, by multiplying the confinement factor by the probability and by substitution of $j_i$ for $j_k$ due to the condition on these indices, we obtain the convergence factor $2^{-3j_i/4}$. This situation is an exception wherein a horseshoe is not constructed for lack of space. It represents an analogue of the confinement that is associated to each nonroot vertex. Note, on the other hand, that additional horseshoes besides those at roots may need to be constructed in general to fill in spaces between the dyadic annuli in $B(n)$. Such is the case because horseshoes are not allowed to overlap. Therefore, if one vertex belongs to a band associated with another, then the horseshoes associated to each can only grow so large (with equal outer radii). Then a larger horseshoe containing both the horseshoes that have grown together must be constructed as though the two vertices had become one (cf. Figure 4). We can summarize this strategy by observing that a maximal number of horseshoes is introduced for a given parametrization of spacing of root vertices. Due to the nesting of two horseshoes with equal outer radii inside a single larger horseshoe, the algebraic factors (e.g. $v^{-2+\varepsilon}$) associated to the confinement of vertices in bands remain always with exponents $(-2+\varepsilon)$, so contribute only convergent terms in our method. Due to this allowance for nesting of horseshoes, the additional dyadic convergence factors that arise from Lemma 5 compensate in exactly the same way a confinement factor would if there were to be no room for a horseshoe.

In conclusion, each vertex in our estimation method for the pivotal sites, be it a root or nonroot, contributes a convergence factor with the same exponent $-\frac{3}{4}$. Thus, each arrangement of the vertices that defines a sub-sum $J$ of $\Sigma_0$ in (4.1)–(4.2) by the above division of cases yields $J \le C_\varepsilon n^{3\tau/4+\tau\varepsilon} \sum_{j_1=0}^{\infty} 2^{-3\tau j_1/4}$. By contrast, when we apply our method to the case of items 1 or 2 in Theorem 1, we omit the construction of horseshoes altogether. We also omit the need for Lemmas 7 and 8. Then by Proposition 1 and Lemma 2 alone, the root vertices and nonroot vertices contribute convergence factors with exponents $-\frac{1}{3}$ and $-\frac{4}{3}$, respectively, in the case of the lowest crossing and exponents $-\frac{3}{4}$ and $-\frac{7}{4}$, respectively, in the case of pioneering sites.

## APPENDIX

Proof of Lemma 7. We proceed by induction on the number of boxes $v$. We establish first the statement of the lemma for $v = 1$. Define $r := \|\mathbf{x} - \mathbf{y}_1\|/4$. We have that

$$B(\mathbf{x}, r) \text{ is disjoint from } B(\mathbf{y}_1, 2^{\lambda_1 - 1}). \tag{A.1}$$

Indeed, $2^{\lambda_1 - 1} < \|\mathbf{x} - \mathbf{y}_1\|/2 = 2r$. Note that $L \geq \|\mathbf{x} - \mathbf{y}_1\| = 4r$, so $l \geq L/2 \geq 2r$. Therefore,

$$B(\mathbf{x}, r) \subset B(\mathbf{x}, l/2) \subset R, \tag{A.2}$$

where the last inclusion follows because $\|\mathbf{x} - \mathbf{x}'\| \leq l/2$ and $B(\mathbf{x}', l) \subset R$. Now we use (A.1), (A.2), the assumption $B(\mathbf{y}_1, 2^{\lambda_1}) \subset R$ and independence to obtain by (2.5) and (5.20) that

$$P(\mathcal{U}_4(\mathbf{x}; R) \cap \mathcal{U}_4(\mathbf{y}_1; R)) \leq P(\mathcal{U}_4(\mathbf{0}, r))P(\mathcal{U}_4(\mathbf{0}, 2^{\lambda_1 - 1})). \tag{A.3}$$

We consider now two cases, $l \leq 24r$ and $l > 24r$. If $l \leq 24r$, then by Lemma 4 with $\kappa = 4$, we have

$$P(\mathcal{U}_4(\mathbf{0}, r)) \leq P(\mathcal{U}_4(\mathbf{0}, l/24)) \leq C_4 P(\mathcal{T}_4(\mathbf{0}, l/24))$$

so we are done by (A.3) in this case. If instead $l > 24r$, then, since $\|\mathbf{x} - \mathbf{y}_1\| = 4r$ implies that

$$B(\mathbf{x}, r), B(\mathbf{y}_1, r) \subset B(\mathbf{x}, 5r)$$

and since also by (A.2) the annulus $\mathbf{x} + A(\mathbf{0}, 5r; l/2) \subset R$, we have that

$$\begin{aligned} &P(\mathcal{U}_4(\mathbf{x}; R) \cap \mathcal{U}_4(\mathbf{y}_1; R)) \\ &\quad \leq P(\mathcal{U}_4(\mathbf{0}, r))P(\mathcal{U}_4(\mathbf{0}, 2^{\lambda_1 - 1}))P(\mathcal{U}_4(\mathbf{0}, 5r; l/2)). \end{aligned} \tag{A.4}$$

We now apply Lemmas 4 and 6 and construct connections across the annulus $\mathbf{x} + A(\mathbf{0}, r; 5r)$ to show first that

$$\begin{aligned} P(\mathcal{U}_4(\mathbf{0}, r))P(\mathcal{U}_4(\mathbf{0}, 5r; l/2)) &\leq C_4 P(\mathcal{T}_4(\mathbf{0}, r))P(\mathcal{U}_4(\mathbf{0}, 5r; l/2)) \\ &\leq C P(\mathcal{U}_4(\mathbf{0}, l/2)). \end{aligned}$$

Then we apply Lemma 4 again to establish that this last probability is at most

$$C' P(\mathcal{T}_4(\mathbf{0}, l/2)).$$

Therefore, by (A.4) and these last two observations, we have established the lemma for $v = 1$ in the case $l > 24r$. Thus, by comparing the two cases, we can take $d = 24$ and $c_1 = 1$ in the statement of the lemma when $v = 1$.

We now proceed to show the inductive step. Assume the statement of the lemma is true with a positive integer $u$ in place of $v$ for some $u < v$ and with $v \geq 2$. Define

$$D_v := \max_{i=1,\ldots,v} \|\mathbf{y}_i - \mathbf{x}\|.$$

Note that $l \geq L/2 \geq D_v/2$, so that $B(\mathbf{x}, \frac{1}{4}D_v) \subset B(\mathbf{x}, l/2) \subset R$. Consider now the shrunken boxes $B_i' := B(\mathbf{y}_i, 2^{\lambda_i - 3v})$, $i = 1, \ldots, v$. Since $D_v > 2^{\lambda_i}$ for all $i$, we have that the sum of the diameters of these boxes satisfies

$$2\sum_{i=1}^{v} 2^{\lambda_i - 3v} \leq 2vD_v 2^{-3v} \leq \tfrac{1}{16}D_v$$

for all $v \geq 2$. Therefore, even if all the shrunken boxes were packed inside the subset $B(\mathbf{x}, \frac{1}{4}D_v)$ of the rectangle $R$, there would be a gap in the $z_1$ coordinates of the vertices $\mathbf{z} \in B(\mathbf{x}, \frac{1}{4}D_v)$ somewhere in the interval $[x_1 + \frac{1}{10}D_v, x_1 + \frac{1}{5}D_v]$ and also in the interval $[x_1 - \frac{1}{5}D_v, x_1 - \frac{1}{10}D_v]$, where we denote $\mathbf{x} = (x_1, x_2)$. Indeed, each of these intervals has width $\frac{1}{10}D_v$ which is strictly greater than the sum of the diameters of the boxes $B_i'$. Similarly, there must be gaps in the $z_2$ coordinates of the vertices $\mathbf{z} \in B(\mathbf{x}, \frac{1}{4}D_v)$ somewhere in corresponding intervals for the second coordinate. Therefore, by constructing a rectangle with sides along some vertical and horizontal lines through the gaps in these intervals, we have that there exists a rectangle $\widetilde{R} := \widetilde{R}(\widetilde{\mathbf{x}}') \subset B(\mathbf{x}, \frac{1}{5}D_v) \subset R$ with shortest and longest half-sides $\widetilde{l}$ and $\widetilde{L}$, respectively, satisfying $\widetilde{L}/\widetilde{l} \leq 2$ such that its center satisfies $\|\widetilde{\mathbf{x}}' - \mathbf{x}\| \leq \frac{1}{20}D_v \leq \widetilde{l}/2$ and such that a certain proper subset $\{B_{i_1}', \ldots, B_{i_u}'\}$ of the set of shrunken boxes lies entirely inside $\widetilde{R}$ and the others lie entirely outside $\widetilde{R}$. To see that the subset will be proper so that the number $u < v$, note that if the index $i_0$ yields the maximum in the definition of $D_v$, then $2^{\lambda_{i_0} - 3v} < \frac{1}{64}D_v$ so that

$$B(\mathbf{y}_{i_0}, 2^{\lambda_{i_0} - 3v}) \cap B(\mathbf{x}, \tfrac{1}{4}D_v) = \varnothing.$$

We now apply the inductive hypothesis with $u < v$. First we note that for all $i$, $B_i' \subset B(\mathbf{x}, 12\widetilde{l})$ since $12\widetilde{l} \geq \frac{12}{10}D_v$. We consider again two cases. If $12\widetilde{l} \leq l/8$, then by this construction, we have that

$$\begin{aligned} &P\Bigg(\mathcal{U}_4(\mathbf{x}; R) \cap \Bigg(\bigcap_{i=1}^{v} \mathcal{U}_4(\mathbf{y}_i; R)\Bigg)\Bigg) \\ &\quad \leq P\Bigg(\mathcal{U}_4(\mathbf{x}; \widetilde{R}) \cap \Bigg(\bigcap_{a=1}^{u} \mathcal{U}_4(\mathbf{y}_{i_a}; \widetilde{R})\Bigg)\Bigg) \qquad \text{(A.5)} \\ &\qquad \times P(\mathcal{U}_4(\mathbf{0}, 12\widetilde{l}; l/2)) \prod_{i \neq i_a} P(\mathcal{U}_4(\mathbf{0}, 2^{\lambda_i - 3v})). \end{aligned}$$

If instead $12\widetilde{l} > l/8$, we omit the factor $P(\mathcal{U}_4(\mathbf{0}, 12\widetilde{l}; l/2))$ in this inequality. By the induction hypothesis, we have that

$$P\left(\mathcal{U}_4(\mathbf{x}; \widetilde{R}) \cap \left(\bigcap_{a=1}^{u} \mathcal{U}_4(\mathbf{y}_{i_a}; \widetilde{R})\right)\right)$$
$$\leq C(u) P(\mathcal{T}_4(\mathbf{0}; \widetilde{l}/d(u))) \prod_{a=1}^{u} P(\mathcal{U}_4(\mathbf{0}, 2^{\lambda_{i_a} - c_1(u)})).$$

We take now the constants $d(u) := 96^u$ and $c_1(u) := 3u$. Thus, if $12\widetilde{l} > l/8$, we are done by this last inequality and (A.5). If, on the other hand, $12\widetilde{l} \leq l/8$, then we estimate by Lemmas 4 and 6 that

$$P(\mathcal{T}_4(\mathbf{0}; \widetilde{l}/d(u))) P(\mathcal{U}_4(\mathbf{0}, 12\widetilde{l}; l/2)) \leq C P(\mathcal{U}_4(\mathbf{0}, l/2)) \leq C' P(\mathcal{T}_4(\mathbf{0}, l/d(v))).$$

This completes the proof of the inductive step and therefore of the lemma. □

**Acknowledgment.** We thank the referee whose comments on the first version of the paper improved it considerably. In particular, the referee is credited with pointing out mistakes in our first attempted proof of Lemma 7. Moreover, the referee's careful reading and mastery in pointing out corrections was at once remarkable and invaluable.

Department of Mathematics
University of Colorado
Colorado Springs, Colorado 80933-7150
USA
e-mail: gjmorrow@math.uccs.edu
e-mail: yzhang@math.uccs.edu